\documentclass[11pt]{amsart}
\usepackage[all]{xy}
\usepackage{amssymb,tikz,color}
\usepackage[colorlinks=true,linkcolor=blue]{hyperref}
\usepackage[mathcal]{euscript}

\newcounter{noindnum}[subsection]
\renewcommand{\thenoindnum}{\roman{noindnum}}
\newcommand{\noindstep}{\refstepcounter{noindnum}{\rm(}\thenoindnum\/{\rm)} }
\newcommand{\stepzero}{\setcounter{noindnum}{0}
}

\renewcommand{\phi}{\varphi}
\renewcommand{\kappa}{\varkappa}
\renewcommand{\epsilon}{\varepsilon}

\newcommand{\divisor}{\mathcal D}

\newcommand{\HH}{\mathbb H}
\newcommand{\C}{\mathbb C}
\renewcommand{\P}{\mathbb P}
\newcommand{\A}{\mathbb A}
\newcommand{\Z}{\mathbb{Z}}
\newcommand{\R}{\mathbb{R}}
\newcommand{\Q}{\mathbb{Q}}
\newcommand{\bL}{\mathbb{L}}
\newcommand{\kk}{\mathbf{k}}
\newcommand{\dd}{\mathbf{d}}

\newcommand{\bE}{\mathbf E}
\newcommand{\bF}{\mathbf F}
\newcommand{\bG}{\mathbf G}

\newcommand{\Symm}{\mathfrak{S}}
\DeclareMathOperator{\ev}{ev}
\DeclareMathOperator{\gr}{gr}

\DeclareMathOperator{\FNF}{FNF}

\DeclareMathOperator{\sMot}{\overline{Mot}^{\lambda-\mathrm{ring}}}
\DeclareMathOperator{\cMot}{\overline{Mot}}
\DeclareMathOperator{\Mot}{Mot}
\DeclareMathOperator{\Spec}{Spec}

\DeclareMathOperator{\diag}{diag}
\DeclareMathOperator{\GL}{GL}
\DeclareMathOperator{\GSp}{GSp}
\DeclareMathOperator{\Lie}{Lie}

\DeclareMathOperator{\cl}{cl}

\DeclareMathOperator{\Id}{Id}

\DeclareMathOperator{\res}{res}

\DeclareMathOperator{\rk}{rk}
\DeclareMathOperator{\Hom}{Hom}
\DeclareMathOperator{\End}{End}

\DeclareMathOperator{\Pow}{\mathrm{Pow}}
\DeclareMathOperator{\Exp}{Exp}
\DeclareMathOperator{\Log}{Log}

\DeclareMathOperator{\tr}{tr}
\DeclareMathOperator{\Sym}{Sym}
\DeclareMathOperator{\Ker}{Ker}

\DeclareMathOperator{\Ext}{Ext}

\newcommand{\cC}{\mathcal C}

\newcommand{\cE}{\mathcal E}

\newcommand{\cH}{\mathcal H}

\newcommand{\cL}{\mathcal L}

\newcommand{\cO}{\mathcal O}
\newcommand{\cP}{\mathcal P}
\newcommand{\cS}{\mathcal S}

\newcommand{\cY}{\mathcal Y}
\newcommand{\cZ}{\mathcal Z}
\newcommand{\cX}{\mathcal X}
\newcommand{\Bun}{\mathcal{B}{un}}
\newcommand{\Nil}{\mathcal{N}{il}}
\newcommand{\Filt}{\mathcal{F}{ilt}}

\newcommand{\Pair}{\mathcal{P}{air}}

\newcommand{\Conn}{\mathcal{C}{onn}}

\newcommand{\END}{\mathcal{E}{nd}}
\newcommand{\HOM}{\mathcal{H}{om}}

\numberwithin{equation}{subsection}

\newtheorem{theorem}{Theorem}[subsection]
\newtheorem{proposition}[theorem]{Proposition}
\newtheorem{corollary}[theorem]{Corollary}
\newtheorem{lemma}[theorem]{Lemma}
\newtheorem{definition}[theorem]{Definition}
\newtheorem{problem}[theorem]{Problem}

{ \theoremstyle{definition}
\newtheorem{remark}[theorem]{Remark}
\newtheorem{remarks}[theorem]{Remarks}
\newtheorem{example}[theorem]{Example}
}

\title[Motivic classes of irregular Higgs bundles and connections]{Motivic classes of irregular Higgs bundles and irregular connections on a curve}

\author{Roman~Fedorov}
\address{Roman Fedorov, University of Pittsburgh, Pittsburgh, PA, USA \emph{and\/}
Max Planck Institute for Mathematics, Bonn, Germany}
\email{fedorov@pitt.edu}

\author{Alexander Soibelman}
\address{Alexander Soibelman, IHES, 35 route de Chartres, Bures-sur-Yvette, F-91440, France}
\email{asoibel@ihes.fr}

\author{Yan Soibelman}
\address{Yan Soibelman, Kansas State University, Manhattan, KS, USA}
\email{soibel@math.ksu.edu}

\begin{document}

\begin{abstract}
Let $X$ be a smooth projective curve over a field of characteristic zero and let~$\divisor$ be an effective divisor on $X$. We calculate motivic classes of various moduli stacks of parabolic vector bundles with irregular connections on~$X$ and of irregular parabolic Higgs bundles on $X$ with poles bounded by~$\divisor$ and with fully or partially fixed formal normal forms. Along the way, we obtain several results about irregular connections and irregular parabolic Higgs bundles. In particular, we give a criterion for the existence of a connection on a higher level parabolic bundle and also develop homological algebra for irregular connections and irregular parabolic Higgs bundles. We also simplify our previous results in the regular case by re-writing the formulas for motivic classes in terms of the HLV generating function.
\end{abstract}


\maketitle

\tableofcontents


\section{Introduction}

\subsection{Overview} Let $\kk$ be a field of characteristic zero and $X$ be a smooth geometrically connected projective curve over~$\kk$ (geometric connectedness means that $X$ remains connected after the base change to an algebraic closure of $\kk$). In~\cite{FedorovSoibelmansParabolic} we calculated motivic classes of moduli stacks of parabolic bundles with connections and parabolic Higgs bundles. In this paper, we extend these results to the irregular case. We also prove various results about irregular connections and irregular Higgs bundles that are of independent interest. We emphasize that while the regular case has been studied extensively, very little is known about motivic classes in the irregular case at the moment.

Let us illustrate in a toy-model example why the irregular case is more complicated. A toy-model for the parabolic Higgs bundles is finite-dimensional vector spaces with endomorphisms preserving flags (these are ``the parabolic Higgs bundles over a zero-dimensional variety''). In the irregular case, the toy-model for Higgs bundles is free modules over the ring $\kk[z]/z^n$ endowed with nilpotent endomorphisms preserving flags of free submodules. These objects are considerably more complicated, in particular, because the homological dimension of the category of modules over the ring $\kk[z]/z^n$ is infinite.

\subsection{Irregular $\epsilon$-connections}\label{sect:IntroEpsConn}
Let $D\subset X(\kk)$ be a finite set of $\kk$-rational points of $X$ (possibly empty). Let $\divisor=\sum_{x\in D}n_xx$ be a divisor on $X$, where $n_x$ are positive integer numbers. Thus $\divisor$ is an~effective divisor with set-theoretic support $D$. Let $\epsilon$ be an element of $\kk$ and $E$ be a vector bundle on $X$. Recall that an \emph{$\epsilon$-connection\/} on $E$ with \emph{poles bounded by $\divisor$\/} is a morphism of sheaves of abelian groups $\nabla:E\to E\otimes\Omega_X(\divisor)$ satisfying the Leibniz rule:
\[
    \nabla(fs)=f\nabla(s)+\epsilon s\otimes df.
\]
Here $f$ and $s$ are any local sections of $\cO_X$ and $E$ respectively. If $\epsilon\ne0$ and $\nabla$ is an $\epsilon$-connection, then $\epsilon^{-1}\nabla$ is a usual singular connection with singularities bounded by $\divisor$. If $\epsilon=0$, then $\epsilon$-connection is the same as a singular Higgs field with singularities bounded by $\divisor$. Thus, the notion of $\epsilon$-connection allows us to talk about connections and Higgs bundles in a uniform way. Trivializing the vector bundle~$E$ locally near $x\in D$, we can write
\begin{equation}\label{eq:EpsConn}
    \nabla=\epsilon d+A,
\end{equation}
where $A$ is a $\rk E\times\rk E$ matrix whose entries are 1-forms with poles of order up to $n_x$. We say that $x\in D$ is a \emph{regular singular point\/} if $n_x=1$ and an \emph{irregular singular point\/} if $n_x\ge2$.

\subsection{Motivic invariants}\label{sect:MotivicInvariants}
We start by recalling the categorical approach to motivic Donaldson--Thomas invariants, which was proposed for the first time in~\cite{KontsevichSoibelman08}. It has been used directly or indirectly in many papers devoted to motivic Donaldson--Thomas invariants, including this one.

Recall that in~\cite{FedorovSoibelmans} we defined a version of the ring of motivic classes of $\kk$-stacks denoted $\Mot(\kk)$. If $\cX$ is an Artin stack of finite type over $\kk$, we have its motivic class $[\cX]\in\Mot(\kk)$. For a short account of the theory of motivic classes we refer the reader to Section~\ref{sect:MotClasses} below (see also~\cite[Sect.~4.1--4.2]{KontsevichSoibelman08}). One should think of $[\cX]$ as the universal enumerative invariant of $\cX$ with values in a commutative ring. Morally, $\Mot(\kk)$ is the ``Grothendieck group of an appropriate category of motives'' though the theory of motivic classes can be developed independently of the theory of motives.

Let $\Conn(\epsilon,X,\divisor,\zeta)$ be the category of $\epsilon$-connections on $X$ with poles bounded by $\divisor$, compatible parabolic structures of level $\divisor$, and normal forms $\zeta$ at singular points (see Section~\ref{sect:ConnFnf} for precise definitions). The objects of this category form an Artin stack locally of finite type, for which we use the same notation. As explained in Section~\ref{sect:GammaD}, this stack decomposes according to the classes of parabolic bundles:
\[
    \Conn(\epsilon,X,\divisor,\zeta)=\bigsqcup_{\gamma\in\Gamma_D,d\in\Z}\Conn_{\gamma,d}(\epsilon,X,\divisor,\zeta),
\]
where $\Gamma_D$ is a commutative monoid defined in Section~\ref{sect:GammaD} below. We can form the ``motivic series''
\[
    \sum_{\gamma\in\Gamma_D,d\in\Z}[\Conn_{\gamma,d}(\epsilon,X,\divisor,\zeta)]e_{\gamma,d}\in
    \cMot[\kk][[\Gamma_D,\Z]].
\]
The abelian group $\cMot[\kk][[\Gamma_D,\Z]]$ is closely related to \emph{quantum tori\/}, see~\cite{KontsevichSoibelman08} and~\cite[Sect.~1.6]{FedorovSoibelmans}. We note that the quantum torus is commutative in our case because $\Conn(\epsilon,X,\divisor,\zeta)$ is a 2-dimensional Calabi--Yau category. If $\epsilon\ne0$, so that we are working with connections, the stacks are of finite type (see Corollary~\ref{cor:StabConn}) and the series above makes sense. If $\epsilon=0$, so that we are working with Higgs fields, we need to add stability condition $\sigma$ and consider the open substacks $\Conn^{\sigma-ss}_{\gamma,d}(0,X,\divisor,\zeta)$.

The starting point of our calculation is considering the category of parabolic bundles $\bE$ with twisted nilpotent endomorphisms $\bE\to\bE\otimes\cL$, where $\cL$ is a fixed line bundle. We denote the category as well as its stack of objects by $\Pair^{nilp}(X,D,\cL)$. The ``components'' $\Pair^{nilp}_{\gamma,d}(X,D,\cL)$ of these stacks have infinite motivic classes, so one needs to ``regularize'' them. The technique for this was developed in~\cite{KontsevichSoibelman08} and consists of looking at subcategory consisting of objects whose Harder--Narasimhan spectrum belongs to a fixed strict sector. Following Schiffmann (see~\cite{SchiffmannIndecomposable}), we use a~somewhat different framework in this paper: a stability condition is imposed on the vector bundle rather than on the pair consisting of a parabolic bundle with an endomorphism (or an $\epsilon$-connection). The natural forgetful functors from $\Conn(\epsilon,X,\divisor,\zeta)$ and $\Pair^{nilp}(X,D,\cL)$ to the category of vector bundles on $X$ allows us to utilize the methods of~\cite{KontsevichSoibelman08}, see~\cite[Sect.~1.6]{FedorovSoibelmans}.

\subsection{Motivic classes of stacks of semistable Higgs bundles on a curve: a short history of the subject}
There are essentially two approaches to calculating the motivic classes of stacks of semistable Higgs bundles. The first approach lowers the homological dimension of the category, while the second one raises it. In the first approach one uses the forgetful map from the stack of Higgs bundles to the stack of vector bundles, and reduces the computations to the homologically 1-dimensional category of coherent sheaves on the curve $X$. This approach goes back to~\cite{GarciaPradaHeinlothSchmitt}. The first general result was obtained by Schiffmann (see~\cite{SchiffmannIndecomposable}) and Mozgovoy and Schiffmann (see~\cite{MozgovoySchiffmann2020}), where it was used for the computations of the number of points over a finite field rather than of motivic classes. For the fields of characteristic zero, the motivic classes were computed in our paper~\cite{FedorovSoibelmans}. This computation follows the strategy of Mozgovoy and Schiffmann but is more involved. In~\cite{FedorovSoibelmans} we also computed the motivic classes of vector bundles with connections.

In the second approach (see~\cite{DiaconescuDonagiPantev}) one uses the description of the category of Higgs bundles in terms of the category of pure 1-dimensional coherent sheaves on $T^\ast X$. The idea is to relate the motivic invariants of semistable objects of this category to some invariants (specifically, Pandharipande--Thomas invariants) of the category of coherent sheaves on the 3-dimensional Calabi--Yau variety $T^\ast X\times {\mathbb A}^1$. In other words, one raises the homological dimension of the category from~2 to 3, with the goal of utilizing the approach of~\cite{KontsevichSoibelman08}.

One of the motivations for the upgrade of the above story to \emph{parabolic\/} Higgs bundles were the conjectures of Hausel, Letellier, and Rodriguez-Villegas (a.k.a.~HLV conjectures), see~\cite{HauselLetellierRodriguezVillegas}. The HLV conjectures deal with the cohomology of character varieties with marked points, but their relation to Donaldson-Thomas theory in the spirit of~\cite{KontsevichSoibelman08} was apparent from the very beginning.

We would like to mention that Mozgovoy formulated in~\cite{mozgovoy2012ADHM} an important conjecture about motivic classes of the twisted Higgs bundles (that is, the canonical bundle $\Omega_X$ is replaced by an arbitrary line bundle $\cL$ such that $\deg\cL\ge\deg\Omega_X$), which was probably the first attempt to introduce rigorously the motivic techniques to HLV conjectures.

In the paper by Chuang, Diaconescu, Donagi, and Pantev~\cite{ChuangDiaconescuDonagiPantev}, a confirmation of the HLV conjectures was obtained with the help of topological string theory (so it is not mathematically rigorous). The authors used the raising of the homological dimension from~2 to~3 and the relation to the Pandharipande--Thomas invariants of the corresponding Calabi--Yau 3-fold.\footnote{Maulik seems to have unpublished results which justify the ``physics part'' of the computations of Chuang--Diaconescu--Donagi--Pantev but without the comparison with Mozgovoy--Schiffmann type formulas.} The answers given in~\cite{SchiffmannIndecomposable,MozgovoySchiffmann2020} are more complicated than the answers proposed in~\cite{HauselLetellierRodriguezVillegas}, which were given in terms of the so-called Hausel--Letellier--Villegas generating function (HLV generating function for short). It was Mellit who finally developed a mathematically rigorous way to re-write the answers of Mozgovoy and Schiffmann in terms of the HLV generating function (see~\cite{MellitIntegrality,MellitNoPunctures,MellitPunctures}). We will revisit his technique in Section~\ref{sect:MellitSimplification}.

In the paper~\cite{FedorovSoibelmansParabolic}, using the technique of~\cite{MellitPunctures}, we calculated the motivic classes of semistable parabolic Higgs bundles as well as of parabolic Higgs bundles with connections over arbitrary fields of characteristic zero. The answers are given in terms similar to~\cite{SchiffmannIndecomposable} rather than in terms of the HLV generating function. In this paper, we will give the answers using an irregular generalization of the HLV generating function. Note that in the case of \emph{irregular\/} parabolic Higgs bundles (or parabolic bundles with irregular connections) nothing was known about motivic invariants except for the conjectures proposed in~\cite{HauselMerebWong} and the mathematically non-rigorous results in~\cite{DiaconescuDonagiPantev}.

\subsection{Parabolic structures of level $\divisor$}\label{sect:ParStr} In the next few subsections we give a detailed formulation of one of our main results. As in~\cite{FedorovSoibelmansParabolic}, we will be interested in $\epsilon$-connections with compatible parabolic structures. One reason for that is that the corresponding moduli spaces are symplectic, see Section~\ref{sect:Twisted}. Recall that $X$ is a smooth geometrically connected projective curve over a field $\kk$ of characteristic zero and $\divisor$ is an effective divisor with set-theoretic support $D\subset X(\kk)$. For a rational point $x\in X(\kk)$, we let $nx$ be the $n$-th infinitesimal neighborhood of $x$, so that $k[nx]\approx\kk[z]/z^n$. The reader might want to look at Section~\ref{sect:Infinitesimal} for more details on infinitesimal neighborhoods.
\begin{definition}\label{def:thick}
Let $E$ be a vector bundle on $X$ and $x\in X(\kk)$. For $n\in\Z_{>0}$, \emph{a level $n$ parabolic structure at $x$ on $E$\/} is a filtration
\begin{equation}\label{eq:ThickPar}
    E|_{nx}=E_0\supset E_1\supset\ldots\supset E_N=0\text{ for }N\gg0
\end{equation}
of the free $k[nx]$-module $E|_{nx}$ by \emph{free\/} submodules. For a divisor $\divisor=\sum_{x\in D}n_xx$ with $n_x>0$, \emph{a level~$\divisor$ parabolic structure on $E$\/} is a collection $E_{\bullet,\bullet}=(E_{x,j})$ where $E_{x,j}$ is a level $n_x$ parabolic structure on~$E$ at $x$ for all $x\in D$. The pair $\bE=(E,E_{\bullet,\bullet})$ is called a \emph{level $\divisor$ parabolic bundle on $X$}.
\end{definition}
We have the stack $\Bun^{par}(X,\divisor)$ classifying level $\divisor$ parabolic bundles on $X$. The forgetful morphism $\Bun^{par}(X,\divisor)\to\Bun(X)$, where $\Bun(X)$ is the stack of vector bundles on $X$, is schematic and locally of finite type.

Given a level $n$ parabolic structure at $x$ on $E$, one can always trivialize $E$ near $x$ so that each $E_i$ is the submodule of $E|_{nx}$ generated by $e_1,\ldots,e_{\rk E_i}$, where $e_1,\ldots,e_{\rk E}$ is the basis provided by the trivialization. We say that such a trivialization is \emph{compatible with the parabolic structure}.

\subsection{Classes of parabolic bundles}\label{sect:GammaD} The stack $\Bun^{par}(X,\divisor)$ decomposes according to the classes of parabolic bundles. In more detail, for any finite set $D$, let $\Gamma_D$ denote the commutative monoid of sequences $\gamma=(r,r_{\bullet,\bullet})$, where $r$ is a nonnegative integer and $r_{\bullet,\bullet}$ is a sequence of nonnegative integers indexed by $D\times\Z_{>0}$ subject to the following condition: \emph{for all $x\in D$ we have $\sum\limits_{j=1}^{\infty}r_{x,j}=r$}. In particular, for all $x\in D$ we have $r_{x,j}=0$ for $j$ large enough. The operation on $\Gamma_D$ is the componentwise addition. For $\gamma=(r,r_{\bullet,\bullet})\in\Gamma_D$ we set $\rk\gamma=r$. Note also that if $D=\emptyset$, then $\Gamma_D=\Z_{\ge0}$.

If $\bE:=(E,E_{\bullet,\bullet})$ is a level $\divisor$ parabolic bundle on $X$, we define its \emph{class\/} by
\[
    \cl(\bE):=(\rk E,(\rk E_{x,j-1}-\rk E_{x,j})_{x\in D,j>0},\deg E)\in\Gamma_D\times\Z.
\]
We have
\begin{equation}\label{Bun_gamma_d}
    \Bun^{par}(X,\divisor)=\bigsqcup_{\substack{\gamma\in\Gamma_D\\d\in\Z}}\Bun_{\gamma,d}^{par}(X,\divisor).
\end{equation}

\begin{definition}\label{def:full}
    A class $\gamma=(r_{\bullet,\bullet})\in\Gamma_D$ is called \emph{full\/} at $x\in D$ if for all $j>0$ we have $r_{x,j}\in\{0,1\}$.
\end{definition}
Let $\cl(E,E_{\bullet,\bullet})=(\gamma,d)$. Then $\gamma$ is full at $x\in D$ if and only if $E_{x,\bullet}$ is a full flag of submodules of $E_{n_xx}$ (with repetitions, that is, for some $j$ we have $E_{x,j-1}=E_{x,j}$).

\subsection{Parabolic $\epsilon$-connections with fixed formal normal forms}\label{sect:ConnFnf} Let $X$, $D$, and $\divisor$ be as before. Define the space of formal normal forms (we will explain the terminology in Remark~\ref{rem:FNF})
\begin{equation}\label{eq:FNF_full}
    \FNF(\divisor):=\prod_{x\in D}(\Omega_X(n_xx)/\Omega_X)^{\Z_{>0}},
\end{equation}
where $\Omega_X(n_xx)/\Omega_X$ is a vector space of polar parts of 1-forms at $x$. Recall that $\epsilon$-connections can be presented in the form~\eqref{eq:EpsConn}. Recall also from Section~\ref{sect:ParStr} the notion of local trivialization compatible with a parabolic structure. We are ready to define our main object of study.
\begin{definition}\label{def:ModSpace}
  Let $X$ be a smooth geometrically connected projective curve over the field~$\kk$. Let $D\subset X(\kk)$ be a finite set of $\kk$-rational points of $X$. Let $\divisor=\sum_{x\in D}n_xx$ be a divisor on $X$, where~$n_x$ are positive integer numbers. Let $\epsilon\in\kk$ and $\zeta=(\zeta_{x,j}\,|\;x\in D,j\in\Z_{>0})\in\FNF(\divisor)$. Then $\Conn(\epsilon,X,\divisor,\zeta)$ is the moduli stack classifying triples $(E,E_{\bullet,\bullet},\nabla)$, where\\
  $\bullet$ $(E,\nabla)$ is an $\epsilon$-connection on $X$ with poles bounded by $\divisor$;\\
  $\bullet$ $E_{\bullet,\bullet}$ is a level $\divisor$ parabolic structure on $E$;\\
  $\bullet$ for $x\in D$ and a local trivialization of $E$ near $x$ compatible with the parabolic structure write
    \[
        \nabla=\epsilon d+A
    \]
    as in~\eqref{eq:EpsConn} and view $A$ as a block matrix with blocks of sizes $\rk(E_{x,i-1}/E_{x,i})\times\rk(E_{x,j-1}/E_{x,j})$. Then the polar part of the $i$-th diagonal block is equal to $\zeta_{x,i}$, and the blocks lying below the diagonal have no pole at $x$.
\end{definition}
It is easy to see that the above compatibility condition is satisfied for every local trivialization of~$E$ near $x$, compatible with the parabolic structure, if it satisfied for a single such trivialization. We also note that for some $i$ we have $\rk(E_{x,i-1}/E_{x,i})=0$, so the corresponding blocks are ``trivial''.

The stack $\Conn(\epsilon,X,\divisor,\zeta)$ decomposes according to the classes $(\gamma,d)\in\Gamma_D\times\Z$ of the parabolic bundle (cf.~\eqref{Bun_gamma_d}). We denote the corresponding closed and open substack by $\Conn_{\gamma,d}(\epsilon,X,\divisor,\zeta)$. We note that when $\epsilon=0$ this definition is similar to~\cite[Sect.~2.1]{DiaconescuDonagiPantev} except that we do not put any stability conditions at the moment. In Section~\ref{sect:Twisted}, we will show that the stacks $\Conn_{\gamma,d}(\epsilon,X,\divisor,\zeta)$ are twisted cotangent bundles of the stacks $\Bun^{par}(X,\divisor)$. Here is the main problem we will be interested in:
\begin{problem}
  Calculate the motivic class of $\Conn_{\gamma,d}(\epsilon,X,\divisor,\zeta)$.
\end{problem}
More precisely, if $\epsilon=0$, the stack is of infinite type and has infinite motivic volume, so we need to add stability conditions.

Note that a local section of $\Omega_X(n_x)$ can have a pole of order strictly less than $n_x$ at $x$.
\begin{definition}\label{def:non-resonant0}
  Let $\gamma=(r,r_{\bullet,\bullet})\in\Gamma_D$ and $\zeta=(\zeta_{x,j})\in\FNF(\divisor)$. We will say that $\zeta$ is
  \emph{non-resonant for $\gamma$ at $x\in D$} if for all $1\le i<j$ such that $r_{x,i}\ne0\ne r_{x,j}$ the polar part $\zeta_{x,i}-\zeta_{x,j}$ has a pole of order exactly $n_x$.
\end{definition}

We will only be able to calculate motivic classes of stacks of parabolic $\epsilon$-connections in the case when its class $\gamma$ is full at all irregular singular points (that is, at all $x\in D$ with $n_x\ge2$) and $\zeta$ is non-resonant for $\gamma$ at each irregular point (see Theorem~\ref{th:IntroMain} below). In this case, the parabolic structures at irregular points are determined by the $\epsilon$-connection uniquely (see Remark~\ref{rem:NotNeeded}).\footnote{Remark~\ref{rem:NotNeeded} is about connections with partially fixed normal forms, see Section~\ref{sect:IntroIrregFNF}. The required result about connections with fully fixed formal normal forms is obtained by taking $\divisor'=\divisor$.} The parabolic structures at regular singular points are still extra data but we note that these structures are ``thin'', that is, of level~1.

\begin{remark}\label{rem:FNF}
  Let $(\bE,\nabla)$ be a point of $\Conn_{\gamma,d}(\epsilon,X,\divisor,\zeta)$ and pick $x\in D$. One can show that if $\zeta$ is non-resonant for $\gamma$ at $x$, then there are lifts $\tilde\zeta_{x,j}$ of the polar parts $\zeta_{x,j}$ to local sections of $\Omega(n_xx)$ such that $\nabla=\epsilon d+\diag(\tilde\zeta_1,\ldots,\tilde\zeta_{\rk E})$. Indeed, one can use an argument similar to the proof of Lemma~\ref{lm:diag} if $\gamma$ is full and similar to~\cite[Appendix~A.1]{DiaconescuDonagiPantev} in the general case.

  Moreover, if $\epsilon\ne0$ and $(\bE',\nabla')$ is a point of $\Conn_{\gamma,d'}(\epsilon,X,\divisor,\zeta')$, then the $\epsilon$-connections are gauge equivalent in the formal neighborhood of $x$ if and only if $\zeta_{x,j}=\zeta'_{x,j}$ whenever $r_{x,j}\ne0$ (here $\gamma=(r,r_{\bullet,\bullet}$)). This motivates the terminology ``formal normal forms''; cf.~\eqref{eq:FNF_full}.
\end{remark}

\subsection{$\lambda$-rings of motivic classes}\label{sect:MotClasses} We are going to define the ring of motivic classes where our final answers will be given. Recall that in~\cite[Sect.~1.1]{FedorovSoibelmans} we defined (following earlier works~\cite[Sect.~1]{Ekedahl09}, \cite{Joyce07}, and~\cite{KontsevichSoibelman08}) the ring of motivic classes of Artin $\kk$-stacks denoted $\Mot(\kk)$. As an abelian group, $\Mot(\kk)$ is generated by symbols $[\cS]$, where~$\cS$ is a stack of finite type over $\kk$ with affine stabilizers, modulo the relations $[\cS]=[\cS-\cS']+[\cS']$, whenever $\cS'$ is a closed substack of $\cS$, and $[\cS]=[\cS'\times\A^r_\kk]$, whenever $\cS\to\cS'$ is a rank $r$ vector bundle. We also defined the dimensional completion $\cMot(\kk)$. In this paper, we prefer to work with a quotient of $\cMot(\kk)$. Recall that a \emph{pre-$\lambda$-ring\/} is a commutative ring~$R$ endowed with a group homomorphism $\lambda\colon R\to(1+zR[[z]])^\times$ such that $\lambda\equiv1+\Id_R\pmod{z^2R[[z]]}$. We write $\lambda=\sum_{i=0}^\infty\lambda_iz^i$, so that the previous conditions imply $\lambda_0=1$, $\lambda_1=\Id_R$. A \emph{$\lambda$-ring\/} is a~pre-$\lambda$-ring $R$ such that $\lambda(1)=1+z$ and for all $x,y\in R$ and $m,n\in\Z_{\ge0}$, $\lambda_n(xy)$ and $\lambda_n(\lambda_m(x))$ can be expressed in terms of $\lambda_i(x)$ and $\lambda_j(x)$ using the ring operations in a standard way (see~\cite[Exp.~5, Def.~2.4, (2.4.1)]{SGA6} and~\cite[Ch.~I, Sect.~1]{knutson1973rings}). In~\cite[Sect.~1.3.1]{FedorovSoibelmans} we defined a structure of a pre-$\lambda$-ring on $\cMot(\kk)$. The characteristic feature of this structure is that for a \emph{quasi-projective $\kk$-variety\/} $Y$ the $\lambda$-operation coincides with Kapranov's motivic $\zeta$-function, that is, we have
\begin{equation}\label{eq:zeta}
    \lambda([Y])=\zeta_Y(z):=\sum_{i=0}^\infty[Y^i/\Symm_i]z^i,
\end{equation}
where $Y^i/\Symm_i$ is the $i$-th symmetric power of $Y$ ($\Symm_i$ is the permutation group). In this paper, we will use the \emph{opposite pre-$\lambda$-ring structure} given by $\lambda([Y]):=(\zeta_Y(-z))^{-1}$. We, however, expect that this is not a $\lambda$-ring structure. We denote by $\sMot(\kk)$ the universal quotient of $\cMot(\kk)$ that is a $\lambda$-ring, see~\cite{LarsenLuntsRational}. We note that in the pre-$\lambda$-ring structure used in~\cite{FedorovSoibelmans} we had $\lambda(1)=(1-z)^{-1}$ so we could not have a non-trivial homomorphism to a $\lambda$-ring. This is why we are working with the opposite pre-$\lambda$-ring structure here (see also Section~\ref{sect:MotHomo}). We need the $\lambda$-ring structure because we would like to make use of Mellit's results~\cite{MellitIntegrality,MellitNoPunctures,MellitPunctures}, which were proved in some universal $\lambda$-ring. We will construct a homomorphism from that $\lambda$-ring to $\sMot(\kk)$ in Proposition~\ref{pr:HomFromUniv} (see also Remark~\ref{rem:NotHomo}). In fact, we will first calculate our motivic classes in $\cMot(\kk)$ (see Theorems~\ref{th:AnswerConnInMot} and~\ref{th:AnswerInCmot}) but the formulas we obtain in $\sMot(\kk)$ are much simpler. We note that instead of $\sMot(\kk)$ we could work in the $\lambda$-ring associated with Chow motives as in, e.g.,~\cite{mozgovoy2006stringy}. For the proof that it is a $\lambda$-ring see, e.g.,~\cite[Lemma~4.1]{FHeinlothFunctEq} or~\cite{Getzler1995mixed}.

\subsection{Partitions, nilpotent linear operators, and motivic modified Macdonald polynomials}\label{sect:Macdonald}
We denote by $\cP$ the set of partitions. For a partition $\mu=\mu_1\ge\mu_2\ldots$, we denote by $\mu'=\mu'_1\ge\mu'_2\ge\ldots$ the conjugate partition. We set $|\mu|:=\sum_i\mu_i$, $n(\mu):=\sum_i(i-1)\mu_i$. For two partitions $\mu$ and $\nu$ we set $\langle\mu,\nu\rangle:=\sum\mu'_i\nu'_i$ so that $\langle\mu,\mu\rangle=2n(\mu)+|\mu|$.

The conjugacy classes of nilpotent linear operators on $r$-dimensional vector spaces are classified by partitions of $r$. More precisely, we associate to a nilpotent linear operator the partition $\mu=\mu_1\ge\mu_2\ldots$, where $\mu_i$ are sizes of the Jordan blocks. This convention agrees with~\cite{SchiffmannIndecomposable,MozgovoySchiffmann2020,OGormanMozgovoy} but differs from that of~\cite{MellitPunctures} by conjugation.

Let $w_\bullet=(w_1,w_2,\dots)$ be an infinite sequence of variables and $R$ be a commutative ring\footnote{By our convention, rings are always unital.}. We denote by $\Sym_R(w_\bullet)$ the ring of symmetric functions with coefficients in $R$. Recall the modified Macdonald polynomials $\tilde H_\mu(w_\bullet;z,q)\in\Sym_{\Z[z,q]}(w_\bullet)$ indexed by partitions $\mu$; see, e.g.,~\cite{HaglundEtAlOnMacdonaldPoly} and references therein. Note that, formally speaking, symmetric functions are not polynomials (they become polynomials upon plugging in $w_{N+1}=w_{N+2}=\dots=0$ for some $N\ge0$).

Let $\bL=\big[\A_\kk^1\big]$ be the motivic class of the affine line (a.k.a.~Lefshetz motive). There is a unique homomorphism $\Sym_{\Z[z,q]}(w_\bullet)\to\Sym_{\cMot(\kk)[z]}(w_\bullet)$ sending $w_\bullet$ to $w_\bullet$, $z$ to $z$, and $q$ to $\bL$. We denote by $\tilde H_\mu^{mot}(w_\bullet;z)$ the image of $\tilde H_\mu(w_\bullet;z,q)$ under this homomorphism. We keep the same notation for the image of $\tilde H_\mu^{mot}(w_\bullet;z)$ in $\Sym_{\sMot(\kk)[z]}(w_\bullet)$.

\begin{remark}
  In~\cite{FedorovSoibelmansParabolic} we used Mellit's parameterization of conjugacy classes of nilpotent operators in the modified motivic Macdonald polynomials, but we used Schiffmann's parameterization in ``Schiffmann's terms'' $J_\lambda^{mot}$, $H_\lambda^{mot}$. That resulted in an inconsistency in the resulting formula. We correct this now by defining the modified motivic Macdonald polynomials slightly differently: namely, we substitute $\bL$ into $H_\mu(w_\bullet;z,q)=H_{\mu'}(w_\bullet;q,z)$ rather than into $H_\mu(w_\bullet;q,z)$. This corresponds to switching from Mellit's parameterization to Schiffmann's.
\end{remark}

\subsection{The motivic generating function}\label{sect:IntroHLV}
We are going to define our main generating function. Recall that $D\subset X(\kk)$ is a set of rational points of $X$. Recall from Section~\ref{sect:GammaD} the monoid $\Gamma_D$. We define the completed ring $\sMot(\kk)[[\Gamma_D]]$ in the natural way, and write its elements as $\sum\limits_{\gamma\in\Gamma_D}A_\gamma e_\gamma$, where $A_\gamma\in\sMot(\kk)$, $e_\gamma$ are basis vectors. It is convenient to identify $e_\gamma$ with a monomial
\begin{equation}\label{eq:MotivicSeries}
    w^\gamma=w^r\prod_{x\in D}\prod_{j=1}^{\infty}w_{x,j}^{r_{x,j}},
\end{equation}
where $\gamma=(r,r_{\bullet,\bullet})$ and $w_{\bullet,\bullet}=(w_{x,j})$ is a sequence of variables indexed by $D\times\Z_{>0}$. We may identify $\sMot(\kk)[[\Gamma_D]]$ with a subring of $\sMot(\kk)[[w_{\bullet,\bullet},w]]$. We will write $\sMot(\kk)[[\Gamma_D,z]]$ for $\sMot(\kk)[[\Gamma_D]][[z]]$ so that the elements of this ring can be written as $\sum\limits_{\gamma\in\Gamma_D,d\in\Z_{\ge0}}A_\gamma w^\gamma z^d$.

\begin{remarks}
(i) We note that these rings are closely related to completed quantum tori considered in~\cite{KontsevichSoibelman08,KontsevichSoibelman10} and mentioned earlier. In our case, they are commutative because we are working with 2-dimensional Calabi--Yau categories (see~\cite[Sect.~1.5]{FedorovSoibelmansParabolic} for the discussion and details).

(ii) When $D\ne\emptyset$ and $(r,r_{\bullet,\bullet})\in\Gamma_D$, the integer $r$ is determined by $r_{\bullet,\bullet}$. Thus, we can equivalently define $\Gamma_D$ as a monoid of sequences $r_{\bullet,\bullet}$ such that for all $x,y\in D$ we have $\sum_jr_{x,\bullet}=\sum_jr_{y,\bullet}$. Also we can plug $w=1$ in~\eqref{eq:MotivicSeries} and in other motivic series without loosing information. It is, however, convenient to keep the variable $w$.
\end{remarks}

Recall that $X$ is a smooth projective geometrically connected curve. We have the motivic $\zeta$-function $\zeta_X(z)$ defined by~\eqref{eq:zeta}. In this paper it will be more convenient for us to work with the motivic L-function, which is actually a polynomial of degree $2g$, where $g:=g(X)$ is the genus of $X$, see~\cite[Prop.~1.3.1]{FedorovSoibelmans}:
\[
    L_X(z):=\zeta_X(z)(1-z)(1-\bL z)\in\sMot[z].
\]
Let $\delta$ be a nonnegative integer (later, we will take $\delta:=\sum_{x\in D}(n_x-1)$ to be the irregularity of $\divisor$). Set
\begin{multline}\label{eq:HLV}
    \Omega^{mot}_{X,D,\delta}\\
    :=\sum_{\mu\in\cP}w^{|\mu|}(-\bL^\frac12)^{(2g+\delta)\langle\mu,\mu\rangle}
    z^{2\delta n(\mu')}\prod_{\Box\in\mu}\frac{L_X(z^{2a(\Box)+1}\bL^{-l(\Box)-1})}
    {(z^{2a(\Box)+2}-\bL^{l(\Box)})(z^{2a(\Box)}-\bL^{l(\Box)+1})}
    \prod_{x\in D}\tilde H^{mot}_\mu(w_{x,\bullet};z^2).
\end{multline}
Here the sum is over all partitions and the product is over all boxes of the partition, $a(\Box)$ and $l(\Box)$ are the arm and the leg of the box respectively. We expand the denominators in positive powers of $z$, so that we have
\[
    \Omega^{mot}_{X,D,\delta}\in\sMot(\kk)[\bL^{\frac12}][[\Gamma_D,z]].
\]
In fact, $\Omega^{mot}_{X,D,\delta}$ is symmetric in each of the sequences of variables $w_{x,\bullet}$. This generating function is closely related to the generating functions considered in~\cite{HauselMerebWong} and~\cite{OGormanMozgovoy}.

\subsection{Plethystic exponent and logarithm}\label{sect:Plethystic}
We need the notions of \emph{plethystic exponent and logarithm}. Let $\sMot(\kk)[[\Gamma_D,z]]^0$ denote the subset of $\sMot(\kk)[[\Gamma_D,z]]$ consisting of series with zero constant terms. Then we have a~bijection
\begin{equation}\label{eq:ExpLog}
\begin{split}
 \Exp\colon&\sMot(\kk)[\bL^\frac12][[\Gamma_D,z]]^0\to1+\sMot(\kk)[\bL^\frac12][[\Gamma_D,z]]^0\colon\\
 &\sum_{\gamma,d} A_\gamma\bL^{n_\gamma}w^\gamma z^d\mapsto\prod_{\gamma,d}\lambda^{-1}_{A_\gamma}(-\bL^{n_\gamma}w^\gamma z^d),
\end{split}
\end{equation}
where we use the $\lambda$-ring structure on $\sMot(\kk)$ (see also Sections~\ref{sect:FullyFixed} and~\ref{sect:ExpLog} below). We note that in~\cite[Sect.~1.3]{FedorovSoibelmansParabolic} we used the opposite pre-$\lambda$-ring structure, but the above plethystic exponent agrees with the one in loc.~cit. Let the plethystic logarithm $\Log$ be the inverse bijection. Similarly, we have the inverse bijections
\[
 \Exp\colon\sMot(\kk)[\bL^\frac12][[\Gamma_D]]^0\leftrightarrows1+\sMot(\kk)[\bL^\frac12][[\Gamma_D]]^0\colon\Log.
\]
Note that $\Omega^{mot}_{X,D,\delta}\in1+\sMot(\kk)[\bL^\frac12][[\Gamma_D,z]]^0$. Set
\begin{equation}\label{eq:HH_mot}
    \HH^{mot}_{X,D,\delta}:=(1-z^2)\Log\Omega^{mot}_{X,D,\delta}.
\end{equation}
We will see in Corollary~\ref{cor:MellitSimpl} that the coefficients of $\HH^{mot}_{X,D,\delta}$ are polynomials in $z$.

\subsection{Explicit formulas for motivic classes of $\Conn_{\gamma,d}(\epsilon,X,\divisor,\zeta)$}\label{sect:ExplAnswers}
We are almost ready for our first main result, but we need some additional notation. Recall that $X$ is a smooth projective geometrically connected curve over $\kk$, and $\divisor$ is an effective divisor on $X$ with the set-theoretic support $D\subset X(\kk)$. Recall from~\eqref{eq:FNF_full} the vector space $\FNF(\divisor)$ of formal normal forms. For $\zeta\in\FNF(\divisor)$, $\epsilon\in\kk$, $\gamma\in\Gamma_D$, and $d\in\Z$, we have Artin stacks $\Conn_{\gamma,d}(\epsilon,X,\divisor,\zeta)$, see Definition~\ref{def:ModSpace} and what follows.

For $\gamma=(r,r_{\bullet,\bullet})\in\Gamma_D$ and $\zeta=(\zeta_{x,j})\in\FNF(\divisor)$ set
\begin{equation}\label{eq:Intro_star}
\gamma\star\zeta:=\sum_{x\in D}\sum_{j>0}r_{x,j}\res\zeta_{x,j}\in\kk.
\end{equation}
The significance of this is that $\Conn_{\gamma,d}(\epsilon,X,\divisor,\zeta)$ is empty unless $\epsilon d+\gamma\star\zeta=0$ (see Corollary~\ref{cor:indecomposable summand}).

Set $\delta=\delta_\divisor:=\sum_{x\in D}(n_x-1)=\deg\divisor-|D|\in\Z_{\ge0}$, where $|D|$ is the number of points in $D$. For $\gamma=(r,r_{\bullet,\bullet})\in\Gamma_D$, we define:
\begin{equation}\label{eq:BetterChi2}
    \chi(\gamma)=\chi(\gamma,g,\divisor):=(2g-2)r^2-\delta r+2\sum_{x\in D,i<j}n_xr_{x,i}r_{x,j}.
\end{equation}

The final piece of notation we need is the following: let $A=\sum_{\gamma\in\Gamma_D}A_\gamma w^\gamma\in\sMot(\kk)[\bL^\frac12][[\Gamma_D]]$ be a motivic series and let $\Xi$ be a predicate on $\Gamma_D$ (that is, a condition, for example, a linear equation). Then $A_\Xi:=\sum_{\gamma\colon\Xi}A_\gamma w^\gamma$ stands for the sum of the monomials whose exponents satisfy $\Xi$.

By Corollary~\ref{cor:StabConn}, the stack $\Conn_{\gamma,d}(\epsilon,X,\divisor,\zeta)$ is of finite type, provided that $\epsilon\ne0$. By Corollary~\ref{cor:MellitSimpl} the coefficients of $\HH^{mot}_{X,D,\delta}$ are polynomials in $z$, so we may evaluate $\HH^{mot}_{X,D,\delta}$ at $z=1$.
\begin{theorem}[see Theorem~\ref{th:Main1}]\label{th:IntroMain}
Assume that $\epsilon\ne0$, that the class $\gamma'\in\Gamma_D$ is full at all $x\in D$ such that $n_x\ge2$, and that $\zeta\in\FNF(\divisor)$ is non-resonant for $\gamma'$ at all such points $x$. Assume $d'\in\Z$ is such that $\epsilon d'+\gamma'\star\zeta=0$. Then the motivic class of $\Conn_{\gamma',d'}(\epsilon,X,\divisor,\zeta)$ in $\sMot(\kk)$ is equal to the coefficient at $w^{\gamma'}$ in
\[
    (-\bL^\frac12)^{\chi(\gamma')} \Exp\left(\bL\left(\HH^{mot}_{X,D,\delta}\Bigl|_{z=1}\right)_{\gamma\star\zeta\in\epsilon\Z}\right),
\]
where the subscript $\gamma\star\zeta\in\epsilon\Z$ stands for the sum of the monomials whose exponents satisfy the condition.
\end{theorem}

\subsection{Explicit formulas for motivic classes of semistable parabolic $\epsilon$-connections}\label{sect:ExplAnswers2}
We can introduce stability conditions and consider motivic classes of the semistable loci $\Conn_{\gamma,d}(\epsilon,X,\divisor,\zeta)$. The semistable loci are of finite type (recall that full stacks have infinite motivic volumes when $\epsilon=0$). If $E$ is a vector bundle on $X$, its subbundle $F$ is called {\it strict}, if $E/F$ is torsion free. Let $(E,E_{\bullet,\bullet},\nabla)$ be a point of $\Conn_{\gamma,d}(\epsilon,X,\divisor,\zeta)$ and let $F\subset E$ be a strict subbundle preserved by $\nabla$. A~new difficulty in the irregular case is that we do not necessarily get an induced parabolic structure on~$F$, see Example~\ref{ex:DoesNotInduce}. The situation improves if we assume that $\gamma$ is full at all $x\in D$ such that $n_x\ge2$, and that $\zeta$ is non-resonant for $\gamma$ at all such $x$. We will see in Corollary~\ref{cor:Subbundle} that under this condition, we get an induced parabolic $\epsilon$-connection $(F,F\cap E_{\bullet,\bullet},\nabla|_F)\in\Conn_{\gamma',d'}(\epsilon,X,\divisor,\zeta)$ for some $(\gamma',d')\in\Gamma_D\times\Z$.

The sequence $\sigma=(\sigma_{\bullet,\bullet})$ of real numbers indexed by $D\times\Z_{>0}$ is called a \emph{sequence of parabolic weights of type $(X,\divisor)$\/} if for all $x\in D$ we have $\sigma_{x,1}\le\sigma_{x,2}\le\ldots$ and for all $j$ we have $\sigma_{x,j}\le\sigma_{x,1}+n_x$ (see Definition~\ref{def:ParWeights}). Similarly to~\eqref{eq:Intro_star}, we set $\gamma\star\sigma:=\sum_{x\in D}\sum_{j>0}r_{x,j}\sigma_{x,j}\in\R$. We call the parabolic $\epsilon$-connection $(E,E_{\bullet,\bullet},\nabla)$ as above \emph{$\sigma$-semistable\/} if for all strict subbundles $F\subset E$ preserved by $\nabla$, we have
\[
    \frac{d'+\gamma'\star\sigma}{\rk F}\le\frac{d+\gamma\star\sigma}{\rk E},
\]
where $(\gamma',d')=\cl(F,F\cap E_{\bullet,\bullet})$ (see Definition~\ref{def:Semistable}(i)). We denote by $\Conn^{\sigma-ss}_{\gamma,d}(\epsilon,X,\divisor,\zeta)$ the open substack of $\Conn_{\gamma,d}(\epsilon,X,\divisor,\zeta)$ classifying $\sigma$-semistable parabolic $\epsilon$-connections. By Corollary~\ref{cor:StabHiggs} these stacks are always of finite type.
\begin{theorem}[see Theorem~\ref{th:Main2}]\label{th:IntroMain2}
Assume that $\gamma'\in\Gamma_D$ is full at all $x\in D$ such that $n_x\ge2$, and that $\zeta\in\FNF(\divisor)$ is non-resonant for $\gamma'$ at all such $x$. Assume that $d'\in\Z$ is such that $\epsilon d'+\gamma'\star\zeta=0$. Let $\sigma$ be a sequence of parabolic weights of type $(X,\divisor)$. Then the motivic class of $\Conn^{\sigma-ss}_{\gamma',d'}(\epsilon,X,\divisor,\zeta)$ in $\sMot(\kk)$ is equal to the coefficient at $w^{\gamma'}$ in
\[
    (-\bL^\frac12)^{\chi(\gamma')}\Exp\left(\bL\left(\HH^{mot}_{X,D,\delta}\Bigl|_{z=1}\right)_{\substack{\gamma\star\zeta\in\epsilon\Z\\ \gamma\star\sigma-\tau\rk\gamma\in\Z }}\right),
\]
where $\tau:=(d'+\gamma'\star\sigma)/\rk\gamma'$.
\end{theorem}

\subsection{Parabolic $\epsilon$-connections with partially fixed normal forms}\label{sect:IntroIrregFNF}
Sometimes one is interested in $\epsilon$-connections with only partially fixed normal forms. The stacks of parabolic $\epsilon$-connections with partially fixed normal forms will also be crucial for us in calculating motivic classes of moduli stacks of parabolic $\epsilon$-connections with fully fixed normal forms; see Section~\ref{sect:IdeaOfCalc} and the proof of Proposition~\ref{pr:ExplicitFull-}. We sketch the definitions sending the reader to Section~\ref{sect:IrregFNF} for details. Let $X$, $D$, and $\divisor=\sum_{x\in D}n_xx$ be as before. Let $\divisor'$ be an effective divisor such that $\divisor'\le\divisor$, that is, $\divisor'=\sum_{x\in D}n'_xx$, where for all $x\in D$ we have $0\le n'_x\le n_x$. Set
\[
    \FNF(\divisor,\divisor'):=\FNF(\divisor)/\FNF(\divisor-\divisor')=\prod_{x\in D}(\Omega_X(n_xx)/\Omega_X((n_x-n'_x)x))^{\Z_{>0}}.
\]
For $\zeta\in\FNF(\divisor,\divisor')$ the stack $\Conn^{prtl}(\epsilon,X,\divisor,\divisor',\zeta)$ classifies triples $(E,E_{\bullet,\bullet},\nabla)$, where $(E,\nabla)$ is an~$\epsilon$-connection on $X$ with poles bounded by $\divisor$, the collection $E_{\bullet,\bullet}$ is a \emph{level $\divisor'$\/} parabolic structure on~$E$, a part of the normal form at $x\in D$ is given by $\zeta_{x,j}$, and we have a compatibility condition between $E_{\bullet,\bullet}$ and $\nabla$ similar to the one in Definition~\ref{def:ModSpace}, see Definition~\ref{def:ModSpacePartial} below. We define non-resonance of $\zeta$ for~$\gamma$ at $x\in D$ similarly to Definition~\ref{def:non-resonant0}. We note that our definition is very flexible: we may put no condition at some $x\in D$ (taking $n'_x=0$), fix the ``irregular part'' (taking $n'_x=n_x-1$), or just fix the leading term of the singularity (taking $n'_x=1$). Let $D':=\{x\in D\colon n'_x>0\}$ be the support of $\divisor'$ and set $\delta:=\deg\divisor-|D'|$. For $\gamma=(r,r_{\bullet,\bullet})\in\Gamma_{D'}$, define:
\begin{equation}\label{eq:BetterChi3}
\chi(\gamma)=\chi(\gamma,g,\divisor,\divisor'):=
(2g-2+\deg(\divisor-\divisor'))r^2+r\sum_{x\in D'}(1-n'_x)+2\sum_{x\in D',i<j}n'_xr_{x,i}r_{x,j}.
\end{equation}
We note that these $\delta$ and $\chi$ coincide with those defined in Section~\ref{sect:ExplAnswers} when $\divisor'=\divisor$. For a sequence~$\sigma$ of parabolic weights we define the open substack
\[
    \Conn^{prtl,\sigma-ss}_{\gamma,d}(\epsilon,\divisor,\divisor',\zeta) \subset\Conn^{prtl}_{\gamma,d}(\epsilon,\divisor,\divisor',\zeta)
\]
similarly to the case of fully fixed formal normal forms.
\begin{theorem}[see Theorem~\ref{th:Main3}]\label{th:IntroMain3}
Assume that $\divisor'<\divisor$, that is, $n'_x<n_x$ for at least one $x\in D$. Assume that $\gamma'\in\Gamma_{D'}$ is full at all $x\in D$ such that $n'_x\ge2$, and that $\zeta\in\FNF(\divisor,\divisor')$ is non-resonant for $\gamma'$ at all such $x$. Let $\sigma$ be a sequence of parabolic weights of type $(X,\divisor')$. Then the motivic class of $\Conn^{prtl,\sigma-ss}_{\gamma',d'}(\epsilon,X,\divisor,\divisor',\zeta)$ in $\sMot(\kk)$ is equal to the coefficient at $w^{\gamma'}$ in
\[
    (-\bL^\frac12)^{\chi(\gamma')}\Exp\left( \left(\HH^{mot}_{X,D,\delta}\Bigl|_{z=1}\right)_{\gamma\star\sigma-\tau\rk\gamma\in\Z}\right),
\]
where $\tau:=(d'+\gamma'\star\sigma)/\rk\gamma'$.
\end{theorem}

\subsection{Existence of connections} It is the 1-dimensional case of a~classical result of Atiyah~\cite{AtiyahConnections} that a vector bundle on a smooth projective curve admits a connection if and only if each its direct summand has degree zero. This was generalized to the parabolic case by Crawley-Boevey, see~\cite[Thm.~7.1]{Crawley-Boevey:Indecomposable}. His generalization is based on a result of Mihai~\cite[Thm.~1]{MihaiMonodromie,MihaiConnexions} (see~\cite[Thm.~7.2]{Crawley-Boevey:Indecomposable} for a~reformulation in the case of curves). In Section~\ref{sect:ExistencePartial} we address this question in the case of {\it irregular} parabolic connections. In other words, given a parabolic bundle $\bE=(E,E_{\bullet,\bullet})\in\Bun^{par}(X,\divisor')(\kk)$, $\zeta\in\FNF(\divisor,\divisor')$, and $\epsilon\in\kk$, we give a criterion for existence of a singular connection $\nabla\colon E\to E\otimes\Omega_X(\divisor)$ such that $(\bE,\nabla)\in\Conn^{prtl}(\epsilon,X,\divisor,\divisor',\zeta)(\kk)$. The answer is given in terms of the ring of endomorphisms of $\bE$. We present here the case $\divisor'=\divisor$ sending the reader to Theorem~\ref{th:Existence} for the general case.
\begin{theorem}\label{th:IntroExistence}
Let $\bE=(E,E_{\bullet,\bullet})\in\Bun^{par}(X,\divisor)(\kk)$ be a level $\divisor$ parabolic bundle on $X$. Let $\zeta\in\FNF(\divisor)$ and $\epsilon\in\kk$. Then $\bE$ can be extended to a $\kk$-point $(\bE,\nabla)$ of $\Conn(\epsilon,X,\divisor,\zeta)$ if and only if for every endomorphism $\Psi$ of $\bE$ we have
    \[
        \epsilon\langle b(E),\Psi\rangle-\sum_{x\in D}\sum_{j>0}\res_x(\zeta_{x,j}[\tr\gr_j\ev_x](\Psi))=0.
    \]
\end{theorem}
Here $b(E)\in H^1(X,\END(E)\otimes\Omega_X)$ is the Atiyah class of $E$, $\langle\bullet,\bullet\rangle$ is the Serre pairing, and $\gr_j\ev_x(\Psi)$ is the induced endomorphism of $E_{x,j-1}/E_{x,j}$, whose trace is an element of $\cO_X(n_xx)/\cO_X=\kk[n_xx]$. The new feature of the irregular situation is that for the existence of the connection it is not enough to assume that every direct summand has degree zero, see Example~\ref{ex:IrregNotSuff}.

\subsection{Interpretation of $\Conn(\epsilon,X,\divisor,\zeta)$ in terms of twisted cotangent bundles}\label{sect:Twisted} For $\gamma=(r,r_{\bullet,\bullet})\in\Gamma_D$ and $d\in\Z$ let $\overline{\Bun_{\gamma,d}^{par}(X,\divisor)}$ denote the moduli stack parameterizing families of collections $(E,E_{\bullet,\bullet},s_\bullet)$, where $(E,E_{\bullet,\bullet})$ is a point of $\Bun_{\gamma,d}^{par}(X,\divisor)$ and $s_x\colon (k[z]/z^{n_x})^r\to E_{n_xx}$ is a~trivialization compatible with the parabolic structure $E_{x,\bullet}$ as in Section~\ref{sect:ParStr} (here $x$ ranges over $D$). Then $\overline{\Bun_{\gamma,d}^{par}(X,\divisor)}\to\Bun_{\gamma,d}^{par}(X,\divisor)$ is a $G$-torsor, where $G=\prod_{x\in D}G_{r_{x,\bullet}}(k[z]/z^{n_x})$ and $G_{r_{x,\bullet}}\subset\GL_r$ is the group of block upper triangular matrices with blocks of sizes $r_{x,i}\times r_{x,j}$. We have
\[
    \Lie(G)_{ab}:=\Lie(G)/[\Lie(G),\Lie(G)]=\bigoplus_{x,j\colon r_{x,j}\ne0}k[z]/z^{n_x}.
\]
The residue pairing identifies the dual of this vector space with $\FNF_\gamma(\divisor):=\bigoplus_{x,j\colon r_{x,j}\ne0}\Omega_X(n_xx)/\Omega_X$. Let $\zeta_\gamma$ be the image of $\zeta$ under the projection $\FNF(\divisor)\to\FNF_\gamma(\divisor)$ forgetting the components $\zeta_{x,j}$ with $r_{x,j}=0$. Then we have the symplectic reduction $\left.\left.T^*\Bigl(\overline{\Bun_{\gamma,d}^{par}(X,\divisor)}\Bigr)\right/\!\!\right/_{\!\!\zeta_\gamma}\!\!G$. It is easy to see that it is identified with $\Conn_{\gamma,d}(0,X,\divisor,\zeta)$.

To take care of the case $\epsilon\ne0$, one has to consider the determinantal line bundle $det\to\Bun(X)$. Let $det^\times$ be the complement of the zero section in the total space of $det$. Note that the fiber product $\overline{\Bun_{\gamma,d}^{par}(X,\divisor)}\times_{\Bun(X)}det^\times$ is a $G\times{\mathbb G}_m$-torsor over $\Bun_{\gamma,d}^{par}(X,\divisor)$. Similarly to the above, we have the symplectic reduction
\[
    \left.\left.T^*\left(\overline{\Bun_{\gamma,d}^{par}(X,\divisor)}\times_{\Bun(X)}det^\times\right)
    \right/\!\!\right/_{\!\!(\zeta_\gamma,\epsilon)}(G\times{\mathbb G}_m)
    \simeq\Conn_{\gamma,d}(\epsilon,X,\divisor,\zeta).
\]

\subsection{Overview of the argument}\label{sect:IdeaOfCalc} The idea of calculation of motivic classes is the following. Assume that $\divisor$ is an irregular divisor, that is, $n_x\ge2$ for some $x\in D$. Similarly to Schiffmann (see~\cite{SchiffmannIndecomposable}), we consider the open substack $\Conn_{\gamma,d}^{prtl,-}(\epsilon,X,\divisor,\divisor',\zeta)$ classifying parabolic $\epsilon$-connections whose underlying vector bundle has no subbundles of positive degree (we call such vector bundles \emph{nonpositive}). The advantage of this stack is that it is of finite type.

Then we proceed as follows. We start with the stack $\Conn^{prtl,-}(\epsilon,X,\divisor,\divisor',\zeta)$ in the case when~$\divisor'$ is a \emph{reduced divisor}. In this case, we show that every parabolic bundle admits a connection (see Lemma~\ref{lm:NoAtiyah}). Thus, it is easy to compare the motivic class of this stack with the motivic class of the stack $\Pair_{\gamma,d}^{nilp,-}(X,\divisor',\cL)$ classifying pairs $(\bE,\Psi)$, where $\bE$ is a parabolic bundle of level $\divisor'$ on $X$ with a nonpositive underlying vector bundle and $\Psi:\bE\to\bE\otimes\cL$ is a twisted nilpotent endomorphism, see Section~\ref{sect:CaseOfReduced}. Here $\cL=\cO_X(\divisor'-\divisor)$.

The motivic classes of the stacks $\Pair_{\gamma,d}^{nilp,-}(X,\divisor',\cL)$ are calculated in Section~\ref{Sect:parabolic twisted} as follows. First of all, the set of equivalence classes of points of the stack $\Pair^{nilp,-}_{\gamma,d}(X,\divisor',\cL)$ decomposes into a disjoint union:
\[
    |\Pair^{nilp,-}_{\gamma,d}(X,\divisor',\cL)|= \bigcup_{\mu\,\vdash\rk\gamma}\Pair^{nilp,-}_{\gamma,d}(X,\divisor',\mu,\cL).
\]
Here $\Pair^{nilp,-}_{\gamma,d}(X,\divisor',\mu,\cL)$ is a constructible subset of $\Pair^{nilp,-}_{\gamma,d}(X,\divisor',\cL)$ corresponding to the pairs with the generic Jordan type $\mu$ (see Section~\ref{sect:Notation} for the notation). Next, we view the disjoint union
\[
    \Pair^{nilp,-}(X,\divisor',\mu,\cL):=\bigsqcup_{\gamma\in\Gamma_{D'},d\in\Z}\Pair^{nilp,-}_{\gamma,d}(X,\divisor',\mu,\cL)
\]
as a $\Gamma_{D'}\times\Z$-graded stack, so that its motivic class is an element of $\cMot(\kk)[[\Gamma_{D'}\times\Z]]$. Then we have the following product formula:
\begin{equation}\label{eq:Product}
  [\Pair^{nilp,-}(X,\divisor',\mu,\cL)]=[\Pair^{nilp,-}(X,\divisor',\emptyset,\cL)]V(\mu)^{|\divisor'|}.
\end{equation}
The significance of this formula is that the first multiple in the RHS corresponds to the case when there are no marked points, while the multiples $V(\mu)$ are independent of the curve $X$ and on the divisor. In the case when $\cL=\cO_X$, this follows from~\cite[Thm.~4.4]{FedorovSoibelmansParabolic}; the general case is proved similarly. The motivic class of $\Pair^{nilp,-}(X,\divisor',\emptyset,\cL)$ is calculated in Section~\ref{sect:TwistedPairs} below (this is a~relatively simple extension of results of~\cite{MozgovoySchiffmann2020}). The local factors are independent of $\cL$; they have been calculated in~\cite[Sect.~4]{FedorovSoibelmansParabolic}. The calculation of the local factors is based on considering the simplest case $X=\P^1_\kk$, $\divisor'=0+\infty$, see~\cite{MellitPunctures}. We note that this calculation has been generalized to arbitrary split reductive groups in~\cite{Singh2021Counting}.

Once the motivic classes of $\Conn^{prtl,-}(\epsilon,X,\divisor,\divisor',\zeta)$ are available for reduced $\divisor'$, we proceed by induction on $\deg\divisor'$ to calculate these motivic classes for any $\divisor'<\divisor$. The idea is to relate the motivic class above to the motivic class of $\Conn^{prtl,-}(\epsilon,X,\divisor,\divisor'-x,\zeta)$, where $x\in D$ is any point with $n'_x\ge2$, see Lemma~\ref{lm:AffineBundle2}. The case $\divisor'=\divisor$ is more difficult. The idea is to relate $\Conn^-(\epsilon,X,\divisor,\zeta)$ to
$\Conn^{prtl,-}(\epsilon,X,\divisor,\divisor-x,\zeta)$. We use Lemma~\ref{lm:existIrregFull} claiming that a parabolic bundle $\bE$ admits an $\epsilon$-connection with formal normal form $\zeta$, if and only if it admits such an $\epsilon$-connection with the normal form fixed to the level $\divisor-x$ \emph{and\/} every direct summand of $\bE$ has degree zero. With this lemma at hand, the argument, given in the second part of Section~\ref{sect:MotClass-}, proceeds similarly to the regular case, see~\cite[Sect.~5 and~8]{FedorovSoibelmansParabolic}.

The rest of the calculation is similar to the regular case. If $\epsilon\ne0$ and $\divisor'=\divisor$, the stack $\Conn(\epsilon,X,\divisor,\zeta)$ is already of finite type, and we calculate its motivic class using a limit argument, see Section~\ref{sect:StabConn}. In the remaining cases, we impose a stability condition first, see Section~\ref{sect:Stability}, and then pass to the limit in Section~\ref{sect:RemainingCase}.

\subsection{E-polynomials and virtual Poincar\'e polynomials} In Section~\ref{sect:EPoly} we use our calculation of motivic classes of stacks of parabolic $\epsilon$-connections to calculate the E-polynomials and the virtual Poincar\'e polynomials of these stacks. As a result (see Remark~\ref{rem:Pantev}), we provide evidence for a conjecture of Diaconescu, Donagi, and Pantev (see~\cite[(1.8)]{DiaconescuDonagiPantev}), which, in turn, is based on a conjecture of Hausel, Mereb, and Wang (see~\cite[Conj.~0.1.1]{HauselMerebWong}).

\subsection{A few words about possible future directions}\label{sect:FutureDirections} Another approach to counting Higgs bundles over finite fields based on orbital integrals and trace formulas was developed in~\cite{LaumonChaudouard}. One can hope that this approach can lead to precise formulas, and moreover that it can be upgraded to the motivic setup as well as to the cases of parabolic and irregular parabolic Higgs bundles and connections.

Recall from Section~\ref{sect:IdeaOfCalc} the stack $\Pair^{nilp}(X,\divisor',\mu,\cL)$ of nilpotent pairs. Although we are considering the irregular case, we only need to consider the above stacks when $\divisor'$ is a reduced divisor. This is related to the fact that we only allow non-resonant formal normal forms $\zeta$. More generally, it is a~challenging problem to calculate the motivic classes of the stacks $\Pair^{nilp,-}(X,\divisor',\mu,\cL)$ when $\divisor'$ is a not necessarily reduced effective divisor. We expect that the factorization formula~\eqref{eq:Product} still holds in this case. However, the local terms depend on the multiplicities of points in $\divisor'$. Since the local terms do not depend on the curve, it seems plausible that, as in the reduced case, we can calculate our motivic classes by considering the simplest case of $X=\P^1_\kk$, $\divisor'=n0+m\infty$, where $n$ and $m$ are positive integers. We do not know how to calculate the motivic classes of the corresponding stacks.

Finally, one can approach the problem using cohomological Hall algebras (COHAs for short) introduced in~\cite{KontsevichSoibelman10} instead of motivic Hall algebras used in~\cite{FedorovSoibelmans}. Indeed, it was proven in~\cite{KontsevichSoibelman10} that at the level of motivic invariants the answers obtained via motivic Hall algebras and via COHAs agree. Concerning the problem we are interested in, one can work either with the 2-dimensional COHA directly associated with the category of Higgs bundles (see~\cite{sala2018cohomological}) or with the 3-dimensional COHA associated with the corresponding Calabi--Yau 3-fold. The two approaches agree, since the 2-dimensional COHA is obtained from 3-dimensional COHA via dimensional reduction (see~\cite[Sect.~4.8]{KontsevichSoibelman10}). The approach via COHAs can give more than just the computation of motivic classes. Indeed many generating functions of enumerative invariants admit interpretations in terms of the representation theory of COHAs (see, e.g.,~\cite{SchiffmannVasserot2013},~\cite{Soibelman2016remarks},%
~\cite{RapcakSoibelmanYangZhao2020},~\cite{MellitMinetsSchiffmannVasserot2023} for a few samples of ideas in that direction). We plan to return to this topic in the future.

\subsection{Notations and conventions}\label{sect:Notation}
We denote by $\kk$ a field of characteristic zero. We denote by~$X$ a smooth projective geometrically connected curve over $\kk$ (recall that geometric connectedness means that $X$ remains connected after the base change to an algebraic closure of $\kk$). We denote by $D$ a~possibly empty set of $\kk$-rational points of $X$ and fix a divisor $\divisor=\sum_{x\in D}n_xx$, where for all $x\in D$ we have $n_x>0$. If $x\in X(\kk)$ is a rational point, we can always find an \'etale morphism $z\colon X'\to\A^1_\kk$, where $X'\subset X$ is a Zariski open neighborhood of $x$ and $z(x)=0$. We call such $z$ \emph{a coordinate near $x$.} In fact, formal coordinates would suffice for most of our purposes as well.

If $E$ is a vector space (resp.~vector bundle), we denote by $E^\vee$ the dual vector space (resp.~vector bundle). We sometimes identify vector bundles with their sheaves of sections. If $E$ is a vector bundle on a scheme $S$ and $s\in S$ is a schematic point, we can trivialize $E$ in a Zariski open neighborhood of~$s$. We call this trivialization a \emph{local trivialization\/} or a \emph{trivialization near $s$}.

We systematically work with Artin stacks that are locally of finite type over a field $\kk$. All our stacks have affine centralizers. If $\cX$ is a $\kk$-stack and $S$ is a $\kk$-scheme, then an $S$-point $\xi$ of $\cX$ is an object of the groupoid $\cX(S)$. In this case, we write $\xi\in\cX(S)$. By ``a point of $\cX$'', we mean an $S$-point for some test $\kk$-scheme $S$.

We denote by $|\cX|$ the set of equivalence classes of $K$-points of a stack $\cX$, where $K$ ranges over all field extensions of $\kk$, see~\cite[Ch.~5]{LaumonMoretBailly} and~\cite[Sect.~2]{FedorovSoibelmans}. Let $\cX$ be a stack of finite type over $\kk$. If $\cS\subset|\cX|$ is a constructible subset, we can define its motivic class $[\cS]$. This is a slight abuse of notation, because $[\cS]$ depends both on $\cS$ and $\cX$. By a stratification we mean a decomposition of $|\cX|$ into the disjoint union of constructible subsets: $|\cX|=\bigcup_{i=1}^n\cS_i$. In this case we have in $\Mot(\kk)$: $[\cX]=\sum_{i=1}^n[\cS_i]$.

\subsection{Acknowledgements} The authors thank E.~Diaconescu, D.~Kazhdan, O.~Kivinen, R.~Li, T.~Pantev, and L.~Rozansky for valuable discussions. The work of R.F.~was partially supported by NSF grants DMS--2001516 and DMS--2402553. The work of A.S. and Y.S. was partially supported by ERC--SyG project, Recursive and Exact New Quantum Theory (ReNewQuantum), which received funding from the European Research Council (ERC) under the European Union's Horizon 2020 research and innovation program
under grant agreement No 810573. A part of this work was done when R.F., A.S., and Y.S. visited IHES. We thank IHES for the hospitality and excellent research conditions.

\section{Parabolic \texorpdfstring{$\epsilon$}{epsilon}-connections and parabolic bundles}
Recall that $X$ is a smooth projective geometrically connected curve over $\kk$, $\divisor=\sum_{x\in D}n_xx$ is an effective divisor with support $D\subset X(\kk)$ and $\epsilon\in\kk$. We are studying $\epsilon$-connections on $X$ with poles bounded by $\divisor$; see Section~\ref{sect:IntroEpsConn}.

\subsection{Generic $\epsilon$-connections}
\begin{definition}\label{def:generic}
  Let $(E,\nabla)$ be an $\epsilon$-connection on $X$ with poles bounded by $\divisor$. We say that $(E,\nabla)$ is \emph{generic\/} at $x\in D$ if for a choice of a coordinate $z$ near $x$ and a local trivialization of~$E$ near $x$ we have
  \[
    \nabla=\epsilon d+\frac A{z^{n_x}}\,dz+\text{ lower order terms},
  \]
  where $A$ is a matrix diagonalizable over $\kk$ with \emph{distinct eigenvalues}.
\end{definition}

This condition is independent of the choice of local coordinate and of the local trivialization of $E$: indeed, different choices will result in conjugating and scaling of $A$.

\begin{remark}
Keeping notation as in the definition, assume that $x$ is an irregular singular point, that is, $n_x\ge2$. Then $\nabla$ is called {\it non-resonant} at $x$ if $A$ is regular semisimple (if $n_x=1$, the definition is slightly different: one requires that no eigenvalues of $A$ differ by a positive integer). Thus, every generic $\epsilon$-connection is non-resonant and the two notions coincide if $\kk$ is algebraically closed, provided that $n_x\ge2$.
\end{remark}

\begin{lemma}\label{lm:diag}
Let an $\epsilon$-connection $(E,\nabla)$ with poles bounded by $\divisor$ be generic at $x\in D$. Then there is a unique up to permutation and scaling basis of the $\kk[n_xx]$-module $E_{n_xx}$ such that after extending this basis to a local trivialization of $E$ near $x$, we obtain
\begin{equation}\label{eq:diag}
    \nabla=\epsilon d+\diag(\tilde\zeta_1,\ldots,\tilde\zeta_{\rk E})+B,
\end{equation}
where $\tilde\zeta_i$ is a local section of $\Omega_X(n_xx)$, and the matrix $B$ has no pole at $x$ (here $\diag(\bullet)$ stands for the diagonal matrix with given diagonal entries).
\end{lemma}
\begin{proof}
  This is well-known but we give a proof. In fact for any $n=1,\ldots,n_x$, we show a similar statement about the $\kk[nx]$-module $E_{nx}$ with the only difference that $B$ is allowed to have a pole of order at most $n_x-n$. When $n=n_x$ we obtain the statement of the lemma.

We proceed by induction on $n$. If $n=1$, this reduces to the claim that if $A$ is a matrix with distinct eigenvalues, then there is a unique up to permutation and scaling basis diagonalizing the matrix. (Here~$A$ as in Definition~\ref{def:generic}.)

For the induction step we choose a coordinate $z$ near $x$. Choose any trivialization of $E$ near~$x$ and write
  \[
    \nabla=\epsilon d+\sum_{j=-n_x}^{\infty}A_jz^j\,dz,
  \]
where $A_j$ are constant $\rk E\times\rk E$ matrices. First we prove the existence. By the induction hypothesis we may assume that $A_j$ is diagonal for $-n_x\le j\le n-1-n_x$ and we need to show that changing a~trivialization we may achieve that $A_{n-n_x}$ is diagonal as well. Since $A=A_{-n_x}$ has distinct diagonal entries, we can find $B$ such that $A_{n-n_x}-[A,B]$ is diagonal. The chosen local trivialization of $E$ gives a trivialization of $E$ over the infinitesimal neighborhood $(n+1)x$. Now we use $\Phi:=\Id+Bz^n\in\GL_n(\kk[[z]]/(z^n))$ to change this trivialization. Then $A_{n-n_x}$ becomes replaced with the diagonal matrix $A_{n-n_x}-[A,B]$, and we are done.

Uniqueness is proved along the same lines using the fact that if $A$ is a diagonal matrix with distinct diagonal entries and $[A,B]$ is diagonal, then $B$ is diagonal as well.
\end{proof}

\subsection{Parabolic structures preserved by generic $\epsilon$-connections}\label{sect:ParStrPreserved}
Let $(E,\nabla)$ be an $\epsilon$-connection with poles bounded by $\divisor=\sum_{x\in D}n_xx$. For $n=1,\ldots,n_x$, we say that a free submodule $F$ of the free $\kk[nx]$-module $E_{nx}$ is \emph{preserved by $\nabla$\/} if for every local section~$s$ of~$E$ such that $s|_{nx}\in F$ we have
\[
    \nabla(s)|_{(E\otimes\Omega_X(n_xx))_{nx}}\in F\otimes_{\kk[nx]}\Omega_X(n_xx)_{nx}.
\]
We say that a level $n$ parabolic structure $E_\bullet=(E_{nx}\supset E_1\supset E_2\supset\ldots)$ at $x$ is preserved by $\nabla$ if each~$E_j$ is preserved (cf.~Definition~\ref{def:thick}).

One can rewrite this definition ``in coordinates'' as follows: consider a local trivialization of $E$ compatible with the parabolic structure. (Recall that this means that each $E_i$ is generated by $e_1,\ldots,e_{\rk E_i}$, where $e_1,\ldots,e_{\rk E}$ is the basis provided by the trivialization, see Section~\ref{sect:ParStr}.) Then $\nabla=\epsilon d+A+B$, where $A$ is a block upper triangular matrix with blocks of sizes $\rk(E_{x,i-1}/E_{x,i})\times\rk(E_{x,j-1}/E_{x,j})$ and $B$ has pole of order at most $n_x-n$.

In the case when $\nabla$ is generic at $x$, there are only finitely many submodules (and thus countably many parabolic structures) of $E_{nx}$ preserved by $\nabla$ in view of the following result.
\begin{lemma}\label{lm:NablaPreservedSubmodules}
Let an $\epsilon$-connection $(E,\nabla)$ with poles bounded by $\divisor$ be generic at $x\in D$. Consider a basis of $E_{n_xx}$, provided by Lemma~\ref{lm:diag}, and let $(e_1,\ldots,e_{\rk E})$ be the induced basis of $E_{nx}$, where $1\le n\le n_x$. Then a free submodule $F\subset E_{nx}$ is preserved by $\nabla$ if and only if it is spanned by a subset of $\{e_1,\ldots,e_{\rk E}\}$.
\end{lemma}
\begin{proof}
Fix a coordinate $z$ near $x$ and use it to trivialize the rank one free $k[nx]$-module $\Omega_X(n_xx)_{nx}$. Then $\nabla$ induces a $k[nx]$-linear operator $\nabla'\colon E_{nx}\to E_{nx}$ and $F$ is preserved by $\nabla$ if and only if $\nabla'(F)\subset F$. Let $(e_1,\ldots,e_{\rk E})$ be a basis of $E_{nx}$ as in the formulation of the lemma, the matrix of~$\nabla'$ in this basis is diagonal with diagonal entries $\zeta'_i$ being the reductions of $z^{n_x}\tilde\zeta_i$ modulo $z^n$. We note that since $\nabla$ is generic at $x$, $\zeta'_i-\zeta'_j$ is invertible in $k[nx]$ whenever $i\ne j$. Now our claim is equivalent to the statement that the invariant submodules of such an endomorphism $\nabla'$ are exactly the submodules generated by elements of the diagonalizing basis. This is proved similarly to the corresponding well-known statement about linear operators in vector spaces (which is a particular case $n=1$ of our statement).
\end{proof}
The following simple corollary will be important for our inductive approach to calculating motivic classes of parabolic $\epsilon$-connections (see Section~\ref{sect:IdeaOfCalc}, Lemma~\ref{lm:AffineBundle2}, and Lemma~\ref{lm:RR}).

\begin{corollary}\label{cor:ExtParStr}
Let an $\epsilon$-connection $(E,\nabla)$ with poles bounded by $\divisor$ be generic at $x\in D$. Assume that we have chosen a level~1 parabolic structure on~$E$ at~$x$ preserved by $\nabla$. Then for $n=2,\ldots,n_x$, there is exactly one way to extend this parabolic structure to a level $n$ parabolic structure at $x$ preserved by~$\nabla$.
\end{corollary}
\begin{proof}
It follows from the above lemma, that a vector subspace of $E_x$ preserved by $\nabla$ can be uniquely extended to a free submodule of $E_{nx}$ preserved by $\nabla$. The statement follows.
\end{proof}

\subsection{Parabolic $\epsilon$-connections with partially fixed normal forms}\label{sect:IrregFNF}
We are going to define the stacks of parabolic $\epsilon$-connections with partially fixed normal forms (cf.~Section~\ref{sect:IntroIrregFNF}). This will be our main object of study. Let $\divisor'$ be an effective divisor such that $\divisor'\le\divisor$, that is, $\divisor'=\sum_{x\in D}n'_xx$, where for all $x\in D$ we have $0\le n'_x\le n_x$. Recall from Section~\ref{sect:IntroIrregFNF} the vector space
\begin{equation}\label{eq:FNF_def}
    \FNF(\divisor,\divisor'):=\FNF(\divisor)/\FNF(\divisor-\divisor')=\prod_{x\in D}(\Omega_X(n_xx)/\Omega_X((n_x-n'_x)x))^{\Z_{>0}}.
\end{equation}
\begin{definition}\label{def:ModSpacePartial}
  Let $\zeta\in\FNF(\divisor,\divisor')$. Then $\Conn^{prtl}(\epsilon,X,\divisor,\divisor',\zeta)$ is the moduli stack classifying triples $(E,E_{\bullet,\bullet},\nabla)$, where\\
  $\bullet$ $(E,\nabla)$ is an $\epsilon$-connection on $X$ with poles bounded by $\divisor$;\\
  $\bullet$ $E_{\bullet,\bullet}$ is a level $\divisor'$ parabolic structure on $E$;\\
  $\bullet$ for every $x\in D$ and every local trivialization of $E$ near $x$ compatible with the parabolic structure, we have in this trivialization:
    \[
        \nabla=\epsilon d+A+B,
    \]
    where $A$ is a block upper triangular matrix with blocks of sizes $\rk(E_{x,i-1}/E_{x,i})\times\rk(E_{x,j-1}/E_{x,j})$ such that the $i$-th diagonal block is equal to $\tilde\zeta_{x,i}\Id$, where $\tilde\zeta_{x,i}$ is any lift of $\zeta_{x,i}$ to a local section of $\Omega_X(n_xx)$, and $B$ has a pole of order at most $n_x-n'_x$ at $x$.
\end{definition}
It is easy to see that the above compatibility condition is satisfied for every local trivialization of~$E$ near~$x$ compatible with the parabolic structure if it is satisfied for a single such trivialization. Note that $\FNF(\divisor,\divisor)=\FNF(\divisor)$ and $\Conn^{prtl}(\epsilon,X,\divisor,\divisor,\zeta)=\Conn(\epsilon,X,\divisor,\zeta)$, where the latter stack is defined in Definition~\ref{def:ModSpace}. We will mostly omit $X$ and $\epsilon$ from the notation as they are now fixed.

Recall that in Section~\ref{sect:GammaD} we defined a monoid $\Gamma_{D'}$ for any finite set $D'$. Let $D':=\{x\in D\colon n'_x>0\}$ be the support of $\divisor'$. Then we have a decomposition into closed and open subsets
\[
    \Conn^{prtl}(\divisor,\divisor',\zeta)=
    \bigsqcup_{\gamma\in\Gamma_{D'},d\in\Z}\Conn^{prtl}_{\gamma,d}(\divisor,\divisor',\zeta).
\]
If $(E,E_{\bullet,\bullet},\nabla)$ is a point of $\Conn^{prtl}(\divisor,\divisor',\zeta)$, then $\nabla$ preserves $E_{\bullet,\bullet}$ in the sense of Section~\ref{sect:ParStrPreserved}. However, the condition in Definition~\ref{def:ModSpacePartial} is in general stronger than just the requirement that the connection preserves the parabolic structure, since the diagonal blocks are required to be scalar.

The following definition is a generalization of Definition~\ref{def:non-resonant0}.
\begin{definition}\label{def:non-resonant}
  Let $\gamma=(r,r_{\bullet,\bullet})\in\Gamma_{D'}$ and $\zeta\in\FNF(\divisor,\divisor')$. We will say that $\zeta$ is
  \emph{non-resonant for~$\gamma$ at $x\in D$\/} if for all $1\le i<j$ such that $r_{x,i}\ne0\ne r_{x,j}$ the polar part $\zeta_{x,i}-\zeta_{x,j}$ has a pole of order exactly $n_x$.
\end{definition}
Recall the definition of full class (Definition~\ref{def:full}). The next lemma follows from the definitions.
\begin{lemma}\label{lm:DiagIffNonResonant}
  Let $(E,E_{\bullet,\bullet},\nabla)\in\Conn_{\gamma,d}^{prtl}(\divisor,\divisor',\zeta)(\kk)$. Let $x\in D$ be such that $n_x\ge2$. Then $(E,\nabla)$ is generic at $x$ if and only if $\gamma$ is full at $x$ and $\zeta$ is non-resonant for $\gamma$ at $x$.
\end{lemma}

We will denote by $\Gamma_{\divisor',\zeta}$ the subset of all $\gamma\in\Gamma_{D'}$ such that $\gamma$ is full at all $x\in D'$ with $n'_x\ge2$ and~$\zeta$ is non-resonant for $\gamma$ at all such $x$. The following lemma also follows from the definitions.

\begin{lemma}\label{lm:Summand}
  Assume that $\gamma_1,\gamma_2\in\Gamma_{D'}$ are such that $\gamma_1+\gamma_2\in\Gamma_{\divisor',\zeta}$. Then $\gamma_1,\gamma_2\in\Gamma_{\divisor',\zeta}$.
\end{lemma}

\begin{remark}\label{rem:NotNeeded}
    Assume that $\gamma\in\Gamma_{\divisor',\zeta}$. Assume that $(E,E_{\bullet,\bullet},\nabla)$ and $(E,E'_{\bullet,\bullet},\nabla)$ are two $\kk$-points of $\Conn_{\gamma,d}^{prtl}(\divisor,\divisor',\zeta)$ such that for all $x\in D$ with $n'_x=1$ we have $E_{x,\bullet}=E'_{x,\bullet}$. Then $(E,E_{\bullet,\bullet},\nabla)=(E,E'_{\bullet,\bullet},\nabla)$ (use Lemma~\ref{lm:NablaPreservedSubmodules}). In other words, the parabolic structures are uniquely determined at all points where $n'_x\ge2$.
\end{remark}
Note that we will calculate our motivic classes only in the case when $\gamma\in\Gamma_{\divisor',\zeta}$ (cf.~Sections~\ref{sect:ExplAnswers} and~\ref{sect:ExplAnswers2}). Thus, according to the above lemma, we will be working with $\epsilon$-connections that are generic at all points with $n'_x\ge2$. In the case $\divisor'=\divisor$, this means that our connections are generic at all irregular points. The problem of calculating the motivic classes of the stacks $\Conn^{prtl}_{\gamma,d}(\divisor,\divisor',\zeta)$ when~$\gamma$ is not full at some points with $n'_x\ge2$ remains a challenge.

\subsection{Stack of $\epsilon$-connections as an ``almost affine bundle''}\label{sect:AlmostAffine} We have a forgetful morphism $\Conn^{prtl}(\divisor,\divisor',\zeta)\to\Bun^{par}(\divisor')$ sending $(E,E_{\bullet,\bullet},\nabla)$ to $(E,E_{\bullet,\bullet})$. The fibers of this morphism are either empty or affine spaces (in general, of different dimensions). The goal of this section is to describe these fibers and to calculate their dimensions.

Let $\bE=(E,E_{\bullet,\bullet})$ be a level $\divisor'$ parabolic bundle. Let $\END(\bE)$ be the sheaf of endomorphism of~$\bE$, that is, of endomorphisms $\Psi$ of $E$ such that for all $x\in D$ and $j>0$ we have $\Psi(E_{x,j})\subset E_{x,j}$. Let $\END_0(\bE)$ be the subsheaf consisting of endomorphisms such that for all $x\in D$ and $j>0$ we have $\Psi(E_{x,j-1})\subset E_{x,j}$.

\begin{lemma}\label{lm:fibers}
  The fiber of the forgetful morphism $\Conn^{prtl}(\divisor,\divisor',\zeta)\to\Bun^{par}(\divisor')$ over a parabolic bundle $\bE\in\Bun^{par}(\divisor')(K)$, where $K\supset\kk$ is a field extension, is either empty or a trivial torsor over the vector space $H^0(X_K,\END_0(\bE)\otimes\Omega_{X_K}(\divisor))$.
\end{lemma}
\begin{proof}
Let $S$ be a test $K$-scheme. Let $(E,E_{\bullet,\bullet},\nabla)\in\Conn^{prtl}(\divisor,\divisor',\zeta)(S)$ and $\Phi\in H^0((X_S,\END_0(\bE)\otimes\Omega_{X_S/S}(\divisor))$, where $X_S:=X_K\times_KS$. Then it follows from Definition~\ref{def:ModSpacePartial} that $(E,E_{\bullet,\bullet},\nabla+\Phi)$ is also in $\Conn^{prtl}(\divisor,\divisor',\zeta)(S)$. Conversely, if $(E,E_{\bullet,\bullet},\nabla)$ and $(E,E_{\bullet,\bullet},\nabla')$ are $S$-points of $\Conn^{prtl}(\divisor,\divisor',\zeta)$, then $\nabla'-\nabla\in H^0(X_S,\END_0(\bE)\otimes\Omega_{X_S/S}(\divisor))$. We see that if the fiber under consideration is non-empty, then it is a torsor under $H^0(X_K,\END_0(\bE)\otimes\Omega_{X_K}(\divisor))$, where we view this vector space as a $K$-group scheme. It remains to note that every torsor over a field under such a group scheme is trivial.
\end{proof}

In particular, this lemma implies that every non-zero fiber of the morphism has a rational point. That is, if the fiber of $\Conn^{prtl}(\divisor,\divisor',\zeta)\to\Bun^{par}(\divisor')$ over $(E,E_{\bullet,\bullet})\in\Bun^{par}(\divisor')(K)$ is non-empty, then $(E,E_{\bullet,\bullet})$ can be extended to $(E,E_{\bullet,\bullet},\nabla)\in\Conn^{prtl}(\divisor,\divisor',\zeta)(K)$.

Finally, we relate the dimensions of the fibers to twisted endomorphisms. Recall that $g$ denotes the genus of $X$.
\begin{lemma}\label{lm:RR-dim}
Let $\bE\in\Bun_{\gamma,d}^{par}(\divisor')(K)$, where $K\supset\kk$ is a field extension and $\gamma=(r,r_{\bullet,\bullet})\in\Gamma_{D'}$. Then
\begin{multline}\label{eq:RR}
    \dim H^0(X_K,\END_0(\bE)\otimes\Omega_{X_K}(\divisor))-\dim H^0(X_K,\END(\bE)\otimes\cO_{X_K}(\divisor'-\divisor))\\=
    (g-1+\deg(\divisor-\divisor'))r^2+\sum_{x\in D}\sum_{i<j}n'_xr_{x,i}r_{x,j}.
\end{multline}
\end{lemma}
\begin{proof}
We have a non-degenerate $\cO_{X_K}$-bilinear trace pairing $\END_0(\bE)\otimes\END(\bE)\to\cO_{X_K}(-\divisor')$. Thus, by Serre duality the vector spaces $H^0(X_K,\END_0(\bE)\otimes\Omega_{X_K}(\divisor))$ and $H^1(X_K,\END(\bE)\otimes\cO_{X_K}(\divisor'-\divisor))$ are dual, so it remains to apply the Riemann--Roch theorem. The calculation is similar to that of~\cite[Lemma~5.5]{FedorovSoibelmansParabolic}.
\end{proof}

\section{Existence of \texorpdfstring{$\epsilon$}{epsilon}-connections}\label{sect:ExistencePartial}
As before, we are fixing a curve $X$, a subset $D\subset X(\kk)$, and effective divisors $\divisor'=\sum_{x\in D}n'_xx$ and $\divisor=\sum_{x\in D}n_xx$ such that for all $x\in D$ we have $0\le n'_x\le n_x$, $0<n_x$. We are also fixing $\epsilon\in\kk$. Recall that in Definition~\ref{def:ModSpacePartial} we defined the moduli stacks $\Conn^{prtl}(\epsilon,X,\divisor,\divisor',\zeta)$ of singular $\epsilon$-connections with partially fixed normal forms. We will now omit $\epsilon$ and $X$ from the notation. In Section~\ref{sect:AlmostAffine} we studied the fibers of the projection $\Conn^{prtl}(\divisor,\divisor',\zeta)\to\Bun^{par}(X,\divisor')$. The main result of the current section is Theorem~\ref{th:Existence}, where we describe the image of this projection. In the regular case, that is, when $n_x=1$ for all $x\in D$, our theorem is just a slight generalization of results of Mihai~\cite[Thm.~1]{MihaiConnexions} and Crawley-Boevey~\cite[Thm.~7.2]{Crawley-Boevey:Indecomposable}. In the irregular case, we are only aware of a simple particular case~\cite[Cor.~4.4]{ArinkinFedorov}.

We will also derive some consequences of Theorem~\ref{th:Existence}. In particular, we will see that existence of an $\epsilon$-connection on a parabolic bundle $\bE$ implies that the generalized degree of every direct summand of $\bE$ is zero (see Corollary~\ref{cor:indecomposable summand}). However, contrary to the regular case~\cite[Thm.~7.1]{Crawley-Boevey:Indecomposable}, this condition is not sufficient (see Example~\ref{ex:IrregNotSuff}).

\subsection{Infinitesimal neighborhoods and residue pairings}\label{sect:Infinitesimal} Let $x\in X$ be a point and $n$ be a positive integer. Then we have the $n$-th infinitesimal neighborhood $nx$ of $x$. This is the closed subscheme of~$X$ corresponding to the ideal sheaf $\cO_X(-nx)$. Choosing a coordinate $z$ near $x$, we can identify its ring of functions $\kk[nx]$ with $\kk(x)[z]/z^n$, though we prefer to work coordinate free (here $\kk(x)$ is the residue field of $x$). Let $F$ be a vector bundle on $X$. Then its restriction to $nx$ is denoted $F_{nx}$. It is a free $\kk[nx]$-module.

For $m\in\Z$ set $F(mx):=F\otimes\cO_X(mx)$. For $n>0$ we identify $F_{nx}$ with $F/F(-nx)$. More generally, for any $n>0$ and $m\in\Z$ we have the $k[nx]$-module $F(mx)_{nx}=F(mx)/F((m-n)x)$. If $n\le m$, we can view this space as the space of partial polar parts of singular sections of $F$. Let $F^\vee$ be the dual vector bundle. Then we have a canonical perfect residue pairing between
the $\kk$-vector spaces $(F\otimes\Omega_X(mx))_{nx}=(F\otimes\Omega_X(mx))/(F\otimes\Omega_X((m-n)x))$ and $F^\vee((n-m)x)_{nx}=F^\vee((n-m)x)/F^\vee(-mx)$.

In particular, let $F=\END(E)=E^\vee\otimes E$ be the sheaf of endomorphisms of a vector bundle $E$. Then~$F$ is self-dual, so we have a \emph{residue-trace pairing}
\begin{equation}\label{eq:ResTr}
  \res_x\tr\colon\bigl(\END(E)\otimes\Omega_X(mx)\bigr)_{nx}\otimes_\kk
  \bigl(\END(E)((n-m)x)\bigr)_{nx}\to\kk.
\end{equation}

The following simple lemma will be used many times.
\begin{lemma}\label{lm:nilpotent}
Let $\cL$ be a line bundle on $X$ of negative degree and $E$ be a vector bundle on $X$. Then every morphism $\Psi\in\Hom(E,E\otimes\cL)$ is globally nilpotent, that is, the composition
\[
    E\xrightarrow{\Psi}E\otimes\cL\xrightarrow{\Psi\otimes\Id_\cL}E\otimes\cL^{\otimes2}
    \xrightarrow{\Psi\otimes\Id_{\cL^{\otimes2}}}\ldots
    \xrightarrow{\Psi\otimes\Id_{\cL^{\otimes\rk\cE-1}}}\Psi\otimes\cL^{\otimes\rk E}
\]
is zero.
\end{lemma}
\begin{proof}
The coefficients of the characteristic polynomial of $\Psi$ belong to $H^0(X,\cL^{\otimes i})$, where $0\le i\le\rk E$. Since $H^0(X,\cL^{\otimes i})=0$ for $i>0$, the statement follows.
\end{proof}

\subsection{$\epsilon$-connections with partially fixed singularities} Recall that we have fixed effective divisors $\divisor'=\sum_{x\in D}n'_xx$ and $\divisor=\sum_{x\in D}n_xx$ such that $\divisor'\le\divisor$. Let $E$ be a vector bundle on $X$. Let $\nabla\colon E\to E\otimes\Omega_X(\divisor)$ be a singular $\epsilon$-connection. Then near $x\in D$ we have $\nabla=\epsilon d+\rho_x+\ldots$, where $\rho_x\in(\END(E)\otimes\Omega_X(n_x))_{n'_xx}$ are leading terms, and the omitted terms have poles of order at most $n_x-n'_x$. We will be interested in the existence of $\epsilon$-connections with prescribed $\rho_x$ for all $x\in D$.

Consider the two-sided ideal of the ring $\End(E)$:
\[
    \End_{\divisor-\divisor'}(E):=H^0(X,\END(E)\otimes\cO_X(\divisor'-\divisor))=\Hom(E,E\otimes\cO_X(\divisor'-\divisor)).
\]
In other words, this is just the subspace of endomorphisms vanishing at each $x\in D$ to order at least $n_x-n'_x$. The restriction of elements of this ideal to the $n'_x$-th infinitesimal neighborhood of $x\in D$ is a homomorphism of vector spaces
\begin{equation}\label{eq:evx}
    \ev_x\colon\End_{\divisor-\divisor'}(E)\to\END(E)((n'_x-n_x)x)_{n'_xx}.
\end{equation}
Recall that for a vector bundle $E$ on $X$ we have its Atiyah class $b(E)\in H^1(X,\END(E)\otimes\Omega_X)$. We denote by $\langle\bullet,\bullet\rangle$ the Serre pairing $H^1(X,\END(E)\otimes\Omega_X)\otimes\End(E)\to\kk$. We have $\langle b(E),\Id_E\rangle=-\deg E$. We refer the reader to~\cite[Sect.~4 and Prop.~18(i)]{AtiyahConnections} and~\cite[Sect.~4.1, (4.1), and Rem.~4.2(a)]{ArinkinFedorov}.

Define $\delta_{\divisor,\divisor'}=0$ if $\divisor'<\divisor$ and $\delta_{\divisor,\divisor'}=1$ if $\divisor'=\divisor$. Recall the residue-trace pairing~\eqref{eq:ResTr}. We are ready to give the criterion of existence for $\epsilon$-connections with partially fixed polar parts.

\begin{proposition}\label{pr:existWithLeading}
Let $E$ be a vector bundle on $X$. For $x\in D$, let $\rho_x\in(\END(E)\otimes\Omega_X(n_xx))_{n'_xx}$ be a partial polar part of an $\epsilon$-connection. Then the following are equivalent: \\
\stepzero\noindstep there is an $\epsilon$-connection $\nabla\colon E\to E\otimes\Omega_X(\divisor)$ such that for all $x\in D$ we have near $x$
\[
    \nabla=\epsilon d+\rho_x+\text{terms with poles of order at most $n_x-n'_x$};
\]
\noindstep for all $\Psi\in\End_{\divisor-\divisor'}(E)$ we have
\begin{equation}\label{eq:Obst}
    \epsilon\delta_{\divisor,\divisor'}\langle b(E),\Psi\rangle-\sum_{x\in D}\res_x\tr(\rho_x\ev_x\Psi)=0,
\end{equation}
  where $b(E)$ is the Atiyah class of $E$.
\end{proposition}
\begin{proof}
  The proof is similar to that of~\cite[Prop.~5.2]{FedorovSoibelmansParabolic}. In more detail, the obstruction to the existence of the $\epsilon$-connection as above lies in $H^1(X,\END(E)\otimes\Omega_X(\divisor-\divisor'))$. As in the proof of~\cite[Prop.~5.2]{FedorovSoibelmansParabolic}, the collection $(\rho_x\,|\;x\in D)$ gives a class $b(\rho)\in H^1(X,\END(E)\otimes\Omega_X(\divisor-\divisor'))$, and the obstruction to existence of the desired $\epsilon$-connection is $\epsilon b(E)-b(\rho)$.

  On the other hand, $H^1(X,\END(E)\otimes\Omega_X(\divisor-\divisor'))$ is Serre dual to $\End_{\divisor-\divisor'}(E)$, so the connection exists if and only if $\langle\epsilon b(E)-b(\rho),\Psi\rangle=0$ for all $\Psi\in\End_{\divisor-\divisor'}(E)$. One checks that the pairing is given by the left hand side of~\eqref{eq:Obst}, using the fact that if $\divisor'<\divisor$, then by Lemma~\ref{lm:nilpotent} $\Psi$ is necessarily nilpotent, so $\langle b(E),\Psi\rangle=0$ by~\cite[Prop.~18(ii)]{AtiyahConnections}.
\end{proof}

\subsection{Existence of parabolic $\epsilon$-connections with a partially or fully fixed formal normal forms and a given underlying parabolic bundle} Recall from Definition~\ref{def:thick} the notion of parabolic bundles of level $\divisor'$. Set $\cL:=\cO_X(\divisor'-\divisor)$. The free rank one $k[n'_xx]$-module $\cL_{n'_xx}=\cO_X((n'_x-n_x)x)_{n'_xx}$, as well as its $\kk$-dual $\Omega_X(n_xx)_{n'_xx}$, will play an important role in the coordinate-free formulations of our results. Recall from~\eqref{eq:evx} the homomorphism
\[
    \ev_x\colon\End_{\divisor-\divisor'}(E)\to
    \END(E)((n'_x-n_x)x)_{n'_xx}=\Hom_{k[n'_xx]\mathrm{-mod}}(E_{n'_xx},E_{n'_xx}\otimes_{k[n'_xx]}\cL_{n'_xx}).
\]
\begin{definition}
Let $\bE=(E,E_{\bullet,\bullet})$ be a level $\divisor'$ parabolic bundle on $X$. We define $\End_{\divisor-\divisor'}(\bE)\subset\End_{\divisor-\divisor'}(E)$ to be the set of endomorphisms $\Psi$ compatible with the parabolic structure in the sense that for all $x\in D$ and $j>0$ we have
\begin{equation}\label{eq:EndPar}
    (\ev_x\Psi)(E_{x,j})\subset E_{x,j}\otimes_{k[n'_xx]}\cL_{n'_xx}.
\end{equation}
\end{definition}
We note that $\End_{\divisor-\divisor'}(\bE)$ is a 2-sided ideal of $\End(\bE)$. For an endomorphism $\Psi\in\End_{\divisor-\divisor'}(\bE)$, $x\in D$, and $j>0$ let us consider the corresponding homomorphism of $k[n'_xx]$-modules $\gr_j\ev_x(\Psi): E_{x,j-1}/E_{x,j}\to(E_{x,j-1}/E_{x,j})\otimes_{k[n'_xx]}\cL_{n'_xx}$. The trace of this homomorphism is an element of the $k[n'_xx]$-module $\cL_{n'_xx}$. We get a linear map of $\kk$-vector spaces:
\[
    \tr\gr_j\ev_x\colon\End_{\divisor-\divisor'}(\bE)
    \xrightarrow{\ev_x\gr_j}\Hom(E_{x,j-1}/E_{x,j},(E_{x,j-1}/E_{x,j})\otimes_{k[n'_xx]}\cL_{n'_xx})
    \xrightarrow{\tr}\cL_{n'_xx}.
\]
Recall from Definition~\ref{def:ModSpacePartial} the stacks $\Conn^{prtl}(\divisor,\divisor',\zeta)$. We are ready for the main result of this section. It generalizes the existence criteria from~\cite{AtiyahConnections,MihaiMonodromie,MihaiConnexions,ArinkinFedorov,FedorovSoibelmansParabolic}.
\begin{theorem}\label{th:Existence}
Let $\divisor'=\sum_{x\in D}n'_xx$ be a divisor such that for all $x\in D$ we have $0\le n'_x\le n_x$. Let $\bE=(E,E_{\bullet,\bullet})\in\Bun^{par}(\divisor')(\kk)$ be a level $\divisor'$ parabolic bundle on $X$. Let $\zeta\in\FNF(\divisor,\divisor')$ be a~partial normal form. Then there is a~$\kk$-point of $\Conn^{prtl}(\divisor,\divisor',\zeta)$ lying over $\bE$ if and only if for every $\Psi\in\End_{\divisor-\divisor'}(\bE)$ we have
    \begin{equation}\label{eq:Obst2}
        \epsilon\delta_{\divisor,\divisor'}\langle b(E),\Psi\rangle-\sum_{x\in D}\sum_{j>0}\res_x(\zeta_{x,j}[\tr\gr_j\ev_x](\Psi))=0.
    \end{equation}
\end{theorem}
We note that $\zeta_{x,j}\in\Omega_X(n_xx)_{n'_xx}$, $[\tr\gr_j\ev_x](\Psi)\in\cL_{n'_xx}$, and there is a perfect residue pairing between these vector spaces defined in Section~\ref{sect:Infinitesimal}. Note also that the sum is finite because for $j\gg0$ we have $E_{x,j-1}=E_{x,j}$, so $[\gr_j\ev_x](\Psi)=0$.

\begin{proof} We are going to follow the proof of~\cite[Thm.~7.1]{Crawley-Boevey:Indecomposable}. For the convenience of the reader we will keep as much notation from the loc.cit.~as possible. We make a convention that the tensor product or $\Hom$ of $k[n'_xx]$-modules is by default taken over $k[n'_xx]$, while the tensor product of $\cO_X$-modules is by default taken over $\cO_X$.

Let $\bE=(E,E_{\bullet,\bullet})$ be a parabolic bundle. We have an exact sequence of vector spaces
\[
\End_{\divisor-\divisor'}(\bE)\to\End_{\divisor-\divisor'}(E)\oplus\bigoplus_{x\in D}\bigoplus_{j=1}^\infty\Hom(E_{x,j}, E_{x,j}\otimes\cL_{n'_xx})\to
\bigoplus_{x\in D}\bigoplus_{j=1}^{\infty}\Hom(E_{x,j},E_{x,j-1}\otimes\cL_{n'_xx})
\]
where the arrows are as follows
\begin{itemize}
\item the morphism $\End_{\divisor-\divisor'}(\bE)\to\End_{\divisor-\divisor'}(E)$ is the obvious embedding;
\item the morphism $\End_{\divisor-\divisor'}(\bE)\to\Hom(E_{x,j}, E_{x,j}\otimes\cL_{n'_xx})$ is induced by $\ev_x$ (see~\eqref{eq:EndPar});
\item the morphism $\End_{\divisor-\divisor'}(E)\to\Hom(E_{x,j},E_{x,j-1}\otimes\cL_{n'_xx})$ is non-zero only when $j=1$
in which case it is induced by $\ev_x$.
\item the morphism $\Hom(E_{x,j},E_{x,j}\otimes\cL_{n'_xx})\to\Hom(E_{x,j},E_{x,j-1}\otimes\cL_{n'_xx})$ is $-e_{x,j}$, where $e_{x,j}: E_{x,j}\to E_{x,j-1}$ is the natural embedding, while the morphism $\Hom(E_{x,j},E_{x,j}\otimes\cL_{n'_xx})\to\Hom(E_{x,j+1},E_{x,j}\otimes\cL_{n'_xx})$ is the restriction. The other components are zero.
\end{itemize}
Applying $\Hom_\kk(\bullet,\kk)$ and using residue-trace pairing~\eqref{eq:ResTr}, we get the following exact sequence similar to~\cite[Lemma~7.4]{Crawley-Boevey:Indecomposable}:
\begin{multline}\label{eq:ExactSeq}
    \bigoplus_{x\in D}\bigoplus_{j=1}^{\infty}\Hom(E_{x,j-1},E_{x,j}\otimes(\Omega_X(n_xx))_{n'_xx})\xrightarrow{F}\\
    \End_{\divisor-\divisor'}(E)^{\vee}\oplus\bigoplus_{x\in D}\bigoplus_{j=1}^{\infty}\Hom(E_{x,j}, E_{x,j}\otimes\Omega_X(n_xx)_{n'_xx})\xrightarrow{G}
    \End_{\divisor-\divisor'}(\bE)^{\vee}.
\end{multline}
A short calculation shows that $F$ sends $\phi=(\phi_{x,j})$ with $\phi_{x,j}\in\Hom(E_{x,j-1},E_{x,j}\otimes\Omega_X(n_xx)_{n'_xx})$ to $(\xi,\eta_{x,j})$, where

\begin{equation}\label{eq:FinExact}
    \xi(\Psi)=\sum_{x\in D}\res_x\tr(\phi_{x,1}\ev_x(\Psi)),\quad
    \eta_{x,j}=\phi_{x,j+1}-\phi_{x,j}|_{E_{x,j}}.
\end{equation}

Assume now that there exists $(\bE,\nabla)\in\Conn^{prtl}(\divisor,\divisor',\zeta)(\kk)$ lying over $\bE$. We need to show that for all $\Psi\in\End_{\divisor-\divisor'}(\bE)$ we have~\eqref{eq:Obst2}. Write near $x\in D$: $\nabla=\epsilon d+\rho_x+\ldots$ so that $\rho_x\in(\END(E)\otimes\Omega_X(n_x))_{n'_xx}$ are as in Proposition~\ref{pr:existWithLeading}. The compatibility condition becomes: for $j>0$ we have $(\rho_x-\zeta_{x,j}\Id)(E_{x,j-1})\subset E_{x,j}\otimes\Omega_X(n_xx)_{n'_x x}$.

Thus, the collection $\phi_{x,j}:=\rho_x|_{E_{x,j-1}}-\zeta_{x,j}\Id$ is an element of the first term of the exact sequence~\eqref{eq:ExactSeq}. We see from~\eqref{eq:FinExact} that $F(\phi_{x,j})=(\xi,(\zeta_{x,j}-\zeta_{x,j+1})\Id)$,
where
\begin{equation}\label{eq:xi0}
    \xi(\Psi)=\sum_{x\in D}\res_x\tr((\rho_x-\zeta_{x,1}\Id)\ev_x(\Psi)).
\end{equation}
By Proposition~\ref{pr:existWithLeading} we have~\eqref{eq:Obst}, so the previous formula can be re-written as
\begin{equation}\label{eq:xi}
    \xi(\Psi)=\epsilon\delta_{\divisor,\divisor'}\langle b(E),\Psi\rangle-\sum_{x\in D}\res_x\tr(\zeta_{x,1}\ev_x(\Psi))= \epsilon\delta_{\divisor,\divisor'}\langle b(E),\Psi\rangle-\sum_{x\in D}\res_x(\zeta_{x,1}\tr\ev_x(\Psi)).
\end{equation}
On the other hand, exactness of~\eqref{eq:ExactSeq} shows that $[G(\xi,(\zeta_{x,j}-\zeta_{x,j+1})\Id)](\Psi)=0$, and a short calculation shows that this together with~\eqref{eq:xi} imply~\eqref{eq:Obst2}.

Conversely, assume that~\eqref{eq:Obst2} is satisfied for all $\Psi\in\End_{\divisor-\divisor'}(\bE)$, and define $\xi\in\End_{\divisor-\divisor'}(E)^{\vee}$ by~\eqref{eq:xi}. Then~\eqref{eq:Obst2} gives $G(\xi,(\zeta_{x,j}-\zeta_{x,j+1})\Id)=0$. The exactness of~\eqref{eq:ExactSeq} shows that there is a collection $\phi_{x,j}\in\Hom(E_{x,j-1},E_{x,j}\otimes\Omega_X(n_xx)_{n'_xx})$ such that $F(\phi_{x,j})=(\xi,(\zeta_{x,j}-\zeta_{x,j+1})\Id)$. Set $\rho_x:=\phi_{x,1}+\zeta_{x,1}\Id$. According to~\eqref{eq:FinExact}, we have~\eqref{eq:xi0}. Subtracting it from~\eqref{eq:xi} (which holds by definition of $\xi$), we see that $\rho_x$ satisfy~\eqref{eq:Obst} for all $\Psi\in\End_{\divisor-\divisor'}(E)$. Thus, Proposition~\ref{pr:existWithLeading} implies that there is an $\epsilon$-connection $\nabla\colon E\to E\otimes\Omega_X(\divisor)$ with leading term $\rho_x$ at each $x\in D$. One checks by induction on $j$ using~\eqref{eq:FinExact} that for all $x$ and $j$ we have $\phi_{x,j}=\rho_x|_{E_{x,j-1}}-\zeta_{x,j}\Id$, so that $(\rho_x-\zeta_{x,j}\Id)(E_{x,j-1})\subset E_{x,j}\otimes\Omega_X(n_xx)_{n'_x x}$. It follows that $\nabla$ is compatible with $\bE$ in the sense of Definition~\ref{def:ModSpacePartial}. The proof of Theorem~\ref{th:Existence} is complete.
\end{proof}

\subsection{Slopes and degrees of parabolic bundles}
Consider $\epsilon\in\kk$, $\zeta=(\zeta_{x,j})\in\FNF(\divisor)$, and $\gamma=(r,r_{\bullet,\bullet})\in\Gamma_D$. Recall that in Section~\ref{sect:ExplAnswers} we defined
\begin{equation}\label{eq:gamma_star_zeta}
    \gamma\star\zeta:=\sum_{x\in D}\sum_{j>0}r_{x,j}\res\zeta_{x,j}\in\kk.
\end{equation}

\begin{definition}\label{def:DegZeta}
For a parabolic bundle $\bE\in\Bun_{\gamma,d}^{par}(\divisor)(\kk)$ we define the \emph{$(\epsilon,\zeta)$-degree of $\bE$} as $\epsilon d+\gamma\star\zeta$. We define the \emph{$(\epsilon,\zeta)$-slope of $\bE$} as the quotient of its $(\epsilon,\zeta)$-degree by $\rk\bE=\rk\gamma$. We say that $\bE$ is \emph{$(\epsilon,\zeta)$-isoslopy\/} if the $(\epsilon,\zeta)$-slope of any of its non-zero direct summands is equal to that of $\bE$.
\end{definition}

\begin{corollary}\label{cor:indecomposable summand}
Let $\bE\in\Bun^{par}(\divisor)(\kk)$ be a level $\divisor$ parabolic bundle. Let $\zeta\in\FNF(\divisor)$ and assume that there is a point of $\Conn(\divisor,\zeta)$ lying over $\bE$. Then every indecomposable direct summand of $\bE$ has zero $(\epsilon,\zeta)$-degree. In other words, $\bE$ is $(\epsilon,\zeta)$-isoslopy of slope zero. In particular, $\bE$ has zero $(\epsilon,\zeta)$-degree.
\end{corollary}
\begin{proof} Take $\Psi\in\End_{\divisor-\divisor'}(\bE)=\End(\bE)$ to be the projector to its direct summand $\bE'$. Then~\eqref{eq:Obst2} yields the result (cf.~\cite[Thm.~7.1]{Crawley-Boevey:Indecomposable}).
\end{proof}

\begin{example}\label{ex:IrregNotSuff}
The converse to Corollary~\ref{cor:indecomposable summand} is false. That is, contrary to the case of regular singular parabolic connections, the vanishing of parabolic degrees of indecomposable direct summands \emph{does not imply\/} the existence of a parabolic connection. Let us give an example. Let $X=\P^1_k$ and $z$ be the standard coordinate. Let $\divisor=2(0)$, $E=\cO\oplus\cO$ be a rank~2 trivial vector bundle. Let the submodule $E_1\subset E_{2(0)}$ be generated by $e:=(1,z)$. Define a level~2 parabolic structure on $E$ at~0 as follows: $E_{0,0}=E_{2(0)}$, $E_{0,1}=E_1$, $E_{0,j}=0$ for $j\ge2$. It is easy to see that $\bE:=(E,E_{0,\bullet})$ is indecomposable.

Let us consider a nilpotent endomorphism of $E$ given by
\[
    \Psi=\begin{pmatrix}
           0 & 1 \\
           0 & 0
         \end{pmatrix}.
\]
Then $\Psi|_{2(0)}e=(z,0)=z(1,z)$, so that $E_1$ is preserved by $\Psi$ and $\Psi$ acts on $E_1$ via multiplication by~$z$. Since $\tr\Psi=0$, $\Psi$ acts on $E_{2(0)}/E_1$ via multiplication by $-z$. It follows from Theorem~\ref{th:Existence} that there is no connection $\nabla$ on $E$ with poles bounded by $2(0)$ preserving $E_{0,\bullet}$ and such that
\[
    \nabla|_{E_1}=d+\frac{dz}{z^2}+0\frac{dz}z+\ldots, \quad
    \nabla|_{E_{2(0)}/E_1}=d+0\frac{dz}{z^2}+0\frac{dz}z+\ldots.
\]
We will have, however, a partial converse to Corollary~\ref{cor:indecomposable summand} in Section~\ref{sect:ExistFullyFixed}.
\end{example}

We just saw that even when $\Psi$ is nilpotent, $\tr\gr_j\ev_x\Psi\ne0$ in general. We have a weaker statement that will be useful in Sections~\ref{sect:ExLeading} and~\ref{sect:ExInd}. Recall the $k[n'_xx]$-module $\cL_{n'_xx}$. Consider its submodule $\cL_{n'_xx}(-x)$ consisting of sections vanishing at $x$. This submodule is not free unless $n'_x=1$, in which case $\cL_{n'_xx}(-x)=0$.

\begin{lemma}\label{lm:nilpotent2}
Let $\Psi\in\End_{\divisor-\divisor'}(\bE)$ be globally nilpotent. Then for all $x\in D$ such that $n'_x>0$ and all $j>0$, the leading term of $[\tr\gr_j\ev_x]\Psi$ vanishes in the sense that $[\tr\gr_j\ev_x]\Psi\in\cL_{n'_xx}(-x)$.
\end{lemma}
We note that if $\divisor'<\divisor$, then $\Psi$ is globally nilpotent by Lemma~\ref{lm:nilpotent}.
\begin{proof}
First of all, we generalize the statement by allowing the curve $X$ to be quasi-projective. The statement is local over $X$, so, trivializing $\cL$ in an open affine neighborhood of $x$ and replacing $X$ with this neighborhood, we may assume that $\cL=\cO_X$. In this case $\Psi$ is just a nilpotent endomorphism of $E$. For a $k[n'_xx]$-module $M$, let $M_x$ denote its closed fiber. Then $(E_{x,j-1}/E_{x,j})_x$ is a subquotient of $E_x$. Since $\Psi_x$ is a nilpotent endomorphism of $E_x$, it induces a nilpotent endomorphism $\gr_j\Psi_x$ of $(E_{x,j-1}/E_{x,j})_x$. Thus, $\tr\gr_j\Psi_x=0$, which is equivalent to our statement.
\end{proof}

Recall from Definition~\ref{def:DegZeta} the notion of isoslopy vector bundles. The isoslopy vector bundles played an important role in calculating motivic classes of connections, see~\cite[Sect.~3.2.1]{FedorovSoibelmans}, of parabolic Higgs bundles, see~\cite[Sect.~5]{FedorovSoibelmansParabolic}, and of parabolic bundles with connections, see~\cite[Sect.~8]{FedorovSoibelmansParabolic}. We will also use them extensively in this paper, see Sections~\ref{sect:ExistFullyFixed} and~\ref{sect:ParPairs}. It is important to know that this notion is stable under field extensions.
\begin{lemma}\label{lm:IsoslopyFieldExt}
  Assume that $K\supset\kk$ is a field extension. Then $\bE\in\Bun^{par}(\divisor)(\kk)$ is $(\epsilon,\zeta)$-isoslopy if and only if its pullback to $X_K$ is $(\epsilon,\zeta)$-isoslopy.
\end{lemma}
\begin{proof}
 The proof is completely similar to that of~\cite[Lemma~3.2.2]{FedorovSoibelmans}, see also~\cite[Sect.~5.4]{FedorovSoibelmansParabolic}. Let us sketch it. The `if' direction is obvious. Since indecomposable direct summands of a parabolic bundle are unique up to an isomorphism, we can reduce the statement to showing that if $\bE$ is indecomposable, then $\bE_K$ is $(\epsilon,\zeta)$-isoslopy. The `if' direction shows that it is enough to consider the two cases: when~$K$ is an algebraic closure of $\kk$, and when $K=\kk(t)$, where $\kk$ is algebraically closed. In the first case we note that the Galois group of $\kk$ acts transitively on the direct summand of $\bE_K$. In the second case, we note that if $\bE_K=\bE'\oplus\bE''$, then we can extend $\bE$, $\bE'$, and $\bE''$ to $X\times_\kk U$, where $U$ is an open subscheme of $\A^1_\kk$. Now, restricting to $X\times_\kk u$, where $u\simeq\Spec\kk$ is a closed point of $U$, we see that $\bE'$ and $\bE''$ have the same slope. (We have used that $\kk$ has characteristic zero so that the algebraic closure coincides with the separable closure. The statement is still valid when $\kk$ has finite characteristic, but one has to consider the case of a purely inseparable extension as well.)
\end{proof}

\subsection{Existence of parabolic $\epsilon$-connections with prescribed leading terms}\label{sect:ExLeading}
\begin{lemma}\label{lm:NoAtiyah}
Assume that $\divisor'<\divisor$ is a reduced divisor, that is, there is a subset $D'\subset D$ such that $\divisor'=\sum_{x\in D'}x$. Then for any $\epsilon\in\kk$, any $\zeta\in\FNF(\divisor,\divisor')$, and any parabolic bundle $\bE\in\Bun^{par}(\divisor')(\kk)$, the fiber of $\Conn^{prtl}(\divisor,\divisor',\zeta)$ over $\bE$ is non-empty.
\end{lemma}
\begin{proof}
    Let $\Psi\in\End_{\divisor-\divisor'}(\bE)$ be any endomorphism. By Lemma~\ref{lm:nilpotent}, it is nilpotent. For $x\in D'$ we have $n'_x=1$, so by Lemma~\ref{lm:nilpotent2} we have $\tr\gr_j\ev_x\Psi\in\cL_{x}(-x)=0$ for any $j>0$. Since $\delta_{\divisor,\divisor'}=0$ as well, all the terms in the obstruction~\eqref{eq:Obst2} vanish. It remains to use Theorem~\ref{th:Existence}.
\end{proof}

\subsection{Inductive existence criteria}\label{sect:ExInd}
Let $x\in D$ be such that $n'_x>0$. We have the natural projection homomorphism
\begin{equation}\label{eq:Px}
    P_x\colon\FNF(\divisor,\divisor')\to\FNF(\divisor,\divisor'-x).
\end{equation}
This homomorphism will play an important role in our inductive calculation of motivic classes of parabolic $\epsilon$-connections (see Section~\ref{sect:IdeaOfCalc}, Lemma~\ref{lm:AffineBundle2}, and Lemma~\ref{lm:RR}).

\begin{lemma}\label{lm:Independent}
  Let $x\in D$ be such that $n'_x>0$. Let $\zeta,\theta\in\FNF(\divisor,\divisor')$ be such that $P_x(\zeta)=P_x(\theta)$. Then there is a~$\kk$-point of $\Conn^{prtl}(\divisor,\divisor',\zeta)$ lying over $\bE\in\Bun^{par}(\divisor')(\kk)$ if and only if there is a~$\kk$-point of $\Conn^{prtl}(\divisor,\divisor',\theta)$ lying over $\bE$.
\end{lemma}
\begin{proof}
    Put $\bE=(E,E_{\bullet,\bullet})$. Take $\Psi\in\End_{\divisor-\divisor'}(\bE)$. In view of Theorem~\ref{th:Existence}, we only need to show that
    \[
        \sum_{y\in D}\sum_{j>0}\res_y(\zeta_{y,j}[\tr\gr_j\ev_y](\Psi))-
        \sum_{y\in D}\sum_{j>0}\res_y(\theta_{y,j}[\tr\gr_j\ev_y](\Psi))=0.
    \]

    Since $P_x(\zeta)=P_x(\theta)$, we have $\zeta_{y,j}=\theta_{y,j}$ if $x\ne y$. Thus, we just need to show that
    \[
        \sum_{j>0}\res_x((\zeta_{x,j}-\theta_{x,j})[\tr\gr_j\ev_x](\Psi))=0.
    \]
    Since $P_x(\zeta)=P_x(\theta)$, we see that for all $j>0$
    \[
        \zeta_{x,j}-\theta_{x,j}\in\Omega_X((n_x-n'_x+1)x)/\Omega_X((n_x-n'_x)x)\subset\Omega_X(n_xx)_{n'_xx}.
    \]
    Combining Lemmas~\ref{lm:nilpotent} and~\ref{lm:nilpotent2}, we see that $[\tr\gr_j\ev_x]\Psi\in\cL_{n'_xx}(-x)$. The subspace
    \[
        \Omega_X((n_x-n'_x+1)x)/\Omega_X((n_x-n'_x)x)\subset\Omega_X(n_xx)/\Omega_X((n_x-n'_x)x)
    \]
        is the annihilator of $\cL_{n'_xx}(-x)$ under the residue pairing, thus,
    \[
        \res_x((\zeta_{x,j}-\theta_{x,j})[\tr\gr_j\ev_x](\Psi))=0
    \]
    and the desired equality holds.
\end{proof}

\subsection{Existence of parabolic $\epsilon$-connections with fully fixed formal normal forms}\label{sect:ExistFullyFixed} We prove a partial converse of Corollary~\ref{cor:indecomposable summand}.

\begin{lemma}\label{lm:existIrregFull}
Let $x\in D$, $\bE\in\Bun^{par}(\divisor)(\kk)$ be a parabolic bundle, and let $\zeta\in\FNF(\divisor)$. Assume that there is $\theta\in\FNF(\divisor)$ such that $P_x(\theta)=P_x(\zeta)$ and the $\bE$-fiber of the projection $\Conn(\divisor,\theta)\to\Bun^{par}(\divisor)$ is non-empty. Then the $\bE$-fiber of the projection $\Conn(\divisor,\zeta)\to\Bun^{par}(\divisor)$ is non-empty if and only if every direct summand of $\bE$ has $(\epsilon,\zeta)$-degree zero.
\end{lemma}
\begin{proof}
The `only if' part follows from Corollary~\ref{cor:indecomposable summand}. For the `if' part, note that in view of Lemma~\ref{lm:IsoslopyFieldExt} we may assume that $\kk$ is algebraically closed. Write $\bE=\bigoplus_i\bE_i$, where $\bE_i$ are indecomposable parabolic bundles (necessarily of $(\epsilon,\zeta)$-degree zero). Let $\Psi\in\End(\bE)$ be an endomorphism. We need to verify condition~\eqref{eq:Obst2} of Theorem~\ref{th:Existence}. We have
\[
    \End(\bE)=\bigoplus_{i,j}\Hom(\bE_i,\bE_j).
\]
Since~\eqref{eq:Obst2} is linear in $\Psi$, we may assume that $\Psi\in\Hom(\bE_i,\bE_j)$ for some $i$ and $j$.

\emph{Case 1: $i\ne j$.} In this case $\Psi$ is globally nilpotent, so the proof repeats that of Lemma~\ref{lm:Independent}.

\emph{Case 2: $i=j$.} In this case $\Psi\in\End(\bE_i)$. Since $\kk$ is algebraically closed, $\End(\bE_i)\subset\End(\bE)$ is linearly generated by its nilpotent radical and the projector to $\bE_i$ (indeed, the quotient of $\End(\bE_i)$ by its nilpotent radical is a division algebra, which must be equal to $\kk$). Thus, it is enough to consider two subcases: $\Psi$ is the projector to $\bE_i$ and $\Psi$ is nilpotent. In the second subcase, we argue as in Case~1 and Lemma~\ref{lm:Independent}. If $\Psi$ is the projector to $\bE_i$, then~\eqref{eq:Obst2} becomes the condition that the $(\epsilon,\zeta)$-degree of $\bE_i$ is zero (cf.~the proof of Corollary~\ref{cor:indecomposable summand}).
\end{proof}

\section{Parabolic bundles with twisted endomorphisms}\label{Sect:parabolic twisted}
In this section we start calculating motivic classes. Let, as before, $X$ be a smooth projective geometrically connected $\kk$-curve. We will calculate the motivic classes of moduli stacks of parabolic bundles with twisted endomorphisms. We only consider parabolic bundles of level $D$, where $D$ is a~reduced divisor. We do not know how to solve a similar problem for higher level parabolic bundles with endomorphisms (see Section~\ref{sect:FutureDirections}). The main result of this section is Corollary~\ref{cor:MotMellitPunctures}. The results are obtained by combining techniques from~\cite{SchiffmannIndecomposable}, \cite{MozgovoySchiffmann2020}, \cite{FedorovSoibelmans}, and~\cite{FedorovSoibelmansParabolic}.

\subsection{Recollection on motivic classes} Recall from Section~\ref{sect:MotClasses} and~\cite[Sect.~1.1]{FedorovSoibelmans} the rings of motivic classes $\Mot(\kk)$ and $\cMot(\kk)$. We also defined in~\cite[Sect.~2.4]{FedorovSoibelmans} the notion of a constructible subset of a stack. If $\cS$ is a constructible subset of a stack $\cX$ of finite type over $\kk$, then we have its motivic class $[\cS]\in\Mot(\kk)$. We denote its image in $\cMot(\kk)$ by $[\cS]$ as well. In the following proposition we collect the facts about motivic classes that we will use in this paper.

\begin{proposition}\label{pr:MotivicClasses}
Let $\cS$ be a constructible subset of finite type of a stack $\cX$ locally of finite type over~$\kk$.

(i) Let $f_i\colon\cY_i\to\cX$, $i=1,2$ be finite type morphisms of $\kk$-stacks and $m$ be an integer. If for every field extension $K\supset\kk$ and every point $\xi\in\cS(K)$ we have $[f_1^{-1}(\xi)]=\bL_K^m[f_2^{-1}(\xi)]$ in $\Mot(K)$ (here $\bL_K:=[\A^1_K]$), then $[f_1^{-1}(\cS)]=\bL^m[f_2^{-1}(\cS)]$ in $\Mot(\kk)$.

(ii) Let $f\colon\cY\to\cX$ be a finite type morphism of $\kk$-stacks and $\cZ$ be a $\kk$-stack of finite type. If for every field extension $K\supset\kk$ and every point $\xi\in\cS(K)$ we have $[f^{-1}(\xi)]=[\cZ\times_\kk K]$, then $[f^{-1}(\cS)]=[\cS][\cZ]$ in $\Mot(\kk)$.
\end{proposition}
\begin{proof}
    (i) We will use the language of motivic functions from~\cite[Sect.~2.1]{FedorovSoibelmans}. Recall the ring $\Mot^{fin}(\cX)$ of motivic functions with finite support. The constructible subsets $f_i^{-1}(\cS)\subset\cY$ are of finite type. Consider the motivic functions $\phi_i:=[f_i^{-1}(\cS)\to\cX]\in\Mot^{fin}(\cX)$, see~\cite[Sect.~2.4]{FedorovSoibelmans}. For every field extension $K\supset\kk$, and every point $\xi\colon\Spec K\to\cX$, we have $\xi^*\phi_1=\bL_K^m\xi^*\phi_2$. Indeed, if $\xi\in\cS(K)$, this follows from the assumption. Otherwise, both sides are equal to zero. Thus, by~\cite[Prop.~2.6.1]{FedorovSoibelmans}, we have $\phi_1=\bL^m\phi_2$. We have the pushforward $pt_*\colon\Mot^{fin}(\cX)\to\Mot(\kk)$, where $pt\colon\cX\to\Spec\kk$ is the structure morphism. Since $[f_i^{-1}(\cS)]=pt_*\phi_i$ and $\phi_1=\bL^m\phi_2$, we are done.

    (ii) Apply part (i) with $m=0$ to $f$ and the projection $\cZ\times\cX\to\cX$.
\end{proof}
We often apply this proposition with $\cS=\cX$. We also note that it is not enough to check the equality of the fibers over $\kk$-points. Indeed, assume that $\kk$ is algebraically closed and consider the morphism $\A^1_\kk-0\to\A^1_\kk-0$ sending $z$ to $z^2$. Every $\kk$-fiber of this morphism is isomorphic to $\Spec\kk\sqcup\Spec\kk$, but $[\A^1_\kk-0]\ne2[\A^1_\kk-0]$.

\subsection{HN-nonpositive and HN-nonnegative vector bundles}\label{sect:nonpositive}
Recall that a vector bundle $E$ on the curve $X$ is called \emph{nonnegative\/} (resp.~\emph{nonpositive}) if it has no quotient bundles of negative degree (resp.~no subbundles of positive degree), see~\cite[Sect.~3.2]{FedorovSoibelmans} and~\cite[Sect.~3.2]{FedorovSoibelmansParabolic} . Equivalently, its Harder--Narasimhan spectrum is nonnegative (resp.~nonpositive). Therefore the corresponding categories coincide with the categories $\cC(V)$ from~\cite[Sect.~3.4 and~6.1]{KontsevichSoibelman08} where $V$ is first (resp. third) quadrant of the plane and $\cC=\Bun(X)$. Recall from~\cite[Lemma~3.2.1]{FedorovSoibelmans} that there is an open substack $\Bun^+_{r,d}(X)\subset\Bun_{r,d}(X)$ classifying nonnegative vector bundles. Moreover,
$\Bun^+_{r,d}(X)$ is of finite type over $\kk$. It follows from the duality for vector bundles that there is also an open substack $\Bun^-_{r,d}(X)\subset\Bun_{r,d}(X)$ classifying nonpositive vector bundles and that it is also of finite type over $\kk$ (see~\cite[Lemma~3.2]{FedorovSoibelmansParabolic}).

If $\cX$ is any stack over $\Bun(X)$, we will set
\begin{equation}\label{eq:X-}
    \cX_{r,d}:=\cX\times_{\Bun(X)}\Bun_{r,d}(X),\quad \cX^-_{r,d}:=\cX\times_{\Bun(X)}\Bun^-_{r,d}(X),\quad\cX^+_{r,d}:=\cX\times_{\Bun(X)}\Bun^+_{r,d}(X).
\end{equation}
We see that if $\cX$ is of finite type over $\Bun(X)$, then $\cX^-_{r,d}$ and $\cX^+_{r,d}$ are of finite type over $\kk$. We have decompositions into disjoint unions of open and closed subsets
\begin{equation}\label{eq:X-decomp}
  \cX=\bigsqcup_{r,d}\cX_{r,d},\quad
  \cX^\pm=\bigsqcup_{r,d}\cX^\pm_{r,d}.
\end{equation}

\subsection{Twisted nilpotent pairs}\label{sect:TwistedPairs}
Fix a line bundle $\cL$ on $X$. Denote by $\Nil(\cL)$ the stack classifying pairs $(E,\Psi)$, where $E$ is a vector bundle on $X$ and $\Psi:E\to E\otimes\cL$ is a nilpotent morphism. We note that we are mostly interested in the case when $\deg\cL<0$, in which case nilpotence is automatic by Lemma~\ref{lm:nilpotent}. The forgetful morphism $\Nil(\cL)\to\Bun(X)$ is schematic and of finite type, so, according to Section~\ref{sect:nonpositive} we have open substacks
$\Nil_{r,d}^+(\cL)\subset\Nil(\cL)$ of finite type over $\kk$. We will now calculate their motivic classes; this is an analogue of~\cite[Thm.~1.4.1]{FedorovSoibelmans}.

As in~\cite[Sect.~4.1]{FedorovSoibelmansParabolic} we see that the stack $\Nil_{r,d}^+(\cL)$ has a locally closed stratification by the generic Jordan type of $\Psi$:
\begin{equation}\label{eq:StratMu}
    |\Nil_{r,d}^+(\cL)|=\bigcup_{\mu\,\vdash r}|\Nil_{r,d}^+(\cL,\mu)|,
\end{equation}
where the strata are parameterized by partitions $\mu$ of $r$. We recall that the notation $|\cX|$ stands for the set of equivalence classes of points of a stack $\cX$.

Recall the motivic series $J_{\mu,X}^{mot}(z), H_{\mu,X}^{mot}(z)\in\cMot(\kk)[[z]]$ defined in~\cite[Sect.~1.3.2]{FedorovSoibelmans} (see also Section~\ref{sect:Schiffmann}, where the universal versions of these series are defined).
\begin{proposition}\label{pr:NilpEnd}
Let $\cL$ be a line bundle on $X$ of degree $-\delta$, where $\delta>0$, and $\mu$ be a partition of a positive integer $r$. We have the following identity in $\cMot(\kk)[\bL^\frac12][[z]]$.
\[
    \bL^{\frac{r^2\delta}2}\sum_{d\ge0}[\Nil_{r,d}^+(\cL,\mu)]z^d=
    \bL^{\frac{(2g-2+\delta)\langle\mu,\mu\rangle}2}
    z^{\delta n(\mu')}J_{\mu,X}^{mot}(z)H_{\mu,X}^{mot}(z),
\]
where $g$ is the genus of $X$, $\mu'$ is the conjugate partition, $\langle\mu,\mu\rangle$ and $n(\mu')$ are defined in Section~\ref{sect:Macdonald}.
\end{proposition}
Note that the powers of $\bL$ are either all integral or all half-integral because $r=|\mu|$ and $\langle\mu,\mu\rangle$ have the same parity. The finite field version of this statement follows from~\cite[Thm.~6.2]{MozgovoySchiffmann2020} (where $\ell=-\delta$). More precisely,~\cite[Thm.~6.2]{MozgovoySchiffmann2020} gives the identity~\eqref{eq:NilCoh}, which implies the proposition by Lemma~\ref{lm:CohVec}. We provide some details, indicating the changes needed for the motivic case. We need a lemma, which is probably well-known, but we were not able to find a reference.
\begin{lemma}\label{lm:LocallyClosed}
  Let $S$ be a locally Noetherian $\kk$-scheme and let $F$ be a coherent sheaf on $X\times_\kk S$. Then the set of $s\in S$ such that $F_s$ has a given degree and a given generic rank is locally closed.
\end{lemma}
\begin{proof}
  First we show the statement about the rank. It is enough to show that the generic rank is upper semicontinuous, that is, if $s\in S$ and $F_s$ has generic rank $r$, then there is a Zariski open neighborhood~$U$ of~$s$ such that the generic rank of $F_{s'}$ is at most $r$ for all $s'\in U$. There is a point $x\in X_s$ such that the vector space $F_x$ is generated by $r$ vectors $f_1$, \ldots, $f_r$. We can lift $f_i$ to sections $\tilde f_i$ of $F$ defined in an open neighborhood $W\subset X\times_\kk S$ of $x$. Shrinking $W$, we may assume that they generate $F$ over $W$. Let $U$ be the image of $W$ under the projection $X\times_\kk S\to S$. Then $U$ is open and for every $s'\in U$, the sheaf $F_{s'}$ is generated by $r$ sections over the non-empty open subset $W\cap X_{s'}$. Thus the generic rank of $F_{s'}$ is at most $r$.

  Having taken care of the ranks, we may now assume that the generic rank of $F$ is constant on the fibers of $X\times_\kk S\to S$. Denote this generic rank by $r$. We will show that the degrees are upper semicontinuous. Let $L$ be a relatively very ample line bundle of relative degree $e$. Then for $s\in S$ and $N\in\Z$ we have $\deg(F_s\otimes L_s^{\otimes N})=\deg F_s+Nre$. We see that it is enough to prove the statement with $F$ replaced by $F\otimes L^{\otimes N}$. Since the statement is local over $S$, we may assume that $S$ is Noetherian. Taking~$N$ large, we may assume that for all $s\in S$ we have $H^1(X_s,F_s)=0$. Thus, by the Riemann--Roch theorem we just need to show that the dimensions $d(s):=\dim H^0(X_s,F_s)$ are upper semicontinuous. Pick $s\in S$. Twisting by a higher power of $L$, we may assume that $F_s$ is generated by $d(s)$ sections $f_1,\ldots,f_{d(s)}$ and that these sections lift to sections $\tilde f_i\in H^0(X\times_\kk U,F_U)$ for an open neighborhood $U$ of $s$. Shrinking $U$, we may assume that the sections generate $F_{s'}$ for all $s'\in U$. Since the pullback functor is right exact, these sections generate $H^0(X_{s'},F_{s'})$ for all $s'\in U$, and we see that $d(s')\le d(s)$.
\end{proof}

\begin{proof}[Proof of Proposition~\ref{pr:NilpEnd}]
\emph{Step 1.} For a coherent sheaf $F$ on $X$ we define its class
\[
    \cl(F):=(\rk F,\deg F)\in\Z_{\ge0}\times\Z.
\]
Let $\Nil^{Coh}=\Nil^{Coh}(\cL)$ be the stack classifying pairs $(F,\Psi)$, where $F$ is a coherent sheaf on $X$ and $\Psi\colon F\to F\otimes\cL$ is a nilpotent morphism. More precisely, if $S$ is a scheme, then an $S$-point of $\Nil^{Coh}$ is a pair $(F,\Psi)$, where $F$ is a coherent sheaf on $X\times S$ flat over $S$ and $\Psi\colon F\to F\otimes p_1^*\cL$ is a nilpotent morphism of coherent sheaves. Here $p_1: X\times S\to X$ is the natural projection. Since the Hilbert polynomials are locally constant for a flat family of coherent sheaves (\cite[Thm.~7.9.4]{EGAIII-2}), there is a~closed and open substack $\Nil_{r,d}^{Coh}\subset\Nil^{Coh}$ classifying pairs $(F,\Psi)$ with $\cl(F)=(r,d)$.

Similarly to~\cite[Sect.~5.1]{MozgovoySchiffmann2020} we stratify this stack as follows. Let $(F,\Psi)\in\Nil^{Coh}(K)$, where $K\supset\kk$ is a field extension. For $i\ge0$ let $F_i$ denote the image of the composition
\begin{equation}\label{eq:Image}
    F\otimes\cL_K^{-i}\xrightarrow{\Psi
    \otimes\Id_{\cL_K^{-i}}}F\otimes\cL_K^{1-i}\xrightarrow{\Psi\otimes\Id_{\cL_K^{1-i}}}\ldots
    \xrightarrow{\Psi\otimes\Id_{\cL_K^{-1}}}F.
\end{equation}
Then $F_0=F$ and for all $i\ge0$ we have $F_{i+1}\subset F_i$; as in~\cite[Sect.~5.1]{MozgovoySchiffmann2020}, we set $F_i'':=F_i/F_{i+1}$. By construction, $\Psi$ induces a chain of surjective homomorphisms
\begin{equation}\label{eq:chain}
F_0''\twoheadrightarrow F''_1\otimes\cL_K\twoheadrightarrow F''_2\otimes\cL_K^{\otimes2}\twoheadrightarrow\ldots.
\end{equation}
Note that for $i\gg0$ we have $F_i=0$ and $F''_i=0$. For $i>0$ set
\[
    \alpha_i:=\cl(F''_{i-1}\otimes\cL_K^{\otimes(i-1)})-\cl(F''_i\otimes\cL_K^{\otimes i})\in\Z_{\ge0}\times\Z.
\]
The sequence $\underline\alpha=(\alpha_i)$, where $\alpha_i=(r_i,d_i)\in\Z_{\ge0}\times\Z$ is called the Jordan type of $(F,\Psi)$. Next, for $(r,d)\in\Z_{\ge0}\times\Z$ and $n\in\Z$, let us set $(r,d)(n):=(r,r+dn)$, so that $\cl(F\otimes L)=\cl(F)(\deg L)$, whenever $L$ is a line bundle. Then $\cl(F_i)-\cl(F_{i+1})=\cl(F_i'')=\sum_{j>i}\alpha_j(i\delta)$, see~\cite[Lemma~5.1]{MozgovoySchiffmann2020}.
\begin{lemma}
  For a sequence $\underline\alpha=(\alpha_i)$ there is a locally closed reduced substack $\Nil^{Coh}(\underline\alpha)\subset\Nil^{Coh}$ classifying nilpotent pairs with Jordan type $\underline\alpha$.
\end{lemma}
\begin{proof}
   Let $(F,\Psi)$ be an $S$-point of $\Nil^{Coh}$, where $S$ is a locally Noetherian test scheme. One can define the sheaves $F_i$ similarly to the case $S=\Spec K$ above. Next, for $s\in S$ fixing the Jordan type of $(F_s,\Psi_s)$ is equivalent to fixing $\cl((F_i)_s)$, which is equivalent to fixing the classes of cokernels of homomorphisms in~\eqref{eq:Image} for all $i$. The formation of cokernels is compatible with restrictions to points because the restriction functor is right exact. Now it follows from Lemma~\ref{lm:LocallyClosed} that the set of $s\in S$ such that $(F_s,\Psi_s)$ has Jordan type $\underline\alpha$ is locally closed in $S$. Thus, the points of $\Nil^{Coh}$ corresponding to nilpotent pairs with Jordan type $\underline\alpha$ form a locally closed subset. By~\cite[Tag~0509]{StacksProject}, we can uniquely assign to this subset a reduced locally closed substack of $\Nil^{Coh}$.
\end{proof}

We return to the proof of Proposition~\ref{pr:NilpEnd}. Similarly to~\eqref{eq:StratMu} we have a stratification of $\Nil^{Coh}$ by the constructible subsets $\Nil^{Coh}_{r,d}(\mu)$ classifying $(F,\Psi)$ with $\Psi$ generically of nilpotent type $\mu$ and $\deg F=d$, where $r\ge0$, $\mu$ is a partition of $r$, and $d\in\Z$. It is easy to see that for a point of $\Nil^{Coh}(\underline\alpha)$, the generic nilpotent type is $\mu=1^{r_1}2^{r_2}\ldots$, where $\alpha_i=(r_i,d_i)$. According to~\cite[Lemma~5.2]{MozgovoySchiffmann2020}, we have $\rk F=\sum_iir_i=|\mu|$ and $\deg F=\sum_iid_i+\delta n(\mu')$. Thus, we see from the above that we have a stratification
\begin{equation}\label{eq:Stratification}
    \Nil^{Coh}_{r,d+\delta n(\mu')}(\mu)=
    \bigcup_{\underline\alpha\in J(\mu,d)}|\Nil^{Coh}(\underline\alpha)|,
\end{equation}
where $J(\mu,d)$ denotes the collection of all $\underline\alpha=((r_i,d_i))$ with $\mu=1^{r_1}2^{r_2}\ldots$, $\sum id_i=d$ (in particular, $r_i=d_i=0$ for $i\gg0$) and such that $d_i\ge0$ if $r_i=0$.

\emph{Step~2.} We define $\Filt$ to be the stack classifying chains of epimorphism of coherent sheaves on the curve $X$: $H_0\twoheadrightarrow H_1\twoheadrightarrow H_2\twoheadrightarrow\ldots$. Let $\underline\alpha=(\alpha_i)$ be as above, let $\Filt(\underline\alpha)$ denote the closed and open substack of $\Filt$ consisting of chains such that for all $i$ we have $\cl(H_i)=\sum_{j>i}\alpha_j(i\delta)$. We would like to compare the stacks $\Nil^{Coh}(\underline\alpha)$ and $\Filt(\underline\alpha)$. Let $\Nil^{Coh,+}$ denote the constructible subset of $\Nil^{Coh}$ consisting of pairs $(F,\Psi)$, where $F$ is nonnegative, that is, has no quotients of negative degree. Set $\Nil_{r,d}^{Coh,+}:=\Nil^{Coh,+}\cap\Nil_{r,d}^{Coh}$ and
$\Nil^{Coh,+}(\underline\alpha):=\Nil^{Coh,+}\cap\Nil^{Coh}(\underline\alpha)$. Similarly, let $\Filt^+(\underline\alpha)$ denote the open substack of $\Filt(\underline\alpha)$ classifying $H_0\twoheadrightarrow H_1\twoheadrightarrow H_2\twoheadrightarrow\ldots$ such that $H_0$ is nonnegative.

For $(r_1,d_1),(r_2,d_2)\in\Z_{\ge0}\times\Z$ set $\langle(r_1,d_1),(r_2,d_2)\rangle:=(1-g)r_1r_2+(r_1d_2-r_2d_1)$. For a sequence $\underline\alpha=(\alpha_i)$, $\alpha_i=(r_i,d_i)$ such that $\alpha_i=(0,0)$ for $i\gg0$ and an integer $l$, set
\[
    \dd_l(\underline\alpha):=-\sum_{k\ge0}\sum_{\substack{i\ge k\\j\ge k+1}}\langle\alpha_{i+1}(kl),\alpha_{j+1}(jl)\rangle.
\]
The following lemma is a motivic version of~\cite[Thm.~5.3]{MozgovoySchiffmann2020}.
\begin{lemma}\label{lm:Projection}
  For any $\underline\alpha\in J(\mu,d)$, we have in $\Mot(\kk)$:
\begin{equation*}
    [\Nil^{Coh,+}(\underline\alpha)]=\bL^{\dd_{\delta}(\underline\alpha)}[\Filt^+(\underline\alpha)].
\end{equation*}
\end{lemma}
\begin{proof}
We construct a morphism $\pi_{\underline\alpha}\colon\Nil^{Coh}(\underline\alpha)\to\Filt(\underline\alpha)$ by sending $(F,\theta)$ to the chain of epimorphisms~\eqref{eq:chain}. We need to explain that this works in families, that is, if $F$ is an $S$-flat sheaf on $X\times S$, then $F''_i$ are also $S$-flat. Since $\Nil^{Coh}(\underline\alpha)$ are reduced and $\kk$ has characteristic zero, it is enough to consider families parameterized by reduced schemes to the effect that we may assume that $S$ is reduced. In this case, the sheaves $F''_i$ are flat because they have constant Hilbert polynomials, see~\cite[Ch.~III, Thm.~9.9]{Hartshorne} (note that Hartshorne assumes that $S$ is integral but essentially the same proof goes through in the reduced case).

Next, similarly to~\cite[Thm.~5.3]{MozgovoySchiffmann2020} we verify that the fibers of $\pi_{\underline\alpha}$ are iterated affine bundles of total rank $\dd_\delta(\underline\alpha)$. Thus, the motivic classes of the fibers are equal to $\bL^{\dd_{\delta}(\underline\alpha)}$. As in~\cite[Prop.~5.6]{MozgovoySchiffmann2020}, one checks that $\pi_{\underline\alpha}^{-1}(\Filt^+(\underline\alpha))=\Nil^{Coh,+}(\underline\alpha)$, and the lemma follows from Proposition~\ref{pr:MotivicClasses}(i).
\end{proof}

\emph{Step 3.} We return to the proof of Proposition~\ref{pr:NilpEnd}. According to~\cite[Lemma~5.5]{MozgovoySchiffmann2020}, we have
\begin{equation}\label{eq:ddelta}
    \dd_\delta(\underline\alpha)=\dd_0(\underline\alpha)-\frac{r^2\delta}2+\frac \delta2\langle\mu(\underline\alpha),\mu(\underline\alpha)\rangle,
\end{equation}
where $\mu(\underline\alpha):=1^{r_1}2^{r_2}\ldots$ and $r=\sum_ir_i$.

Set $J(\mu):=\bigsqcup_{d\ge0}J(\mu,d)$ and for $\underline\alpha=((r_i,d_i))\in J(\mu)$ set $\deg\underline\alpha:=\sum_{i>0}id_i$. From~\eqref{eq:Stratification}, Lemma~\ref{lm:Projection}, and~\eqref{eq:ddelta}, we obtain
\begin{multline*}
    \bL^{\frac{r^2\delta}2}\sum_{d\ge0}[\Nil^{Coh,+}_{r,d}(\cL,\mu)]z^d=
    \bL^{\frac{r^2\delta}2}\sum_{e\ge0}\sum_{\underline\alpha\in J(\mu,e)}[\Nil^{Coh,+}(\underline\alpha)]z^{e+\delta n(\mu')}=\\
    \bL^{\frac{r^2\delta}2}z^{\delta n(\mu')}
    \sum_{\underline\alpha\in J(\mu)}[\Nil^{Coh,+}(\underline\alpha)]z^{\deg\underline\alpha}=
    \bL^{\frac{r^2\delta}2}z^{\delta n(\mu')}
    \sum_{\underline\alpha\in J(\mu)}\bL^{\dd_\delta(\underline\alpha)}[\Filt^+(\underline\alpha)]
    z^{\deg\underline\alpha}\\=
    \bL^{\frac\delta2\langle\mu,\mu\rangle}z^{\delta n(\mu')}
    \sum_{\underline\alpha\in J(\mu)}
    \bL^{\dd_0(\underline\alpha)}
    [\Filt^+(\underline\alpha)]z^{\deg\underline\alpha}=
    \bL^{\frac\delta2\langle\mu,\mu\rangle}z^{\delta n(\mu')}
    \sum_{d\ge0}[\Nil^{Coh,+}_{r,d}(\cO_X,\mu)]z^d.
\end{multline*}
We note that the last sum is exactly the generating function for motivic classes of stacks of non-twisted nilpotent pairs. Following~\cite{SchiffmannIndecomposable}, we denote it by $\Xi_\mu(z)$. By~\cite[(5.18)]{SchiffmannIndecomposable} (whose motivic version is valid by results from~\cite[Sect.~5 and~6]{FedorovSoibelmans}), we obtain
\[
    \Xi_\mu(z)=\bL^{(g-1)\langle\mu,\mu\rangle}J_{\mu,X}^{mot}(z)H_{\mu,X}^{mot}(z)
    \Exp\left(\frac{[X]}{\bL-1}\cdot\frac z{1-z}\right),
\]
where $\frac z{1-z}$ is expanded in positive powers of $z$ and $\Exp\colon\cMot(\kk)[[z]]^0\to\cMot(\kk)[[z]]$ is the plethystic exponent, see Section~\ref{sect:ExpLog} below. Note that this formula is based on interpreting the motivic classes $[\Filt^+(\underline\alpha)]$ in terms of Eisenstein series in motivic Hall algebras, see~\cite[Sect.~6]{FedorovSoibelmans}. Combining the last two displayed formulas, we obtain
\begin{equation}\label{eq:NilCoh}
    \bL^{\frac{r^2\delta}2}\sum_{d\ge0}[\Nil^{Coh,+}_{r,d}(\cL,\mu)]z^d=
    \bL^{\frac{2g-2+\delta}2\langle\mu,\mu\rangle}z^{\delta n(\mu')}
    J_{\mu,X}^{mot}(z)H_{\mu,X}^{mot}(z)
    \Exp\left(\frac{[X]}{\bL-1}\cdot\frac z{1-z}\right).
\end{equation}

\emph{Step 4.} We need one more lemma.
\begin{lemma}\label{lm:CohVec}
\[
    \sum_{d\ge0}[\Nil^{Coh,+}_{r,d}(\cL,\mu)]z^d=\left(\sum_{d\ge0}[\Nil^+_{r,d}(\cL,\mu)]z^d\right)
    \left(\sum_{d\ge0}[\Nil^{Coh,+}_{0,d}(\cL)]z^d\right).
\]
\end{lemma}
\begin{proof}
  The proof is similar to that of~\cite[Lemma~6.3.1]{FedorovSoibelmans}, so we just sketch it. Let $F$ be a coherent sheaf on $X$. Let $F_{tor}\subset F$ denote its torsion subsheaf, then $F\simeq F_{tor}\oplus F_{vec}$, where $F_{vec}$ is a vector bundle. For $e\in[0,d]$ consider the morphism $\pi_e\colon\Nil^{Coh,+}_{0,e}(\cL)\times\Nil^+_{r,d-e}(\cL,\mu)\to\Nil^{Coh,+}_{r,d}(\cL,\mu)$ sending $((F',\Psi'),(F'',\Psi''))$ to $(F'\oplus F'',\Psi'\oplus\Psi'')$. Denote the constructible image of this morphism by $\Nil^{Coh,+}_{r,d,e}$. It consists of pairs $(F,\Psi)$ such that $F_{tor}$ has degree $e$, so we get a stratification $\Nil^{Coh,+}_{r,d}(\cL,\mu)=\bigcup_{e=0}^d\Nil^{Coh,+}_{r,d,e}$. One shows that $[\Nil^{Coh,+}_{r,d,e}]=[\Nil^{Coh,+}_{0,e}(\cL)\times\Nil^+_{r,d-e}(\cL,\mu)]$ by comparing the fibers over a coherent sheaf $F_{tor}\oplus F_{vec}$, where $F_{tor}$ is a torsion sheaf of degree $e$, $F_{vec}$ is a~vector bundle of degree $d-e$. Thus,
  \[
   [\Nil^{Coh,+}_{r,d}(\cL,\mu)]=\sum_{e=0}^d[\Nil^{Coh,+}_{r,d,e}]=
    \sum_{e=0}^d[\Nil^{Coh,+}_{0,e}(\cL)][\Nil^+_{r,d-e}(\cL,\mu)],
  \]
  which is equivalent to our statement.
\end{proof}
Next,~\eqref{eq:NilCoh} with $\mu$ being the trivial partition gives $\sum_{d\ge0}[\Nil^{Coh,+}_{0,d}(\cL)]z^d=\Exp\left(\frac{[X]}{\bL-1}\cdot\frac z{1-z}\right)$. Combining this with the previous lemma and~\eqref{eq:NilCoh}, we get the desired result. The proof of Proposition~\ref{pr:NilpEnd} is complete.
\end{proof}

\subsection{Twisted parabolic pairs}\label{sect:ParabolicPairs}
Recall that $D$ is a finite set of rational points on~$X$ and $\cL$ is a line bundle on $X$. \emph{A parabolic bundle of level $D$ on $X$\/} is a pair $(E,E_{\bullet,\bullet})$, where for each $x\in D$ we have a~filtration
\[
    E|_x=E_{x,0}\supset E_{x,1}\supset\ldots\supset E_{x,N}=0\text{ for }N\gg0.
\]
In fact, if we view $D$ as a reduced divisor $D=\sum_{x\in D}x$, then this is the same as a parabolic bundle of level $D$ defined in Section~\ref{sect:ParStr}. An \emph{$\cL$-twisted nilpotent parabolic pair\/} $(\bE,\Psi)$ consists of a parabolic bundle $\bE=(E,E_{\bullet,\bullet})$ of level $D$ on $X$ and a globally nilpotent morphism $\Psi:E\to E\otimes\cL$ compatible with the parabolic structure in the sense that for all $x\in D$ and $j>0$ we have $\Psi(E_{x,j})\subset E_{x,j}\otimes\cL_x$. The nilpotence of $\Psi$ is automatic by Lemma~\ref{lm:nilpotent}, provided that $\deg\cL<0$.

We denote the stack of $\cL$-twisted nilpotent parabolic pairs by $\Pair^{nilp}(X,D,\cL)$. We suppress $X$ from the notation. We have a decomposition according to the class of the parabolic bundle (cf.~ Section~\ref{sect:ParStr})
\[
    \Pair^{nilp,-}(D,\cL)=\bigsqcup_{\gamma\in\Gamma_D,d\le0}\Pair^{nilp,-}_{\gamma,d}(D,\cL),
\]
where, as usual, the superscript ``$-$'' denotes pairs with non-positive underlying vector bundles. The stacks $\Pair^{nilp,-}_{\gamma,d}(D,\cL)$ are Artin stacks of finite type because the forgetful morphisms $\Pair^{nilp}_{\gamma,d}(D,\cL)\to\Bun_{\rk\gamma,d}(X)$ are schematic morphisms of finite type (see Section~\ref{sect:nonpositive}). Further, as in Section~\ref{sect:TwistedPairs}, we have a stratification by the generic Jordan type of the nilpotent endomorphism: $|\Pair^{nilp,-}_{\gamma,d}(D,\cL)|=\bigcup_{\mu\,\vdash\rk\gamma}\Pair^{nilp,-}_{\gamma,d}(D,\mu,\cL)$.

Recall from Section~\ref{sect:Macdonald} the modified Macdonald polynomials $\tilde H_\mu^{mot}(w_\bullet;z)\in\Sym_{\cMot(\kk)[z]}(w_\bullet)$ (which should not be confused with the series $H_{\mu,X}^{mot}(z)$). Recall also that for $\gamma=(r,r_{\bullet,\bullet})\in\Gamma_D$ we have set $w^\gamma=w^r\prod_{x\in D}\prod_{j=1}^{\infty}w_{x,j}^{r_{x,j}}$. The following is a twisted analogue of~\cite[Thm.~4.6.1]{FedorovSoibelmansParabolic} (notations as in Proposition~\ref{pr:NilpEnd}).

\begin{proposition}\label{pr:MotMellitPunctures} Let $\cL$ be a line bundle on $X$ of negative degree $-\delta$ and $\mu$ be a partition of an integer number $r$. Then we have in $\cMot(\kk)[\bL^{\frac12}][[\Gamma_D,z]]$
\[
    \bL^{\frac{r^2\delta}2}
    \sum_{\substack{\gamma\in\Gamma_D\\ \rk\gamma=r,\,d\le0}}[\Pair_{\gamma,d}^{nilp,-}(D,\mu,\cL)]w^\gamma z^{-d}=
    \bL^{\frac{(2g-2+\delta)\langle\mu,\mu\rangle}2}z^{\delta n(\mu')}
    w^rJ_{\mu,X}^{mot}(z)H_{\mu,X}^{mot}(z)\prod_{x\in D}\tilde H_\mu^{mot}(w_{x,\bullet};z).
\]
\end{proposition}
\begin{proof}
  Same as the proof of~\cite[Thm.~4.6.1]{FedorovSoibelmansParabolic} using Proposition~\ref{pr:NilpEnd} instead of~\cite[Thm.~1.4.1]{FedorovSoibelmans}.
\end{proof}

\subsection{Schiffmann's generating function}\label{sect:SchKernel} The series $\sum_{\mu\in\cP} q^{(g-1)\langle\mu,\mu\rangle}J_\mu(z)H_\mu(z)$, where $\cP$ stands for the set of partitions, first appeared in~\cite[(1.3)]{SchiffmannIndecomposable}. Its twisted version was considered in~\cite[Thm.~6.2]{MozgovoySchiffmann2020}. Its parabolic motivic version was used in~\cite[Sect.~1.3]{FedorovSoibelmansParabolic}. Now we define its twisted (a.k.a.~irregular) parabolic motivic version. For $\delta\in\Z_{\ge0}$ and a finite set $D$, put
\begin{equation}\label{eq:OmegaSchiffmann}
    \Omega^{Sch,mot}_{X,D,\delta}:=
    \sum_{\mu\in\cP}w^{|\mu|}(-\bL^\frac12)^{(2g-2+\delta)\langle\mu,\mu\rangle}z^{\delta n(\mu')}J_{\mu,X}^{mot}(z)
    H_{\mu,X}^{mot}(z)\prod_{x\in D}\tilde H_\mu^{mot}(w_{x,\bullet};z).
\end{equation}
This belongs to $\cMot(\kk)[\bL^\frac12][[\Gamma_D,z]]$.

Multiplying the formulas of Proposition~\ref{pr:MotMellitPunctures} by $(-1)^{|\mu|^2\delta}$, recalling that $|\mu|$ and $\langle\mu,\mu\rangle$ have the same parity, and summing over all partitions, we get the following corollary giving explicit formulas for the generating functions of motivic classes $[\Pair_{\gamma,d}^{nilp,-}(D,\cL)]$.
\begin{corollary}\label{cor:MotMellitPunctures} Let $\cL$ be a line bundle on $X$ of negative degree $-\delta$. Then in $\cMot(\kk)[\bL^{\frac12}][[\Gamma_D,z]]$ we have
\[
    \sum_{\gamma\in\Gamma_D,d\le0}(-\bL^\frac12)^{(\rk\gamma)^2\delta}[\Pair_{\gamma,d}^{nilp,-}(D,\cL)]w^\gamma z^{-d}=\Omega^{Sch,mot}_{X,D,\delta}.
\]
\end{corollary}
\begin{remark}
Replacing $-\bL^\frac12$ by $\bL^\frac12$ in the above formula and in $\Omega^{Sch,mot}_{X,D,\delta}$, one still gets a valid formula. The significance of the minus sign in the definition of $\Omega^{Sch,mot}_{X,D,\delta}$ will be clear in Section~\ref{sect:MellitSimplification}, see the proof of Lemma~\ref{lm:Admissible}.
\end{remark}

\section{Motivic classes of stacks of \texorpdfstring{$\epsilon$}{epsilon}-connections with underlying HN-nonpositive vector bundles}\label{sect:MotClass-} In Definition~\ref{def:ModSpacePartial} we defined the stacks $\Conn^{prtl}(\epsilon,X,\divisor,\divisor',\zeta)$ of parabolic $\epsilon$-connections with partially fixed normal forms, where $\divisor'=\sum_{x\in D'}n'_xx\le\divisor=\sum_{x\in D}n_xx$ are effective divisors on~$X$, the support of $\divisor$ is a subset $D$ of the set of rational points of $X$, and $\zeta\in\FNF(\divisor,\divisor')$ is a~partial normal form, see~\eqref{eq:FNF_def}. In this section, we will omit $\epsilon$ and $X$ from the notation if it does not lead to confusion.

The closed and open substacks $\Conn^{prtl}_{\gamma,d}(\divisor,\divisor',\zeta)\subset\Conn^{prtl}(\divisor,\divisor',\zeta)$, where $\gamma\in\Gamma_{D'}$ and $d\in\Z$ (see~\eqref{eq:X-}), are of infinite type unless $\divisor=\divisor'$ and $\epsilon\ne0$ (cf.~Corollary~\ref{cor:StabConn}). When $\Conn^{prtl}_{\gamma,d}(\divisor,\divisor',\zeta)$ is of infinite type, it is also of infinite motivic volume. Thus, we first calculate the motivic classes of the finite type stacks $\Conn_{\gamma,d}^{prtl,-}(\divisor,\divisor',\zeta)$ in $\cMot(\kk)$ classifying parabolic $\epsilon$-connections with nonpositive underlying vector bundles. The forgetful morphism $\Conn_{\gamma,d}^{prtl}(\divisor,\divisor',\zeta)\to\Bun(X)$ is of finite type, so, as explained in Section~\ref{sect:nonpositive}, the stacks $\Conn_{\gamma,d}^{prtl,-}(\divisor,\divisor',\zeta)$ are of finite type over $\kk$. Later, in Section~\ref{sect:Stabilization}, we will calculate the motivic classes of the semistable loci of $\Conn^{prtl}(\divisor,\divisor',\zeta)$ (see Theorem~\ref{th:AnswerInCmot}) as well as the motivic classes of $\Conn(\epsilon,\divisor,\zeta)$ when $\epsilon\ne0$ (see Theorem~\ref{th:AnswerConnInMot}).

The strategy of the present calculation is as follows. First, we assume that $\divisor'$ is a reduced divisor, in which case we relate the motivic classes under consideration to the motivic classes of twisted parabolic pairs calculated in Corollary~\ref{cor:MotMellitPunctures}. The comparison involves the Riemann--Roch theorem and Serre duality, cf.~\cite[Sect.~5.4]{FedorovSoibelmansParabolic}.

Then we proceed by induction on $\deg\divisor'$, comparing the motivic class of $\Conn_{\gamma,d}^{prtl,-}(\divisor,\divisor',\zeta)$ with the motivic class of $\Conn_{\gamma,d}^{prtl,-}(\divisor,\divisor'-x,\zeta)$, where $x\in D$ is any point with $n'_x\ge2$. We can only do this if $\gamma$ is full at $x$ (see Definition~\ref{def:full}) and $\zeta$ is non-resonant for $\gamma$ at $x$ (see Definition~\ref{def:non-resonant}). When $\divisor'<\divisor$, the calculation, roughly speaking, boils down to showing that motivic classes do not depend on $\zeta$ (this is the ``irregular isomonodromy'' property of motivic classes), which is done with the help of Lemma~\ref{lm:Independent}. In Section~\ref{sect:ExplNonpositive}, assuming that $\divisor'<\divisor$, we obtain explicit formulas for the motivic classes valid in $\cMot(\kk)$ using Corollary~\ref{cor:MotMellitPunctures}

The case $\divisor'=\divisor$, considered in Section~\ref{sect:FullyFixed}, is somewhat more difficult and is based on Lemma~\ref{lm:existIrregFull}; the explicit formulas are given in Proposition~\ref{pr:ExplicitFull-}.

\subsection{Stacks of parabolic $\epsilon$-connections with partially fixed normal forms and nonpositive underlying vector bundles} Recall that we denote the support of the divisor $\divisor'$ by $D'$, so that $D'=\{x\in D\colon n'_x>0\}$. In Section~\ref{sect:ParabolicPairs} we defined the stacks $\Pair^{nilp}(D',\cL)$, whenever $D'\subset X(\kk)$ is a set of rational points and $\cL$ is a line bundle on $X$. Their open substacks $\Pair^{nilp,-}_{r,d}(D',\cL)$, where $r\ge0$, $d\le0$, are of finite type by Section~\ref{sect:nonpositive}.

Let $g$ be the genus of $X$. Recall also from Section~\ref{sect:IrregFNF} that we have a subset $\Gamma_{\divisor',\zeta}\subset\Gamma_{D'}$ consisting of the classes $\gamma$ such that $\gamma$ is full at all $x\in D'$ with $n'_x\ge2$ and $\zeta$ is non-resonant for $\gamma$ at all such $x$.

\begin{proposition}\label{pr:main}
  Assume that for at least one $x\in D$ we have $n'_x<n_x$. Let $\gamma=(r,r_{\bullet,\bullet})\in\Gamma_{\divisor',\zeta}$. Then for all $d\in\Z_{\le0}$ we have in $\Mot(\kk)$:
\[
    [\Conn^{prtl,-}_{\gamma,d}(\divisor,\divisor',\zeta)]=\bL^m[\Pair_{\gamma,d}^{nilp,-}(D',\cO_X(-\divisor''))],
\]
where
\[
\begin{split}
    \divisor'':=\divisor-D'=\sum_{x\in D\colon n'_x>0}(n_x-1)x+\sum_{x\in D\colon n'_x=0}n_xx,\\m=(g-1+\deg\divisor'')r^2+\sum_{x\in D'}\left((1-n'_x)r+
    \sum_{i<j}r_{x,i}r_{x,j}\right).
\end{split}
\]
\end{proposition}
The proposition will be proved by induction on $\deg\divisor'$. The base case is when $\divisor'$ is a reduced divisor, that is, $\divisor'=\sum_{x\in D'}x$.

\subsection{Proof of Proposition~\ref{pr:main}: the case of a reduced divisor}\label{sect:CaseOfReduced}
Assume that $\divisor'$ is a reduced divisor. Then, the fibers of the forgetful morphism $\Conn^{prtl,-}_{\gamma,d}(\divisor,\divisor',\zeta)\to\Bun^{par,-}_{\gamma,d}(\divisor')$ are non-empty, see Lemma~\ref{lm:NoAtiyah}, . Thus, according to Lemma~\ref{lm:fibers}, for a field extension $K\supset\kk$ the fiber over $\bE\in\Bun^{par}(\divisor')(K)$ is a trivial torsor over $H^0(X_K,\END_0(\bE)\otimes\Omega_{X_K}(\divisor))$. Thus, its motivic class is equal to $F_\bE=\bL_K^{\dim H^0(X_K,\END_0(\bE)\otimes\Omega_{X_K}(\divisor))}$. On the other hand, since $\divisor'$ is reduced, we have a forgetful morphism $\Pair_{\gamma,d}^{nilp,-}(D',\cO_X(-\divisor''))\to\Bun^{par,-}_{\gamma,d}(\divisor')$. Its fiber over $\bE$ is equal to
$H^0(X_K,\END(\bE)\otimes\cO_{X_K}(-\divisor''))$, so its motivic class is equal to $F'_\bE=\bL_K^{\dim H^0(X_K,\END(\bE)\otimes\cO_{X_K}(-\divisor''))}$. Since $\divisor'$ is reduced, we have $\divisor''=\divisor-\divisor'$, so~\eqref{eq:RR} shows that $F_\bE=\bL_K^m F'_\bE$. Now the statement follows from Proposition~\ref{pr:MotivicClasses}(i) applied to $\cS=\cX=\Bun^{par,-}_{\gamma,d}(\divisor')$. \qed

\subsection{End of proof of Proposition~\ref{pr:main}: the induction step}
If $\divisor'$ is not a reduced divisor, then there is an $x\in D$ with $n'_x\ge2$. If we replace $\divisor'$ with $\divisor'-x$, then the number $m$ in the formulation of the proposition gets replaced with $m-r$. Thus, by induction on $\deg\divisor'$, the proposition is an obvious corollary of the base case and the following lemma.
\begin{lemma}\label{lm:AffineBundle2}
Assume that $x\in D$ is a point with $n'_x\ge2$. Recall from~\eqref{eq:Px} the projection
\[
    P_x\colon\FNF(\divisor,\divisor')\to\FNF(\divisor,\divisor'-x).
\]
Assume that $\gamma=(r,r_{\bullet,\bullet})\in\Gamma_{D'}$ is full at $x$ and that $\zeta_{x,\bullet}$ is non-resonant for $\gamma$ at $x$. Then we have in $\Mot(\kk)$
\begin{equation}\label{eq:ind}
        [\Conn_{\gamma,d}^{prtl,-}(\divisor,\divisor'-x,P_x(\zeta))]=
        \bL^r[\Conn_{\gamma,d}^{prtl,-}(\divisor,\divisor',\zeta)].
\end{equation}
\end{lemma}
\begin{proof}
Choose a coordinate $z$ near $x$. We construct a morphism
\[
    f\colon\Conn_{\gamma,d}^{prtl,-}(\divisor,\divisor'-x,P_x(\zeta))\to\A^r_\kk
\]
defined as the coefficient at $z^{n'_x}$ of the formal normal form of the $\epsilon$-connection at $x$ (cf.~Remark~\ref{rem:FNF}). More precisely, if $(E,E_{\bullet,\bullet},\nabla)$ is a $\kk$-point of $\Conn_{\gamma,d}^{prtl,-}(\divisor,\divisor'-x,P_x(\zeta))$, then, by our assumption on~$x$, $\gamma$, and $\zeta$ and Lemma~\ref{lm:DiagIffNonResonant}, $(E,\nabla)$ is generic at $x$. Thus by Corollary~\ref{cor:ExtParStr} we can uniquely extend the level $n'_x-1$ parabolic structure $E_{x,\bullet}$ to a level $n'_x$ parabolic structure $E'_{x,\bullet}$ at $x$ preserved by $\nabla$ in the sense of Section~\ref{sect:ParStrPreserved}. In other words, in a trivialization compatible with this parabolic structure, we have
\[
    \nabla=\epsilon d+A+B,
\]
where $A$ is an upper triangular matrix with entries being 1-forms with poles of order at most $n_x$, and~$B$ has a pole of order at most $n_x-n'_x$ at $x$. Let $A^{diag}$ denote the diagonal part of $A$, write
\[
    A^{diag}=\sum_{j=-n_x}^{n'_x-n_x-1}A_jz^j\,dz+\ldots,
\]
where $A_j\in\kk^r$. Then $f(E,E_{\bullet,\bullet},\nabla):=A_{n'_x-n_x-1}$. It is clear that $A^{diag}$ and therefore $A_j$ do not depend on the choice of the trivialization compatible with $E'_{x,\bullet}$. One checks that Corollary~\ref{cor:ExtParStr} works in families, so we indeed get a morphism $f$. We are going to show that the motivic classes of the fibers of $f$ are equal to $[\Conn_{\gamma,d}^{prtl,-}(\divisor,\divisor',\zeta)]$ and then apply Proposition~\ref{pr:MotivicClasses}(ii). In more detail, using the coordinate $z$ we get a projection $\FNF(\divisor,\divisor')\to\kk^{\Z_{>0}}$ sending $\theta_{\bullet,\bullet}$ to the sequence $(\theta'_j)$, where $\theta'_j$ is the coefficient of $\theta_{x,j}$ at $z^{n'_x-n_x-1}\,dz$. Using this projection we get an isomorphism $\FNF(\divisor,\divisor')=\FNF(\divisor,\divisor'-x)\times\kk^{\Z_{>0}}$. Now let $K\supset\kk$ be a field extension and let $\xi\in\A^r_\kk(K)=K^r$ be a point. View $P_x(\zeta)$ as an element of $\FNF(\divisor,\divisor'-x)\otimes_\kk K$ using the embedding $\kk\hookrightarrow K$. Under the previous isomorphism the pair $(P_x(\zeta),\xi)\in\FNF(\divisor,\divisor'-x)\otimes_\kk K\oplus K^{\Z_{>0}}$ corresponds to some $\theta\in\FNF(\divisor,\divisor')\otimes_\kk K$. By Corollary~\ref{cor:ExtParStr} applied to $X_K$ and the construction, we see that the $\xi$-fiber of $f$ is identified with $\Conn_{\gamma,d}^{prtl,-}(X_K,\divisor,\divisor',\theta)$.

We claim that
\begin{equation}\label{eq:MotFibers}
    [\Conn_{\gamma,d}^{prtl,-}(X_K,\divisor,\divisor',\theta)]=
    [\Conn_{\gamma,d}^{prtl,-}(X_K,\divisor,\divisor',\zeta)]=
    [\Conn_{\gamma,d}^{prtl,-}(X,\divisor,\divisor',\zeta)\times_\kk K].
\end{equation}
To show the first equality, consider the forgetful morphism
\[
    \Conn_{\gamma,d}^{prtl,-}(X_K,\divisor,\divisor',\theta)\to\Bun^{par,-}_{\gamma,d}(X_K,\divisor').
\]
According to Lemma~\ref{lm:Independent} applied to $X_K$, the image of this morphism depends only on $P_x(\theta)$. On the other hand, the motivic classes of the non-empty fibers are independent of $\theta$ by Lemma~\ref{lm:fibers}. The claim now follows from Proposition~\ref{pr:MotivicClasses}(i) applied with $m=0$. The second equality in~\eqref{eq:MotFibers} is clear. Thus we can apply Proposition~\ref{pr:MotivicClasses}(ii) with $\cS=\cX=\A_\kk^r$ to complete the proof of the lemma.
\end{proof}

The lemma completes the proof of Proposition~\ref{pr:main}. \qed

We remark that $A^{diag}$ in the proof above coincides with $(\tilde\zeta_1,\ldots,\tilde\zeta_{\rk E})$ in Lemma~\ref{lm:diag} up to terms with poles of order at most $n_x-n'_x$ and that $A_j$ with $j<n'_x-n_x-1$ are equal to the corresponding coefficients of $P_x(\zeta)$.

\subsection{Explicit formulas for the motivic classes of parabolic $\epsilon$-connections with nonpositive underlying vector bundles: case of partially fixed normal forms}\label{sect:ExplNonpositive} Recall that we have divisors $\divisor=\sum_{x\in D}n_xx$ and $\divisor'=\sum_{x\in D}n'_xx$ such that for all $x\in D$ we have $0<n_x$ and $n'_x\le n_x$. For $\zeta\in\FNF(\divisor,\divisor')$ we defined a subset $\Gamma_{\divisor',\zeta}$ of the monoid $\Gamma_{D'}$ in Section~\ref{sect:IrregFNF}. Our next goal is to give explicit formulas for the motivic classes of stacks $[\Conn^{prtl,-}_{\gamma,d}(\divisor,\divisor',\zeta)]$, where $\divisor'<\divisor$ and $\gamma\in\Gamma_{\divisor',\zeta}$. Recall that $g$ is the genus of $X$ and in Section~\ref{sect:SchKernel} we defined the irregular motivic Schiffmann's generating function $\Omega^{Sch,mot}_{X,D',\delta}$, where $D'$ is any finite set and $\delta\ge0$.

\begin{proposition}\label{pr:ExplicitPartial-}
Assume that $\divisor'<\divisor$ and $\zeta\in\FNF(\divisor,\divisor')$. Let $D':=\{x\in D\colon n'_x>0\}$ be the support of $\divisor'$ and set
\[
    \delta:=\deg\divisor-|D'|\in\Z_{>0}.
\]
Next, for $\gamma\in\Gamma_{D'}$ recall from~\eqref{eq:BetterChi3}
\begin{equation}\label{eq:chi}
    \chi(\gamma)=\chi(\gamma,g,\divisor,\divisor'):=
    (2g-2+\deg(\divisor-\divisor'))r^2+r\sum_{x\in D'}(1-n'_x)+2\sum_{x\in D',i<j}n'_xr_{x,i}r_{x,j}.
\end{equation}
Then for $\gamma=(r,r_{\bullet,\bullet})\in\Gamma_{\divisor',\zeta}$ and $d\le0$ the motivic class of $\Conn^{prtl,-}_{\gamma,d}(\divisor,\divisor',\zeta)$ in $\cMot(\kk)$ is equal to the coefficient of $(-\bL^\frac12)^{\chi(\gamma)}\Omega^{Sch,mot}_{X,D',\delta}$ at $w^\gamma z^{-d}$.
\end{proposition}
\begin{proof}
Since $\gamma$ if full at all points with $n'_x>1$, we have
\begin{multline*}
    \chi(\gamma)=(2g-2+\deg(\divisor-\divisor'))r^2+r\sum_{x\in D'}(1-n'_x)+
    2\sum_{\substack{x\in D'\\i<j}}r_{x,i}r_{x,j}+2\sum_{\substack{x\colon n'_x>1\\i<j}}(n'_x-1)r_{x,i}r_{x,j}\\=
    (2g-2+\deg(\divisor-\divisor'))r^2+r\sum_{x\in D'}(1-n'_x)+
    2\sum_{\substack{x\in D'\\i<j}}r_{x,i}r_{x,j}+\sum_{x\colon n'_x>1}(n'_x-1)r(r-1)\\=
    (2g-2+\delta)r^2+2r\sum_{x\in D'}(1-n'_x)+2\sum_{\substack{x\in D'\\i<j}}r_{x,i}r_{x,j}=2m-r^2\delta,
\end{multline*}
where $m$ is from Proposition~\ref{pr:main}. Thus, the proposition gives in $\Mot(\kk)$:
\[
    [\Conn^{prtl,-}_{\gamma,d}(\divisor,\divisor',\zeta)]=
    (-\bL^\frac12)^{\chi(\gamma)}(-\bL^\frac12)^{r^2\delta}[\Pair_{\gamma,d}^{nilp,-}(D',\cO_X(-\divisor''))].
\]
It remains to apply the homomorphism $\Mot(\kk)\to\cMot(\kk)$ and use Corollary~\ref{cor:MotMellitPunctures}.
\end{proof}

\subsection{Stacks of parabolic $\epsilon$-connections with fully fixed formal normal forms}\label{sect:FullyFixed} We now consider the case $\divisor'=\divisor$, that is, the case of fully fixed formal normal forms. We note that if we set $\divisor'=\divisor$ in the previous formulas, then $\delta$ becomes $\deg\divisor-|D|\in\Z_{>0}$, that is, the irregularity of the divisor $\divisor$, while~$\chi$ becomes~\eqref{eq:BetterChi2}. However, as we will see momentarily, the final formula is slightly more complicated than the one obtained via this naive substitution. It involves the mutually inverse plethystic exponent and logarithm:
\[
 \Exp\colon\cMot(\kk)[\bL^\frac12][[\Gamma_D,z]]^0\rightleftarrows1+\cMot(\kk)[\bL^\frac12][[\Gamma_D,z]]^0
 \colon\Log.
\]
Here the superscript ``$^0$'' stands for series with zero constant term. The plethystic exponent is defined as follows (cf.~Section~\ref{sect:Plethystic} and~\cite[Sect.~3.4]{FedorovSoibelmansParabolic}). If $Y$ is a quasi-projective $\kk$-scheme, $n\in\frac12\Z$, $\gamma\in\Gamma_D$, and $j\in\Z_{\ge0}$, we set $\Exp([Y]\bL^nw^\gamma z^j):=\zeta_Y(\bL^nw^\gamma z^j)$. We extend this to $\Mot(\kk)[\sqrt\bL][[\Gamma_D,z]]^0$ using the formula $\Exp(A\pm B)=\Exp(A)\Exp(B)^{\pm1}$ and continuity (see also Section~\ref{sect:ExpLog} below).

Recall that for $\gamma\in\Gamma_D$ and $\zeta\in\FNF(\divisor)$ we defined in~\eqref{eq:gamma_star_zeta} the product $\gamma\star\zeta\in\kk$.
\begin{proposition}\label{pr:ExplicitFull-}
Assume that $\zeta\in\FNF(\divisor)$, $\gamma'=(r',r'_{\bullet,\bullet})\in\Gamma_{\divisor,\zeta}$, and $d'\le0$. Define the irregularity $\delta:=\deg\divisor-|D|\in\Z_{>0}$ and recall from~\eqref{eq:BetterChi2} the function
\[
    \chi(\gamma)=\chi(\gamma,g,\divisor):=(2g-2)r^2-\delta r+2\sum_{x\in D,i<j}n_xr_{x,i}r_{x,j}.
\]
Then the motivic class of $\Conn^-_{\gamma',d'}(\divisor,\zeta)$ in $\cMot(\kk)$ is equal to the coefficient at $w^{\gamma'}z^{-d'}$ in
\[
    (-\bL^\frac12)^{\chi(\gamma')}\Exp\left(\Bigl(\bL\Log\Omega^{Sch,mot}_{X,D,\delta}\Bigr)_{-\epsilon d+\gamma\star\zeta=0}\right),
\]
where the subscript $-\epsilon d+\gamma\star\zeta=0$ stands for the sum of the monomials whose exponents satisfy this condition.
\end{proposition}

\subsection{Proof of Proposition~\ref{pr:ExplicitFull-}: compatible parabolic bundles}\label{sect:ParPairs} Fix a point $x$ with $n_x\ge2$ and $\zeta\in\FNF(\divisor)$. Then we have $P_x(\zeta)\in\FNF(\divisor,\divisor-x)$, see~\eqref{eq:Px}. We will also need the truncation morphisms $T_x\colon\Bun^{par}(\divisor)\to\Bun^{par}(\divisor-x)$. Recall from Section~\ref{sect:ParStrPreserved} the notion of a parabolic structure preserved by an $\epsilon$-connection.

\begin{definition}\label{def:Bun_zeta}
Let $\Bun^{par}(\divisor,\zeta)$ denote the constructible subset of $\Bun^{par}(\divisor)$ consisting of parabolic bundles $(E,E_{\bullet,\bullet})$ such that there is an $\epsilon$-connection $\nabla$ on $E$ with poles bounded by $\divisor$ satisfying two conditions:

(i) $(T_x(E,E_{\bullet,\bullet}),\nabla)\in\Conn^{prtl}(\divisor,\divisor-x,P_x(\zeta))$;

(ii) $E_{x,\bullet}$ is preserved by $\nabla$.
\end{definition}
We note that (ii) is weaker than the compatibility in Definition~\ref{def:ModSpace}, since we do not require the diagonal blocks to be scalar (which is, however, automatic if $\gamma$ is full at $x$ because every $1\times1$ matrix is scalar). Note also that $\Bun^{par}(\divisor,\zeta)$ decomposes according to the class of $\bE$:
\[
    \Bun^{par}(\divisor,\zeta)=\bigsqcup_{\gamma\in\Gamma_D,d\in\Z}\Bun_{\gamma,d}^{par}(\divisor,\zeta).
\]

\begin{lemma}\label{lm:ImageCs}
Assume that $\gamma\in\Gamma_D$ if full at $x$ and $d\in\Z$ is such that $\epsilon d+\gamma\star\zeta=0$. Then the image of the forgetful morphism $\Conn_{\gamma,d}(\divisor,\zeta)\to\Bun^{par}(\divisor)$ is the intersection of $\Bun_{\gamma,d}^{par}(\divisor,\zeta)$ with the constructible subset of $(\epsilon,\zeta)$-isoslopy parabolic bundles.
\end{lemma}
\begin{proof}
Assume that $\bE=(E,E_{\bullet,\bullet})\in\Bun_{\gamma,d}^{par}(\divisor,\zeta)(K)$, where $K\subset\kk$ is a field extension, is $(\epsilon,\zeta)$-isoslopy of degree zero, we need to show that $\bE$ is in the image of $\Conn_{\gamma,d}(\divisor,\zeta)$. Base changing, we may assume that $K=\kk$.

Since $\gamma$ is full at $x$, by definition of $\Bun_{\gamma,d}^{par}(\divisor,\zeta)$ there is $\theta\in\FNF(\divisor)$ with $P_x(\theta)=P_x(\zeta)$ and an $\epsilon$-connection $\nabla$ on $E$ such that $(E,E_{\bullet,\bullet},\nabla)\in\Conn(\divisor,\theta)(\kk)$ (see the remark after Definition~\ref{def:Bun_zeta}). Now, by Lemma~\ref{lm:existIrregFull}, $\bE$ is in the image of $\Conn_{\gamma,d}(\divisor,\zeta)$.

The other inclusion follows from the definition of $\Bun^{par}(\divisor,\zeta)$ and Corollary~\ref{cor:indecomposable summand}.
\end{proof}
We introduce some notation. Let $\cS\subset\cY$ be a constructible subset of a stack $\cY$ and $f\colon\cX\to\cY$ be a~morphism of stacks. Then $\cX\times_\cY\cS$ stands for the constructible subset $f^{-1}(\cS)$ of the stack $\cX$.
We denote by $\Bun^{par,(\epsilon,\zeta)-iso}(\divisor,\zeta)\subset\Bun^{par}(\divisor,\zeta)$ the constructible subset of $(\epsilon,\zeta)$-isoslopy parabolic bundles (not necessarily of slope zero). Denote by $\Pair(\divisor)$ the stack of pairs $(\bE,\Psi)$, where $\bE\in\Bun^{par}(\divisor)$, and $\Psi\in\End(\bE)$ (cf.~Section~\ref{sect:ParabolicPairs}). Let $\Pair^{nilp}(\divisor)\subset\Pair(\divisor)$ be the substack of pairs $(\bE,\Psi)$, where the endomorphism $\Psi$ is globally nilpotent. The following lemma will allow us to relate the motivic classes of stacks $[\Conn^-_{\gamma,d}(\divisor,\zeta)]$ with the already calculated motivic classes of $[\Conn^{prtl,-}_{\gamma,d}(\divisor,\divisor-x,P_x(\zeta))]$ (cf.~Lemma~\ref{lm:AffineBundle2}).
\begin{lemma}\label{lm:RR}
Assume that $\gamma=(r,r_{\bullet,\bullet})\in\Gamma_D$ is full at $x\in D$ and that $\zeta$ is non-resonant for $\gamma$ at~$x$. Set
\[
    \beta(\gamma):=\frac{\chi(\gamma)+\delta r}2=(g-1)r^2+\sum_{y\in D}\sum_{i<j}n_yr_{y,i}r_{y,j}.
\]
We have in $\Mot(\kk)$:

(i)
\[
    [\Conn^-_{\gamma,d}(\divisor,\zeta)]=
    \begin{cases}
    \bL^{\beta(\gamma)}[\Pair^-(\divisor)\times_{\Bun^{par}(\divisor)}\Bun_{\gamma,d}^{par,(\epsilon,\zeta)-iso}(\divisor,\zeta)] &\text{ if }
    \epsilon d+\gamma\star\zeta=0\\
    0 &\text{ otherwise }.
    \end{cases}
\]
(ii)
\[
    [\Conn^{prtl,-}_{\gamma,d}(\divisor,\divisor-x,P_x(\zeta))]=
    \bL^{r+\beta(\gamma)}[\Pair^{nilp,-}(\divisor)\times_{\Bun^{par}(\divisor)}\Bun_{\gamma,d}^{par}(\divisor,\zeta)].
\]
\end{lemma}
\begin{proof}
  (i) The case $\epsilon d+\gamma\star\zeta\ne0$ follows from Corollary~\ref{cor:indecomposable summand}. Assume that $\epsilon d+\gamma\star\zeta=0$. By Lemma~\ref{lm:ImageCs}, the image of $\Conn^-_{\gamma,d}(\divisor,\zeta)$ in $\Bun^{par}(\divisor)$ is exactly $\Bun_{\gamma,d}^{par,(\epsilon,\zeta)-iso,-}(\divisor,\zeta)$. Thus, by Proposition~\ref{pr:MotivicClasses}(i), we just need to check that for a field extension $K\supset\kk$ and a parabolic bundle $\bE\in\Bun_{\gamma,d}^{par,(\epsilon,\zeta)-iso}(\divisor,\zeta)(K)$, the $\bE$-fiber of $\Conn_{\gamma,d}(\divisor,\zeta)$ differs from the fiber of $\Pair(\divisor)$ by a factor of $\bL_K^{\beta(\gamma)}$.

  By Lemma~\ref{lm:fibers}, the $\bE$-fiber of $\Conn_{\gamma,d}(\divisor,\zeta)$ is a trivial torsor over $H^0(X_K,\END_0(\bE)\otimes\Omega_{X_K}(\divisor))$, where $\END_0(\bE)$ stands for the sheaf of endomorphisms of $\bE$ inducing zero on $\gr E_{y,\bullet}$ for all $y\in D$. Hence, the motivic class of the fiber is $\bL^m$, where
  \[
    m=\dim H^0(X_K,\END_0(\bE)\otimes\Omega_{X_K}(\divisor)).
  \]
  On the other hand, the fiber of the projection $\Pair(\divisor)\to\Bun^{par}(\divisor)$ over $\bE$ is $H^0(X_K,\END(\bE))$ and the dimensions differ by $\beta(\gamma)$ by~\eqref{eq:RR} applied with $\divisor'=\divisor$. The statement follows.

  (ii) Consider the morphism $\Conn^{prtl}_{\gamma,d}(\divisor,\divisor-x,P_x(\zeta))\to\Bun^{par}(\divisor)$ assigning to $(E,E_{\bullet,\bullet},\nabla)$ the parabolic bundle $(E,E'_{\bullet,\bullet})$, where $E'_{\bullet,\bullet}$ is the unique extension of the level $\divisor-x$ parabolic structure~$E_{\bullet,\bullet}$ to a level $\divisor$ parabolic structure compatible with $\nabla$, see Corollary~\ref{cor:ExtParStr} (one checks that this extension works in families). Then the image of this morphism is $\Bun_{\gamma,d}^{par}(\divisor,\zeta)$.

  Thus, by Proposition~\ref{pr:MotivicClasses}(i), we just need to check that for a field extension $K\supset\kk$ and a parabolic bundle $\bE\in\Bun_{\gamma,d}^{par}(\divisor,\zeta)(K)$, the $\bE$-fiber of $\Conn^{prtl}_{\gamma,d}(\divisor,\divisor-x,\zeta)$ differs from the fiber of $\Pair^{nilp}(\divisor)$ by a factor of $\bL_K^{r+\beta(\gamma)}$.

  Let $\END_{0,x}(\bE)\subset\END(\bE)$ be the subsheaf consisting of endomorphisms inducing $0$ on $\gr\bE_{y,\bullet}$ when $y\ne x$ and inducing zero to order $n_x-1$ on $\gr\bE_{x,\bullet}$ (thus, $\END_{0,x}(\bE)$ in general strictly contains $\END_0(\bE)$). Let $(E_{x,\bullet})_{red}$ be the level~1 parabolic structure at $x$ obtained by restricting $E_{x,\bullet}$ from $n_xx$ to $x$. Denote by $\END_x(\bE)\subset\END(\bE)$ the subsheaf consisting of endomorphisms inducing $0$ on the associated graded of $(E_{x,\bullet})_{red}$. We note for the future that we have an exact sequence
  \[
  0\to\END_x(\bE)\to\END(\bE)\to\kk_x^{\rk\gamma}\to0,
  \]
  where $\kk_x$ is the skyscrapper sheaf supported at $x$. (We are using that $\gamma$ is full at $x$.) Thus, using Lemma~\ref{lm:RR-dim} with $\divisor'=\divisor$, we get
  \begin{equation}\label{eq:beta_x}
    \chi(\END_x(\bE))=\chi(\END(\bE))-\rk\gamma=-r-\beta(\gamma),
  \end{equation}
  where $\chi$ stands for the Euler characteristic.

  Similarly to part~(i) and to Lemma~\ref{lm:fibers}, the $\bE$-fiber of $\Conn^{prtl}_{\gamma,d}(\divisor,\divisor-x,\zeta)$ is a trivial torsor over $H^0(X_K,\END_{0,x}(\bE)\otimes\Omega_{X_K}(\divisor))$. Consider the projection $\Pair^{nilp}(\divisor)\to\Bun^{par}(\divisor)$. Its fiber over $\bE\in\Bun^{par}(\divisor)(K)$ consists of all nilpotent endomorphisms of $\bE$. We claim that this fiber is equal to $H^0(X_K,\END_x(\bE))$. Indeed, every nilpotent endomorphism induces zero on $\gr(E_{x,\bullet})_{red}$ because $\gamma$ is full at $x$. Conversely, every endomorphism contained in $H^0(X_K,\END_x(\bE))$ is globally nilpotent because the $i$-th coefficient of its characteristic polynomials is a section of $H^0(X_K,\cO_{X_K}(-ix))=0$.

  Next, we have a non-degenerate trace pairing between the sheaves $\END_{0,x}(\bE)\otimes\cO_{X_K}(\divisor)$ and $\END_x(\bE)$. Thus, by Serre duality, the Riemann--Roch theorem, and~\eqref{eq:beta_x}, we have
  \[
    \dim H^0(X_K,\END_{0,x}(\bE)\otimes\Omega_{X_K}(\divisor))-\dim H^0(X_K,\END_x(\bE))=-\chi(\END_x(\bE))=r+\beta(\gamma).
  \]
This completes the proof of Lemma~\ref{lm:RR}.
\end{proof}

\subsection{Proof of Proposition~\ref{pr:ExplicitFull-}: relating parabolic pairs with nilpotent parabolic pairs} The goal of this section is to relate the motivic classes of moduli stacks of parabolic pairs with the motivic classes of moduli stacks of nilpotent parabolic pairs. This is similar to~\cite[Prop.~3.8.1]{FedorovSoibelmans} and~\cite[Cor.~4.21]{FedorovSoibelmansParabolic} and uses the notion of power structure:
\[
    \Pow:(1+\cMot(\kk)[[\Gamma_D,z]]^0)\times\cMot(\kk)\to1+\cMot(\kk)[[\Gamma_D,z]]^0:
    (f,A)\mapsto\Exp(A\Log(f)).
\]
The pairs are related to nilpotent pairs via the Jordan decomposition. Thus, we need some results concerning direct sums of parabolic bundles. For every field extension $K\supset\kk$, $\Bun^{par}(\divisor)(K)$ is naturally an additive category. Consider its full subcategory $\Bun^{par}(\divisor,\zeta)(K)$.

\begin{lemma}\label{lm:cS}
\stepzero\noindstep\label{cS:i} The category $\Bun^{par}(\divisor,\zeta)(K)$ is closed under direct sums.

\noindstep\label{cS:ii} The category $\Bun^{par}(\divisor,\zeta)(K)$ is closed under taking a direct summand.

\noindstep\label{cS:iii} Let $K'\supset K$ be a finite field extension and $\bE\in\Bun^{par}(\divisor)(K')$. Then $\bE$ is in $\Bun^{par}(\divisor,\zeta)(K')$ if and only if $R_{K'/K}(\bE)$ is in $\Bun^{par}(\divisor,\zeta)(K)$, where $R_{K'/K}$ is the pushforward along the morphism $\Spec K'\to\Spec K$.
\end{lemma}
\begin{proof}
\eqref{cS:i} is clear because $\Conn^{prtl}(\divisor,\divisor-x,P_x(\zeta))$ is closed under direct sums.

To prove~\eqref{cS:ii}, let $\bF=(F,F_{\bullet,\bullet})$ be a direct summand of $\bE=(E,E_{\bullet,\bullet})$ and let $\nabla\colon E\to E\otimes\Omega_{X_K}(\divisor)$ be an $\epsilon$-connection such that $(\bE,\nabla)\in\Conn^{prtl}(\divisor,\divisor-x,P_x(\zeta))(K)$ and condition~(ii) of Definition~\ref{def:Bun_zeta} is satisfied. Define the $\epsilon$-connection $\nabla|_\bF$ as the composition
\[
    F\hookrightarrow E\xrightarrow{\nabla}E\otimes\Omega_{X_K}(\divisor)\twoheadrightarrow F\otimes\Omega_{X_K}(\divisor).
\]
Then $(\bF,\nabla|_\bF)\in\Conn^{prtl}(\divisor,\divisor-x,P_x(\zeta))(K)$ and it is also cleat that condition~(ii) of Definition~\ref{def:Bun_zeta} is satisfied.

\eqref{cS:iii} Let $\bE\in\Bun^{par}(\divisor,\zeta)(K')$ and assume that $R_{K'/K}(\bE)\in\Bun^{par}(\divisor,\zeta)(K)$. Let $R_{K'/K}(\bE)_{K'}$ denote the pullback of $R_{K'/K}(\bE)$ to $X_{K'}$. It is clear that $R_{K'/K}(\bE)_{K'}\in\Bun^{par}(\divisor,\zeta)(K')$. However, since $K$ is a direct summand of $K'\otimes_KK'$, $\bE$ is a direct summand of $R_{K'/K}(\bE)_{K'}$. Now by part~\eqref{cS:ii}, $\bE\in\Bun^{par}(\divisor,\zeta)(K')$. The other direction is clear.
\end{proof}

Let $\cX$ be a $\kk$-stack graded by $\Gamma_D\times\Z_{\le0}$, that is, we are given a decomposition into closed and open substacks $\cX=\bigsqcup_{\gamma\in\Gamma_D,d\in\Z_{\le0}}\cX_{\gamma,d}$. Fix and integer $\delta$. We set
\begin{equation}\label{eq:gradedStacks}
  [\cX]:=\sum_{\gamma\in\Gamma_D,d\le0}(-1)^{\delta\rk\gamma}[\cX_{\gamma,d}]w^\gamma z^{-d}\in\Mot(\kk)[[\Gamma_D,z]].
\end{equation}
This convention is similar to one in~\cite[Sect.~3.3]{FedorovSoibelmansParabolic} except for the $(-1)^{\delta\rk\gamma}$ multiple and replacing $z$ with $z^{-1}$.
\begin{proposition}\label{pr:NilpPow} We have in $\cMot(\kk)[[\Gamma_D,z]]$
\[
    [\Pair^-(\divisor)\times_{\Bun^{par}(\divisor)}\Bun^{par}(\divisor,\zeta)]=
    \Pow([\Pair^{nilp,-}(\divisor)\times_{\Bun^{par}(\divisor)}\Bun^{par}(\divisor,\zeta)],\bL).
\]
\end{proposition}
We note that the identity is valid in $\Mot(\kk)[[\Gamma_D,z]]$ as well but this coarser version suffices for our purposes.
\begin{proof}
  This is similar to~\cite[Prop.~3.8.1]{FedorovSoibelmans} in view of the previous lemma. In more detail, let $(\bE,\Psi)$ be a~$K$-point of $\Pair^-(\divisor)\times_{\Bun^{par}(\divisor)}\Bun^{par}(\divisor,\zeta)$. As in~\cite[Prop.~3.8.1]{FedorovSoibelmans}, we can uniquely decompose
  \[
    (\bE,\Psi)\xrightarrow{\simeq}\bigoplus_i R_{\kk(x_i)/K}(\bE_i,x_i\Id+\Psi_i),
  \]
  where $x_i$ are distinct closed points of $\A_K^1$ (the eigenvalues of $\Psi$), $(\bE_i,\Psi_i)$ are $\kk(x_i)$-points of the stack $\Pair^{nilp,-}(\divisor)$, $\kk(x_i)\supset K$ are the residue field of $x_i$. By Lemma~\ref{lm:cS}(\ref{cS:ii},\ref{cS:iii}), $\bE_i\in\Bun^{par}(\divisor,\zeta)(k(x_i))$ for all $i$.

  Conversely, if $x_i$ are distinct closed points of $\A_K^1$ and $(\bE_i,\Psi_i)$ are $\kk(x_i)$-points of the stack
  \[
    \Pair^{nilp,-}(\divisor)\times_{\Bun^{par}(\divisor)}\Bun^{par}(\divisor,\zeta),
  \]
  then the pair $\bigoplus_i R_{\kk(x_i)/K}(\bE_i,x_i\Id+\Psi_i)$ is a $K$-point of $\Pair^-(\divisor)\times_{\Bun^{par}(\divisor)}\Bun^{par}(\divisor,\zeta)$ according to Lemma~\ref{lm:cS}(\ref{cS:i},\ref{cS:iii}). It remains to use a version of~\cite[Lemma~3.8.2]{FedorovSoibelmans}.
\end{proof}

\subsection{Proof of Proposition~\ref{pr:ExplicitFull-}: a product formula for isoslopy pairs} It remains to relate the isoslopy parabolic bundles with all parabolic bundles. This is accomplished via a product formula similar to~\cite[Prop.~3.5.1]{FedorovSoibelmans} and~\cite[Prop.~5.6]{FedorovSoibelmansParabolic}. For this, we need a lemma.

\begin{lemma}
    For a field extension $K\supset\kk$, idempotents in the category $\Bun^{par}(\divisor,\zeta)(K)$ are split.
\end{lemma}
\begin{proof}
  Note that idempotents are split in the category of vector bundles on $X$. It follows that they are split in $\Bun^{par}(\divisor)$. It remains to use Lemma~\ref{lm:cS}\eqref{cS:ii}.
\end{proof}

\begin{proposition}\label{pr:SlopeProd} We have in $\Mot(\kk)[[\Gamma_D,z]]$
\begin{multline*}
    [\Pair^-(\divisor)\times_{\Bun^{par}(\divisor)}\Bun^{par}(\divisor,\zeta)]\\=\prod_{\tau\in\kk}\left(
        \sum_{\substack{\gamma\in\Gamma_D,d\le0\\ \epsilon d+\gamma\star\zeta=\tau\rk\gamma}}(-1)^{\delta\rk\gamma}[\Pair^-(\divisor)\times_{\Bun^{par}(\divisor)} \Bun_{\gamma,d}^{par,(\epsilon,\zeta)-iso}(\divisor,\zeta)]w^\gamma z^{-d}
    \right).
\end{multline*}
\end{proposition}
\begin{proof}
It follows from the previous lemma and a version of~\cite[Thm.~3]{Atiyah-KrullSchmidt} (cf.~also~\cite[Prop.~3.1.2]{FedorovSoibelmans}) that if $K$ is an extension of $\kk$, then any $\bE\in\Bun^{par}(\divisor,\zeta)(K)$ can be written as the direct sum of indecomposable objects: $\bE=\oplus_i\bE_i$, where $\bE_i$ are unique up to isomorphism. There rest is as in~\cite[Prop.~3.5.1]{FedorovSoibelmans}.
\end{proof}

\subsection{End of proof of Proposition~\ref{pr:ExplicitFull-}}
Our first goal is to combine Lemma~\ref{lm:RR}(ii) with Proposition~\ref{pr:ExplicitPartial-} to compute the motivic classes $[\Pair^{nilp,-}(\divisor)\times_{\Bun^{par}(\divisor)}\Bun_{\gamma,d}^{par}(\divisor,\zeta)]$. Recall that $\delta=\deg\divisor-|D|$.

\begin{lemma}\label{lm:ExplPair}
Assume that $\gamma\in\Gamma_{\divisor,\zeta}$. Then the motivic class of $\Pair^{nilp,-}(\divisor)\times_{\Bun^{par}(\divisor)}\Bun_{\gamma,d}^{par}(\divisor,\zeta)$ in $\cMot(\kk)$ is equal to the coefficient of $(-\bL^\frac12)^{-\delta\rk\gamma}\Omega^{Sch,mot}_{X,D,\delta}$ at $w^\gamma z^{-d}$.
\end{lemma}
\begin{proof}
Set $r:=\rk\gamma$. We want to apply Proposition~\ref{pr:ExplicitPartial-} with $\divisor'=\divisor-x$. Note that
\begin{multline*}
\chi(\gamma,\divisor,\divisor-x)=(2g-2+1)r^2+r(1-\delta)+2\sum_{\substack{y\ne x\\i<j}}n_yr_{y,i}r_{y,j}
+2\sum_{i<j}(n_x-1)r_{x,i}r_{x,j}\\
=\chi(\gamma)+r^2+r-r(r-1)=\chi(\gamma)+2r,
\end{multline*}
where we used that $\gamma$ is full at $x$. Denote by $A\in\cMot(\kk)$ the coefficient of $\Omega^{Sch,mot}_{X,D,\delta}$ at $w^{\gamma}z^{-d}$. Using Lemma~\ref{lm:RR}(ii) and Proposition~\ref{pr:ExplicitPartial-} we obtain
\begin{multline*}
    [\Pair^{nilp,-}(\divisor)\times_{\Bun^{par}(\divisor)}\Bun_{\gamma,d}^{par}(\divisor,\zeta)]=
    \bL^{-r-\beta(\gamma)}[\Conn^{prtl,-}_{\gamma,d}(\divisor,\divisor-x,P_x(\zeta))]\\=
    \bL^{-r-\beta(\gamma)}(-\bL^\frac12)^{\chi(\gamma,\divisor,\divisor-x)}A=
    (-\bL^\frac12)^{\chi(\gamma)-2\beta(\gamma)}=(-\bL^\frac12)^{-\delta r}A.
\end{multline*}
\end{proof}

We need a simple lemma, which will be used a few times.
\begin{lemma}\label{lm:GoodExp}
Fix $\zeta\in\FNF(\divisor)$. Let $H_1(w,z),H_2(w,z)\in\cMot(\kk)[[\Gamma_D,z]]^0$ be two motivic series without constant terms. Then the coefficients of $H_1$ and $H_2$ coincide at all $w^\gamma z^d$ with $\gamma\in\Gamma_{\divisor,\zeta}$ if and only if the same condition is satisfied for $\Exp(H_1)$ and $\Exp(H_2)$.
\end{lemma}
\begin{proof}
     By Lemma~\ref{lm:Summand}, any summand of $\gamma\in\Gamma_{\divisor,\zeta}$ belongs to $\Gamma_{\divisor,\zeta}$. Now our lemma (in both directions) follows by induction from the fact that the coefficient at $w^\gamma z^d$ in $\Exp(H_i)$ is equal to the sum of the coefficient at $w^\gamma z^d$ in $H_i$ and a term depending only on coefficients at $w^{\gamma'}z^{d'}$ where $\gamma'$ is a summand of $\gamma$.
\end{proof}

We return to the proof of Proposition~\ref{pr:ExplicitFull-}. Applying the plethystic logarithm to the identity of Proposition~\ref{pr:NilpPow}, we obtain in $\cMot(\kk)[[\Gamma_D,z]]$:
\[
    \Log[\Pair^-(\divisor)\times_{\Bun^{par}(\divisor)}\Bun^{par}(\divisor,\zeta)]=
    \bL\Log[\Pair^{nilp,-}(\divisor)\times_{\Bun^{par}(\divisor)}\Bun^{par}(\divisor,\zeta)].
\]
Consider the identity of Proposition~\ref{pr:SlopeProd}. Applying the homomorphism $\Mot(\kk)\to\cMot(\kk)$ and the plethystic logarithm to this identity and restricting to sum of the monomials $A_{\gamma,d}w^\gamma z^d$ with $-\epsilon d+\gamma\star\zeta=0$ (that is, $\tau=0$), we obtain
\begin{multline}\label{eq:Logs}
    \Log\left(\sum_{\substack{\gamma\in\Gamma_D,d\le0\\ \epsilon d+\gamma\star\zeta=0}}
    (-1)^{\delta\rk\gamma} [\Pair^-(\divisor)\times_{\Bun^{par}(\divisor)}\Bun_{\gamma,d}^{par,(\epsilon,\zeta)-iso}(\divisor,\zeta)]w^\gamma z^{-d}\right)\\=
    \Bigl(\Log([\Pair^-(\divisor)\times_{\Bun^{par}(\divisor)}\Bun^{par}(\divisor,\zeta)])\Bigr)_{-\epsilon d+\gamma\star\zeta=0}\\=
    \bL\Bigl(\Log(\Pair^{nilp,-}(\divisor)\times_{\Bun^{par}(\divisor)}\Bun^{par}(\divisor,\zeta))\Bigr)_{-\epsilon d+\gamma\star\zeta=0}.
\end{multline}
Consider the series $\Omega^{Sch,mot}_{X,D,\delta}(z,\bL^{-\delta/2}w,w_{\bullet,\bullet})$ obtained by substituting $w$ with $\bL^{-\delta/2}w$, so that its coefficient at $w^\gamma z^d$ differs from the coefficient of $\Omega^{Sch,mot}_{X,D,\delta}$ at $w^\gamma z^d$ by the factor of $(\bL^\frac12)^{-\delta\rk\gamma}$. It follows from Lemma~\ref{lm:ExplPair} that the coefficients at $w^\gamma z^{-d}$ of the motivic series
\[
    \Omega^{Sch,mot}_{X,D,\delta}(z,\bL^{-\delta/2}w,w_{\bullet,\bullet})\text{ and }
    [\Pair^{nilp,-}(\divisor)\times_{\Bun^{par}(\divisor)}\Bun^{par}(\divisor,\zeta)]
\]
coincide whenever $\gamma\in\Gamma_{\divisor,\zeta}$ (recall the presence of $(-1)^{\delta\rk\gamma}$ in convention~\eqref{eq:gradedStacks}). Noting that the~substitution $w\mapsto \bL^{\delta/2}w$ commutes with Log, and applying Lemma~\ref{lm:GoodExp} and~\eqref{eq:Logs}, we see that the coefficients at $w^{\gamma'}z^{-d'}$ of the series
\begin{equation*}
    \Bigl(\bL\Log\Omega^{Sch,mot}_{X,D,\delta}(z,\bL^{-\delta/2}w,w_{\bullet,\bullet})\Bigr)_{-\epsilon d+\gamma\star\zeta=0}
\end{equation*}
and
\[
    \Log\left(\sum_{\substack{\gamma\in\Gamma_D,d\le0\\ \epsilon d+\gamma\star\zeta=0}}
    (-1)^{\delta\rk\gamma}
    [\Pair^-(\divisor)\times_{\Bun^{par}(\divisor)}\Bun_{\gamma,d}^{par,(\epsilon,\zeta)-iso}(\divisor,\zeta)]w^\gamma z^{-d}\right)
\]
coincide whenever $\gamma'\in\Gamma_{\divisor,\zeta}$. Exponentiating and using Lemma~\ref{lm:RR}(i) and Lemma~\ref{lm:GoodExp}, we see that for $\gamma'\in\Gamma_{\divisor,\zeta}$ and $d'\le0$, the motivic class of $\Conn^-_{\gamma',d'}(\divisor,\zeta)$ is equal to the coefficient at $w^{\gamma'}z^{-d'}$ in
\[
    (-1)^{\delta\rk\gamma}\bL^{\beta(\gamma')}
    \Exp\Bigl(\bL\Log\Omega^{Sch,mot}_{X,D,\delta}(z,\bL^{-\delta/2}w,w_{\bullet,\bullet})\Bigr)_{-\epsilon d+\gamma\star\zeta=0}.
\]
Since $2\beta(\gamma')-\delta\rk\gamma'=\chi(\gamma')$ and the substitution $w\mapsto \bL^{-\delta/2}w$ commutes with $\Exp$, we are done with the proof of Proposition~\ref{pr:ExplicitFull-}. \qed

\section{Semistability of parabolic \texorpdfstring{$\epsilon$}{epsilon}-connections}\label{sect:Stability}
Recall that we have divisors $\divisor'=\sum_{x\in D'}n'_xx\le\divisor=\sum_{x\in D}{n_xx}$ on the $\kk$-curve $X$, and we are studying the moduli stacks of parabolic $\epsilon$-connections with partially or fully fixed normal forms denoted $\Conn_{\gamma,d}^{prtl}(\epsilon,\divisor,\divisor',\zeta)$, where $(\gamma,d)\in\Gamma_{D'}\times\Z$. These stacks are of infinite type (and, in fact, of infinite motivic volume) unless $\divisor=\divisor'$ and $\epsilon\ne0$. To obtain substacks of finite type, we will add stability conditions. The relation between the motivic classes of parabolic $\epsilon$-connections with HN-nonpositive underlying vector bundles and the motivic classes of semistable loci of these stacks is given by a~version of Kontsevich--Soibelman factorization formula, see Proposition~\ref{pr:KontsevichSoibelman}. The explicit formulas in terms of the irregular Schiffmann's generating function $\Omega^{Sch,mot}_{X,D,\delta}$ are provided by Corollary~\ref{cor:ExplNonnegativeSemistable}. To derive these formulas, we will study the homological algebra of level $\divisor'$ parabolic bundles and of parabolic bundles with $\epsilon$-connections; see Section~\ref{sect:Homological}. These results are of independent interest. We are still working with parabolic $\epsilon$-connections whose underlying vector bundle is nonpositive, which results in a subtle notion of nonpositive-semistable parabolic $\epsilon$-connections, see Definition~\ref{def:Semistable}(ii) and~\cite[Rem.~3.3.1]{FedorovSoibelmans}.

\subsection{Parabolic subbundles and quotient bundles} Recall that if $E$ is a vector bundle on $X$, its subbundle $F$ is called {\it strict}, if $E/F$ is torsion free. The new difficulty in the irregular case is that, given a parabolic $\epsilon$-connection $(E,E_{\bullet,\bullet},\nabla)$ and a strict vector subbundle $F\subset E$ preserved by $\nabla$, the intersections $F\cap E_{\bullet,\bullet}$ do not necessarily give rise to a parabolic structure on $F$. Indeed, it might happen that $F\cap E_{x,j}$ is not a free submodule of $E_{x,j}$ as the following example shows.

\begin{example}\label{ex:DoesNotInduce}
  As in Example~\ref{ex:IrregNotSuff}, let $X=\P^1_k$ and $z$ be the standard coordinate. Take $\divisor=\divisor'=2(0)$, and let $E:=\cO\oplus\cO$ be the trivial rank~2 vector bundle. Let the submodule $E_1\subset E_{2(0)}$ be generated by $e:=(1,z)$. Define a level~2 parabolic structure on $E$ at~0 as follows: $E_{0,0}=E_{2(0)}$, $E_{0,1}=E_1$, $E_{0,j}=0$ for $j\ge2$. Take $\nabla=d$. Of course, $\nabla$ is a connection with poles bounded by $\divisor$. It preserves any parabolic structure at 0, in particular, it preserves $E_{0,\bullet}$. On the other hand, consider the subbundle $F=\cO\oplus0\subset E$. Clearly, this subbundle is preserved by $\nabla$. On the other hand, $F\cap E_{0,1}$ is the submodule of $E_{2(0)}$ generated by $(z,0)$, so it is not free.
\end{example}

The situation improves when the $\epsilon$-connection is generic as the following lemma shows (see Section~\ref{sect:ParStrPreserved} for the definition of parabolic structures preserved by a connection).

\begin{lemma}\label{lm:Subbundle}
Let $(E,\nabla)$ be an $\epsilon$-connection on $X$ with poles bounded by $\divisor$. Let $x\in D$. Assume that~$E_{x,\bullet}$ is a level $n'_x$ parabolic structure at $x$ preserved by $\nabla$, where $n'_x\le n_x$. Assume that~$\nabla$ is generic at $x$ or $n'_x=1$. Then for any strict subbundle $F\subset E$ preserved by $\nabla$, the intersections $F\cap E_{x,\bullet}$ give rise to a level $n'_x$ parabolic structure on $F$ at $x$ preserved by $\nabla|_F$, while $E_{x,\bullet}/(F\cap E_{x,\bullet})$ is a level $n'_x$ parabolic structure on $E/F$ preserved by the $\epsilon$-connection induced by $\nabla$.
\end{lemma}
\begin{proof}
  The case $n'_x=1$ is obvious (cf.~\cite[Sect.~3.5]{FedorovSoibelmansParabolic}). Thus we may assume that $n'_x\ge2$ so that~$\nabla$ is generic. Let $r=\rk E$ and $e_1,\ldots,e_r$ be a basis of $E_{n'_xx}$ that diagonalizes $\nabla$ as in Lemma~\ref{lm:diag}. Let us call a $\kk[n'_xx]$-submodule of $E_{n'_xx}$ {\it standard}, if it is generated by a subset of $\{e_1,\ldots,e_r\}$. Note that each $E_{x,j}$ and $F_{n_xx}$ are preserved by $\nabla$, see Section~\ref{sect:ParStrPreserved}. By Lemma~\ref{lm:NablaPreservedSubmodules}, all $E_{x,j}$ and $F_{n'_xx}$ are standard. It follows that $F\cap E_{x,j}$ are standard submodules of the free module $E_{n'_xx}$, so, in particular, they are free. Thus, $F\cap E_{x,\bullet}$ is a level $n'_x$ parabolic structure on $F$, while $E_{x,\bullet}/(F\cap E_{x,\bullet})$ is a level $n'_x$ parabolic structure on $E/F$. Clearly, they are preserved by the respective induced $\epsilon$-connections.
\end{proof}

For $\gamma=(r,r_{\bullet,\bullet}),\gamma'=(r',r'_{\bullet,\bullet})\in\Gamma_{D'}$, we write $\gamma'\le\gamma$, if $r'_{x,j}\le r_{x,j}$ for all $x$ and $j$. Equivalently, $\gamma-\gamma'\in\Gamma_{D'}$. Recall that when $\zeta\in\FNF(\divisor,\divisor')$, the set $\Gamma_{\divisor',\zeta}$ consists of all classes $\gamma\in\Gamma_{D'}$ such that~$\gamma$ is full at all $x\in D'$ with $n'_x\ge2$ and $\zeta$ is non-resonant for $\gamma$ at all such $x$.

\begin{corollary}\label{cor:Subbundle}
Let $(E,E_{\bullet,\bullet},\nabla)\in\Conn^{prtl}_{\gamma,d}(\divisor,\divisor',\zeta)(\kk)$. Assume that $\gamma\in\Gamma_{\divisor',\zeta}$. Then for any strict subbundle $F\subset E$ preserved by $\nabla$, we have

(i) $(F,F\cap E_{\bullet,\bullet},\nabla|_F)\in\Conn^{prtl}_{\gamma',d'}(\divisor,\divisor',\zeta)(\kk)$ for some $\gamma'\le\gamma$ and $d'\in\Z$;

(ii) $\gamma',\gamma-\gamma'\in\Gamma_{\divisor',\zeta}$;

(iii) $(E/F,E_{\bullet,\bullet}/(F\cap E_{\bullet,\bullet}),\nabla|_{E/F})\in\Conn^{prtl}_{\gamma-\gamma',d-d'}(\divisor,\divisor',\zeta)(\kk)$, where $\nabla|_{E/F}$ is the $\epsilon$-connection induced by $\nabla$ on $E/F$.
\end{corollary}
\begin{proof}
  By Lemma~\ref{lm:DiagIffNonResonant}, $\nabla$ is generic at all $x$ with $n'_x\ge2$. Thus, by Lemma~\ref{lm:Subbundle}, $F\cap E_{\bullet,\bullet}$ and $E_{\bullet,\bullet}/(F\cap E_{\bullet,\bullet})$ are level $\divisor'$ parabolic structure on $F$ and $E/F$ respectively, preserved by the induced $\epsilon$-connections.
  We also note that, by definition of $\Conn^{prtl}_{\gamma,d}(\divisor,\divisor',\zeta)$, the connection $\nabla$ induces multiplication by $\zeta_{x,j}$ on $E_{x,j-1}/E_{x,j}$. Thus, it induces multiplication by $\zeta_{x,j}$ on $(F\cap E_{x,j-1})/(F\cap E_{x,j})$ to the effect that $(F,F\cap E_{\bullet,\bullet},\nabla|_F)\in\Conn^{prtl}(\divisor,\divisor',\zeta)(\kk)$. Similarly, $(E/F,E_{\bullet,\bullet}/(F\cap E_{\bullet,\bullet}),\nabla|_{E/F})\in\Conn^{prtl}(\divisor,\divisor',\zeta)(\kk)$. The rest of the corollary follows from the definitions.
\end{proof}

\subsection{Definition of Semistability}\label{sect:DefStab} The definition below is a version of~\cite[Def.~3.8]{FedorovSoibelmansParabolic}.
\begin{definition}\label{def:ParWeights}
We say that a sequence $\sigma=\sigma_{\bullet,\bullet}$ of real numbers indexed by $D'\times\Z_{>0}$ is a \emph{sequence of parabolic weights of type $(X,\divisor')$,\/} if for all $x\in D'$ we have
\begin{equation}\label{eq:StabCond}
    \sigma_{x,1}\le\sigma_{x,2}\le\ldots
\end{equation}
and for all $x\in D'$ and all $j>0$ we have $\sigma_{x,j}\le\sigma_{x,1}+n'_x$.
\end{definition}

Similarly to~\cite[Sect.~1.2.2]{FedorovSoibelmansParabolic}, for a sequence $\sigma$ of parabolic weights, we define the \emph{$\sigma$-degree\/} of a~parabolic bundle $\bE\in\Bun_{\gamma,d}^{par}(\divisor')(\kk)$ by
\[
    \deg_{1,\sigma}\bE:=d+\gamma\star\sigma=d+\sum_{x\in D'}\sum_{j>0}r_{x,j}\sigma_{x,j}\in\R,
\]
where $\gamma=(r,r_{\bullet,\bullet})$. We define the \emph{$\sigma$-slope of $\bE$\/} as $\deg_{1,\sigma}\bE/\rk\bE$. We remark that, more generally, for $\kappa\in\R$, we can define the $(\kappa,\sigma)$-degree by $\deg_{\kappa,\sigma}\bE:=\kappa d+\gamma\star\sigma$ but we will not need this except for $\kappa=1$. The following lemma explains why we need the above assumptions on $\sigma_{x,j}$ (cf.~\cite[Lemma~3.12]{FedorovSoibelmansParabolic}).
\begin{lemma}\label{lm:GenericIso}
  Let $\bE\to\bF$ be a morphism in $\Bun_{\gamma,d}^{par}(\divisor')(\kk)$ such that the underlying morphism of vector bundles is generically an isomorphism. Then $\deg_{1,\sigma}\bE\le\deg_{1,\sigma}\bF$.
\end{lemma}
\begin{proof}
  Write $\bE=(E,E_{\bullet,\bullet})$ and $\bF=(F,F_{\bullet,\bullet})$. Consider a morphism $\bE\to\bF$ with the underlying morphism of vector bundles $f\colon E\to F$. Assume that $f$ is generically an isomorphism. For $x\in D'$, let $f_x\colon E_{n'_xx}\to F_{n'_xx}$ be the induced morphism of restrictions to infinitesimal neighborhoods of $x$. Let $d_x$ be length of the kernel of $f_x$. It is equal to the length of the cokernel of $f_x$. Since the cokernel of $f$ is a torsion sheaf, we obtain $\deg E\le\deg F-\sum_{x\in D'}d_x$. On the other hand,
  \[
    d_x\ge\dim_\kk\Ker(f_x|_{E_{x,j}})\ge n'_x\rk E_{x,j}-n'_x\rk F_{x,j}
  \]
  for all $j>0$. Thus,
  \begin{multline*}
        \deg_{1,\sigma}\bE=\deg E+\sum_{x\in D',j>0}\sigma_{x,j}(\rk E_{x,j-1}-\rk E_{x,j})=\\
        \deg E+\sum_{x\in D'}\left(\sigma_{x,1}\rk E+\sum_{i>0}(\sigma_{x,i+1}-\sigma_{x,i})\rk E_{x,i}\right)\le\\
        \deg F-\sum_{x\in D'} d_x+\sum_{x\in D'}\left(\sigma_{x,1}\rk F+\sum_{i>0}(\sigma_{x,i+1}-\sigma_{x,i})
        \left(\rk F_{x,i}+\frac{d_x}{n'_x}\right)\right)=\\
        \deg_{1,\sigma}\bF+\sum_{x\in D'}d_x\left(-1+\frac1{n'_x}\sum_{i>0}(\sigma_{x,i+1}-\sigma_{x,i})\right)\le\deg_{1,\sigma}\bF.
  \end{multline*}
\end{proof}

The following definition is similar to~\cite[Def.~6.1]{FedorovSoibelmansParabolic}.
\begin{definition}\label{def:Semistable}
Let $(E,E_{\bullet,\bullet},\nabla)\in\Conn^{prtl}_{\gamma,d}(\divisor,\divisor',\zeta)(\kk)$. Assume that $\gamma\in\Gamma_{\divisor',\zeta}$.

\stepzero\noindstep We say that $(E,E_{\bullet,\bullet},\nabla)$ is \emph{$\sigma$-semistable\/} if for all strict subbundles $F$ preserved by $\nabla$ we have
\begin{equation}\label{eq:Semistable}
    \frac{\deg_{1,\sigma}(F,F\cap E_{\bullet,\bullet})}{\rk F}\le\frac{\deg_{1,\sigma}(E,E_{\bullet,\bullet})}{\rk E}.
\end{equation}

\noindstep We say that $(E,E_{\bullet,\bullet},\nabla)$ is \emph{$\sigma$-nonpositive-semistable}, if $E$ is nonpositive and~\eqref{eq:Semistable} is satisfied for all strict subbundles $F$ such that $\nabla$ preserves $F$ and $E/F$ is a nonpositive vector bundle.
\end{definition}

Using an argument similar to~\cite[Lemma~3.7]{Simpson1}, we see that these conditions are open, so we get open substacks
\[
\Conn^{prtl,\sigma-ss,-}_{\gamma,d}(\divisor,\divisor',\zeta)\subset\Conn_{\gamma,d}^{prtl,-}(\divisor,\divisor',\zeta)\text{ and }
\Conn_{\gamma,d}^{prtl,\sigma-ss}(\divisor,\divisor',\zeta)\subset\Conn_{\gamma,d}^{prtl}(\divisor,\divisor',\zeta),
\]
whenever $\gamma\in\Gamma_{\divisor',\zeta}$.

\begin{lemma}\label{lm:MorphismSS}
Let $\gamma_1,\gamma_2\in\Gamma_{\divisor',\zeta}$ and $d_1,d_2\in\Z$.

  (i) Let $(\bE_i,\nabla_i)\in\Conn_{\gamma_i,d_i}^{prtl,\sigma-ss}(\divisor,\divisor',\zeta)(\kk)$ be parabolic $\epsilon$-connections $(i=1,2)$ such that the $\sigma$-slope of $\bE_1$ is greater than that of $\bE_2$. Then $\Hom((\bE_1,\nabla_1),(\bE_2,\nabla_2))=0$.

  (ii) A similar statement holds for $(\bE_i,\nabla_i)\in\Conn_{\gamma_i,d_i}^{prtl,\sigma-ss,-}(\divisor,\divisor',\zeta)(\kk)$.
\end{lemma}
\begin{proof}
  (i) Let $E_i$ be the underlying vector bundle of $\bE_i$ ($i=1,2$). Assume that $f\colon(\bE_1,\nabla_1)\to(\bE_2,\nabla_2)$ is a non-zero morphism. Denote the morphism of underling vector bundles by $f$ as well. Let $F:=\Ker f$ and $F'$ be the saturation of the image of $f$ (that is, $F'$ is the unique strict subbundle of $E_2$ such that $f(E_1)\subset F'$ and $f(E_1)$ coincides with $F'$ generically). Since $F$ and $F'$ are preserved by $\nabla$, we can use Lemma~\ref{lm:Subbundle} to promote them to parabolic $\epsilon$-connections
  $(\bF,\nabla),(\bF',\nabla')\in\Conn^{prtl}(\divisor,\divisor',\zeta)(\kk)$. Then $f$ induces a morphism $\bE_1/\bF\to\bF'$, which is generically an isomorphism, so we can apply Lemma~\ref{lm:GenericIso}. Combining this with semistability of $(\bE_i,\nabla_i)$, we get the following contradiction:
  \[
    \frac{\deg_{1,\sigma}\bE_1}{\rk\bE_1}\le
    \frac{\deg_{1,\sigma}\bE_1/\bF}{\rk(\bE_1/\bF)}\le
    \frac{\deg_{1,\sigma}\bF'}{\rk\bF'}\le
    \frac{\deg_{1,\sigma}\bE_2}{\rk\bE_2}.
  \]
  (ii) Keeping the notation from part~(i), we note that $F$ and $F'$ are nonpositive, as they are subbundles of nonpositive vector bundles $E_1$ and $E_2$ respectively. The rest of the proof is the same as in part~(i).
\end{proof}

\subsection{Harder--Narasimhan filtrations} In Proposition~\ref{pr:KontsevichSoibelman} we will prove a version of Kontsevich--Soibelman factorization formula. The proof of this formula uses Harder--Narasimhan filtrations.

Let $(E,E_{\bullet,\bullet},\nabla)\in\Conn^{prtl}_{\gamma,d}(\divisor,\divisor',\zeta)(\kk)$. Assume that $\gamma\in\Gamma_{\divisor',\zeta}$. Let $0=E_0\subset E_1\subset\ldots\subset E_m=E$ be a filtration by strict subbundles preserved by $\nabla$. For all $i$ we get an induced $\epsilon$-connection on $E_i/E_{i-1}$, which we denote by $\nabla|_{E_i/E_{i-1}}$. By Corollary~\ref{cor:Subbundle}, we also get an induced parabolic structure $(E_i\cap E_{\bullet,\bullet})/(E_{i-1}\cap E_{\bullet,\bullet})$ on $E_i/E_{i-1}$ preserved by the connection.
\begin{proposition}\label{pr:HarderNarasimhan}
(i) Let $(E,E_{\bullet,\bullet},\nabla)\in\Conn^{prtl}_{\gamma,d}(\divisor,\divisor',\zeta)(\kk)$. Assume that $\gamma\in\Gamma_{\divisor',\zeta}$. Then there is a unique filtration $0=E_0\subset E_1\subset\ldots\subset E_m=E$ by strict subbundles preserved by $\nabla$ such that all the quotients
\begin{equation}\label{eq:quotients}
    (E_i/E_{i-1},(E_i\cap E_{\bullet,\bullet})/(E_{i-1}\cap E_{\bullet,\bullet}),\nabla|_{E_i/E_{i-1}})
\end{equation}
are $\sigma$-semistable, and the sequence of slopes
\begin{equation}\label{eq:slopes}
    \tau_i:=\frac{\deg_{1,\sigma}(E_i/E_{i-1},(E_i\cap E_{\bullet,\bullet})/(E_{i-1}\cap E_{\bullet,\bullet}))}
    {\rk(E_i/E_{i-1})}
\end{equation}
is strictly decreasing.

(ii) Let $(E,E_{\bullet,\bullet},\nabla)\in\Conn^{prtl,-}_{\gamma,d}(\divisor,\divisor',\zeta)(\kk)$. Assume that $\gamma\in\Gamma_{\divisor',\zeta}$. Then there is a unique filtration $0=E_0\subset E_1\subset\ldots\subset E_m=E$ by strict subbundles preserved by $\nabla$ such that all the quotients~\eqref{eq:quotients} are $\sigma$-nonpositive-semistable, and the sequence of slopes~\eqref{eq:slopes} is strictly decreasing.
\end{proposition}
We will call the filtration $0=\bE_0\subset\bE_1\subset\ldots\subset\bE_m=\bE$, where $\bE_i:=(E_i,E_i\cap E_{\bullet,\bullet})$, the \emph{Harder--Narasimhan filtration} and the decreasing sequence $\tau_i$ as above the \emph{Harder--Narasimhan spectrum of $(E,E_{\bullet,\bullet},\nabla)$}. We note that a $\sigma$-nonpositive-semistable parabolic bundle is not necessarily $\sigma$-semistable, so the filtrations in parts~(i) and~(ii) are in general different even if $E$ is nonpositive, see~\cite[Rem.~3.3.1]{FedorovSoibelmans}.
\begin{proof}
  The proof is similar to that of~\cite[Prop.~3.11, Prop.~6.2]{FedorovSoibelmansParabolic}; we sketch it for the reader's convenience. We will only prove (ii) as (i) is completely similar. We proceed by induction on $\rk E$. Let $F\subset E$ be a strict vector subbundle (automatically nonpositive) preserved by $\nabla$ and such that $E/F$ is nonpositive. Then by Corollary~\ref{cor:Subbundle} there is $\gamma'\le\gamma$ such that $(F,F\cap E_{\bullet,\bullet},\nabla|_F)\in\Conn^{prtl,-}_{\gamma',d'}(\divisor,\divisor',\zeta)(\kk)$ and $(E/F,E_{\bullet,\bullet}/(F\cap E_{\bullet,\bullet}),\nabla|_{E/F})\in\Conn^{prtl,-}_{\gamma-\gamma',d-d'}(\divisor,\divisor',\zeta)(\kk)$. Note also that $\gamma',\gamma-\gamma'\in\Gamma_{\divisor',\zeta}$. We choose $F$ such that $\deg_{1,\zeta}(F,F\cap E_{\bullet,\bullet})$ is maximal. It is clear that $(F,F\cap E_{\bullet,\bullet},\nabla|_F)$ is nonpositive-semistable. We apply the induction hypothesis to $(E/F,E_{\bullet,\bullet}/(F\cap E_{\bullet,\bullet}),\nabla|_{E/F})$. The proof is completed in a standard way using Lemma~\ref{lm:MorphismSS}(ii) (see~\cite[Sect.~1.3]{HarderNarasimhan}).
\end{proof}

\subsection{Homological algebra of parabolic bundles and parabolic $\epsilon$-connections}\label{sect:Homological} The category $\Bun^{par}(\divisor')(\kk)$ is an exact category in the sense of Quillen. Precisely, a sequence
\[
    0\to(E^1,E^1_{\bullet,\bullet})\to(E^2,E^2_{\bullet,\bullet})\to(E^3,E^3_{\bullet,\bullet})\to0
\]
is exact if the corresponding sequence of vector bundles $0\to E^1\to E^2\to E^3\to0$ is exact and for all $x\in D'$ and $j>0$ the induced sequence of $\kk[n'_xx]$-modules
\[
    0\to E^1_{x,j}\to E^2_{x,j}\to E^3_{x,j}\to0
\]
is exact. In principle, one can embed $\Bun^{par}(\divisor')(\kk)$ into an abelian category of parabolic sheaves and talk about resolutions as in~\cite{GothenKing05}. However, there are certain problems with this approach. For example, the category of $\kk[nx]$-modules has infinite homological dimension whenever $n>1$. To avoid these difficulties, we take an ad hoc approach. We note that our results are similar to that of~\cite{GothenKing05} when $\divisor'=0$ and $\epsilon=0$.

\begin{lemma}\label{lm:Ext=H1}
  Let $\bE,\bF\in\Bun^{par}(\divisor')(\kk)$ be parabolic bundles and let $\HOM(\bE,\bF)$ be the sheaf of homomorphisms. Then the vector space $H^1(X,\HOM(\bE,\bF))$ classifies extensions of $\bE$ by $\bF$.
\end{lemma}
\begin{proof}
  We claim that every extension $0\to\bF\to\bG\to\bE\to0$ is split locally in the Zariski topology. Indeed, it is enough to check the statement in an open neighborhood of $x\in D'$. Write $\bE=(E,E_{\bullet,\bullet}),\bF=(F,F_{\bullet,\bullet}),\bG=(G,G_{\bullet,\bullet})$. Since the $\kk[n'_xx]$-module $E_{x,j}$ is free for all $j$, we can split the exact sequence $0\to F_{\kk[n'_xx]}\to G_{\kk[n'_xx]}\to E_{\kk[n'_xx]}\to0$ in a way compatible with the filtrations (consider a basis of $E_{\kk[n'_xx]}$ compatible with the filtration and lift this basis to $G_{\kk[n'_xx]}$). We can extend this splitting to a Zariski open neighborhood of $x$. Now, using the \v{C}ech description of $H^1(X,\HOM(\bE,\bF))$, we obtain the statement.
\end{proof}

In view of the previous lemma we write $\Ext^1(\bE,\bF)$ for $H^1(X,\HOM(\bE,\bF))$.
\begin{lemma}\label{lm:ExactSequence0}
  Let $\bE=(E,E_{\bullet,\bullet}),\bF=(F,F_{\bullet,\bullet})\in\Bun^{par}(\divisor')(\kk)$ be parabolic bundles. Then we have a canonical exact sequence
  \begin{multline}\label{eq:ExtHomExact}
    0\to\Hom(\bE,\bF)\to\Hom(E,F)\to
    \bigoplus_{x\in D'}\Hom_{\kk[n'_xx]}(E_{n'_xx},F_{n'_xx})/\Hom_{\kk[n'_xx]}(E_{x,\bullet},F_{x,\bullet})
    \\ \to\Ext^1(\bE,\bF)\to\Ext^1(E,F)\to0.
  \end{multline}
\end{lemma}
\begin{proof}
    The sheaf of homomorphisms $\HOM(\bE,\bF)$ fits into the exact sequence
    \[
        0\to\HOM(\bE,\bF)\to\HOM(E,F)\to\bigoplus_{x\in D'}\Hom_{\kk[n'_xx]}(E_{n'_xx},F_{n'_xx})/\Hom_{\kk[n'_xx]}(E_{x,\bullet},F_{x,\bullet})\to0,
    \]
    where we view the rightmost term as a torsion sheaf supported on $\divisor'$. It remains to consider the corresponding long exact sequence of cohomology groups and use Lemma~\ref{lm:Ext=H1}.
\end{proof}

\begin{corollary}\label{cor:EulerChar0}
   Keeping the notations of the lemma, let us write $\cl(\bE)=(r^\bE,r^\bE_{\bullet,\bullet},d^\bE)$ and $\cl(\bF)=(r^\bF,r^\bF_{\bullet,\bullet},d^\bF)$. Then we have
  \[
      \dim\Hom(\bE,\bF)-\dim\Ext^1(\bE,\bF)=\chi(\HOM(\bE,\bF))=(1-g)r^\bE r^\bF+r^\bE d^\bF-r^\bF d^\bE-\sum_{\substack{x\in D'\\i<j}}n'_xr^\bE_{x,i}r^\bF_{x,j}.
  \]
\end{corollary}
\begin{proof}
  Apply the Euler characteristic to the exact sequence~\eqref{eq:ExtHomExact} and note that
  \[
    \dim\left(\Hom_{\kk[n'_xx]}(E_{n'_xx},F_{n'_xx})/\Hom_{\kk[n'_xx]}(E_{x,\bullet},F_{x,\bullet})\right)= n'_x\sum_{i<j}r^\bE_{x,i}r^\bF_{x,j}.
  \]
\end{proof}

Let $\HOM_0(\bE,\bF)\subset\HOM(\bE,\bF)$ be the subsheaf of homomorphisms $\Psi$ such that for all $x$ and $j$ we have $\Psi(E_{x,j})\subset F_{x,j+1}$, where $\bE=(E,E_{\bullet,\bullet})$ and $\bF=(F,F_{\bullet,\bullet})$ (in other words, these are homomorphisms inducing zero on the associated graded, cf.~Section~\ref{sect:AlmostAffine}). Let $\Hom_0(\bE,\bF)=H^0(X,\HOM_0(\bE,\bF))$ be the corresponding space of homomorphisms. Similarly to the proof of Lemma~\ref{lm:RR-dim}, we have a non-degenerate $\cO_X$-bilinear trace pairing $\HOM(\bE,\bF)\otimes\HOM_0(\bF,\bE)\to\cO_X(-\divisor')$. Applying Serre duality to $\Ext^1(\bE,\bF)=H^1(X,\HOM(\bE,\bF))$ (see Lemma~\ref{lm:Ext=H1}), we get the following lemma.
\begin{lemma}\label{lm:SerreDuality}
The vector space $\Ext^1(\bE,\bF)$ is dual to $\Hom_0(\bF,\bE)\otimes\Omega_X(\divisor')$.
\end{lemma}

Similarly, the category $\Conn^{prtl}(\divisor,\divisor',\zeta)$ is an exact category. Consider parabolic $\epsilon$-connections $(\bE,\nabla_\bE),(\bF,\nabla_\bF)\in\Conn^{prtl}(\divisor,\divisor',\zeta)(\kk)$, and the sheaf of homomorphisms $\HOM(\bE,\bF)$. It follows from Definition~\ref{def:ModSpacePartial} that $\nabla_\bE$ and $\nabla_\bF$ induce an $\epsilon$-connection
\[
    \HOM(\bE,\bF)\xrightarrow{\nabla}\HOM_0(\bE,\bF)\otimes\Omega_X(\divisor).
\]
We denote the above 2-term complex of sheaves of vector spaces by $\HOM^\bullet((\bE,\nabla_\bE),(\bF,\nabla_\bF))$. We normalize the complex so it is concentrated in cohomological degrees 0 and 1. Let $\Ext^i((\bE,\nabla_\bE),(\bF,\nabla_\bF))$ be the $i$-th hypercohomology group of this complex. It is clear that $\Ext^0((\bE,\nabla_\bE),(\bF,\nabla_\bF))=\Hom((\bE,\nabla_\bE),(\bF,\nabla_\bF))$.

\begin{lemma}\label{lm:ExactSequence}
  Let $(\bE,\nabla_\bE),(\bF,\nabla_\bF)\in\Conn^{prtl}(\divisor,\divisor',\zeta)(\kk)$ be parabolic $\epsilon$-connections. Then we have an exact sequence
  \begin{multline}\label{eq:ExactSequence}
    0\to\Hom((\bE,\nabla_\bE),(\bF,\nabla_\bF))\\ \to
    \Hom_0(\bE,\bF\otimes\Omega_X(\divisor))\to
    \Ext^1((\bE,\nabla_\bE),(\bF,\nabla_\bF))\to\Ext^1(\bE,\bF)\\
    \to H^1(X,\HOM_0(\bE,\bF\otimes\Omega_X(\divisor)))\to\Ext^2((\bE,\nabla_\bE),(\bF,\nabla_\bF))\to0.
  \end{multline}
\end{lemma}
\begin{proof}
    We have an exact triangle
  \[
    \HOM_0(\bE,\bF\otimes\Omega_X(\divisor))[-1]\to\HOM((\bE,\nabla_\bE),(\bF,\nabla_\bF))\to\HOM(\bE,\bF)
    \to\HOM_0(\bE,\bF\otimes\Omega_X(\divisor)).
  \]
  Applying the cohomology functor, we get the required statement. (Note that the $i$-th cohomology group of any coherent sheaf vanishes for $i\ge2$ because $\dim X=1$.)
\end{proof}

Let $(\bE,\nabla_\bE)\in\Conn^{prtl}(\divisor,\divisor',\zeta)(\kk)$. Let $\divisor''$ be a divisor on $X$ whose set-theoretic support is contained in $D$. Since $E|_{X-D}=E(\divisor'')|_{X-D}$, the $\epsilon$-connection $\nabla$ gives rise to a meromorphic $\epsilon$-connection~$\nabla''$ on $E(\divisor'')$. It is easy to see that $\nabla''$ is an $\epsilon$-connection on $E(\divisor'')$ with poles bounded by $\divisor$. We warn the reader that $(\bE(\divisor''),\nabla'')\in\Conn^{prtl}(\divisor,\divisor',\zeta')(\kk)$, where in general $\zeta'\ne\zeta$. We send the reader to Section~\ref{sect:TwistConn} for details. For now, we will take $\divisor''=\divisor'-\divisor$ and abuse the notation by denoting the induced connection~$\nabla''$ simply by $\nabla_\bE$, so that $(\bE(\divisor'-\divisor),\nabla_\bE)$ is a parabolic $\epsilon$-connection with poles bounded by $\divisor$.
\begin{lemma}\label{lm:Ext}
Let $(\bE,\nabla_\bE),(\bF,\nabla_\bF)\in\Conn^{prtl}(\divisor,\divisor',\zeta)(\kk)$ be parabolic $\epsilon$-connections.

  (i) Write $\cl(\bE)=(\gamma_1,d_1)$ and $\cl(\bF)=(\gamma_2,d_2)$. Assume that $\gamma_1+\gamma_2\in\Gamma_{\divisor',\zeta}$. Then the vector space $\Ext^1((\bE,\nabla_\bE),(\bF,\nabla_\bF))$ classifies extensions of $(\bE,\nabla_\bE)$ by $(\bF,\nabla_\bF)$.

  (ii) The vector space $\Ext^2((\bE,\nabla_\bE),(\bF,\nabla_\bF))$ is dual to $\Hom((\bF,\nabla_\bF),(\bE(\divisor'-\divisor),\nabla_\bE))$.
\end{lemma}
\begin{proof}
  For (i), we note that similarly to the proof of Lemma~\ref{lm:ExactSequence}, we can show that every exact sequence $0\to(\bF,\nabla_\bF)\to(\bG,\nabla_\bG)\to(\bE,\nabla_\bE)\to0$ is split in a formal neighborhood of $x$. Indeed, write $\bG=(G,G_{\bullet,\bullet})$. Since $\cl(\bG)=(\gamma_1+\gamma_2,d_1+d_2)$ and $\gamma_1+\gamma_2\in\Gamma_{\divisor',\zeta}$, the $\epsilon$-connection $\nabla_\bG$ is generic at $x$ by Lemma~\ref{lm:DiagIffNonResonant}. By Lemma~\ref{lm:diag}, we can find a formal trivialization of $G$ at $x$ compatible with $G_{x,\bullet}$ and diagonalizing $\nabla_\bG$. Now the claim follows from Lemma~\ref{lm:NablaPreservedSubmodules}. The union of the formal neighborhoods of point of $D$ with $X-D$ is an fppf cover of $X$, so we can use it to calculate the hypercohomology (cf.~\cite[Tag03DW]{StacksProject}). The statement follows.

  For (ii) note that, similarly to the proof of Lemma~\ref{lm:RR-dim}, the locally free sheaf $\HOM(\bE,\bF)$ is dual to $\HOM_0(\bF,\bE(\divisor'))$, while $\HOM_0(\bE,\bF\otimes\Omega_X(\divisor))$ is dual to $\HOM(\bF\otimes\Omega_X(\divisor),\bE(\divisor'))$. Applying Serre duality to the exact sequence~\eqref{eq:ExactSequence}, we get an exact sequence
  \[
    0\to\Ext^2((\bE,\nabla_\bE),(\bF,\nabla_\bF))^\vee\to
    H^0(X,\HOM(\bF,\bE(\divisor'-\divisor)))\to
    H^0(X,\HOM_0(\bF,\bE\otimes\Omega_X(\divisor'))
  \]
  and the statement follows from~\eqref{eq:ExactSequence}.
\end{proof}

\begin{lemma}\label{lm:EulerChar}
Assume that $(\bE,\nabla_\bE)\in\Conn_{\gamma_1,d_1}^{prtl}(\divisor,\divisor',\zeta)(\kk)$,
$(\bF,\nabla_\bF)\in\Conn_{\gamma_2,d_2}^{prtl}(\divisor,\divisor',\zeta)(\kk)$, and that $\gamma_1+\gamma_2\in\Gamma_{\divisor',\zeta}$. Denote the Harder--Narasimhan spectra of $\bE$ and $\bF$ by
$\tau_1>\tau_2>\ldots>\tau_m$ and $\tau'_1>\tau'_2>\ldots>\tau'_{m'}$ respectively. If $\tau_1<\tau'_{m'}$, then
\[
    \dim\Ext^1((\bE,\nabla_\bE),(\bF,\nabla_\bF))-\dim\Hom((\bE,\nabla_\bE),(\bF,\nabla_\bF))=
    \frac{\chi(\gamma_1)+\chi(\gamma_2)-\chi(\gamma_1+\gamma_2)}2,
\]
where $\chi$ is from~\eqref{eq:BetterChi3}.
\end{lemma}
\begin{proof}
  Since the Harder--Narasimhan spectrum of $\bE(\divisor'-\divisor)$ is obtained from that of $\bE$ by subtracting $\deg(\divisor-\divisor')\rk\bE$, it follows from Lemma~\ref{lm:MorphismSS} that $\Hom((\bF,\nabla_\bF),(\bE(\divisor'-\divisor),\nabla_\bE))=0$.
  Thus, by Lemma~\ref{lm:Ext}(ii), we have $\Ext^2((\bE,\nabla_\bE),(\bF,\nabla_\bF))=0$. Now we obtain from~\eqref{eq:ExactSequence}:
  \begin{multline*}
    \dim\Ext^1((\bE,\nabla_\bE),(\bF,\nabla_\bF))-\dim\Hom((\bE,\nabla_\bE),(\bF,\nabla_\bF))\\=
    \chi(\HOM_0(\bE,\bF\otimes\Omega_X(\divisor)))-\chi(\HOM(\bE,\bF))=
    -\chi(\HOM(\bF,\bE(\divisor'-\divisor))-\chi(\HOM(\bE,\bF)),
  \end{multline*}
  where $\chi$ stands for the Euler characteristic and we used Serre duality as in the proof of Lemma~\ref{lm:Ext}(ii). Now applying Corollary~\ref{cor:EulerChar0} twice, we can calculate the last expression:
  \begin{multline*}
        -\left((1-g)r^\bF r^\bE+r^\bF\deg(\bE(\divisor'-\divisor))-r^\bE d^\bF-\sum_{\substack{x\in D'\\i<j}}n'_xr^\bF_{x,i}r^\bE_{x,j}\right)\\
        -\left((1-g)r^\bE r^\bF+r^\bE d^\bF-r^\bF d^\bE-\sum_{\substack{x\in D'\\i<j}}n'_xr^\bE_{x,i}r^\bF_{x,j}\right)
        \\=(2g-2)r^\bE r^\bF-r^\bF(\deg\bE(\divisor'-\divisor)-\deg\bE)+2\sum_{\substack{x\in D'\\i<j}}n'_xr^\bE_{x,i}r^\bF_{x,j}
        \\=(2g-2+\deg(\divisor-\divisor'))r^\bE r^\bF+2\sum_{\substack{x\in D'\\i<j}}n'_xr^\bE_{x,i}r^\bF_{x,j}.
  \end{multline*}
  Lemma~\ref{lm:EulerChar} follows.
\end{proof}

\subsection{Kontsevich--Soibelman factorization formula} Recall that we have divisors $\divisor'\le\divisor$ on $X$ with supports $D'\subset D\subset X(\kk)$. Consider any series $H(w,w_{\bullet,\bullet},z)=\sum_{\gamma\in\Gamma_{D'},d\in\Z}A_{\gamma,d}w^\gamma z^d$, where $A_{\gamma,d}$ are elements of $\Mot(\kk)$ (or, more generally, of any commutative ring). Assume that $H(0,0,z)=1$. Then, using an inductive argument, we see that we can uniquely factorize
\[
H(w,w_{\bullet,\bullet},z)=\prod_{\tau\in\R}\left(1+\sum_{\substack{\gamma\in\Gamma_{D'},d\in\Z\\ \gamma\ne0,d+\gamma\star\sigma=\tau\rk\gamma}}C_{\gamma,d} w^\gamma z^d\right).
\]
We drop $\tau$ from the notation for $C_{\gamma,d}$ because $\tau$ is uniquely determined by $\gamma$ and $d$. Let, as before, $\delta:=\deg\divisor-|D'|\in\Z_{\ge0}$ and $\chi\colon\Gamma_{D'}\to\Z$ be defined by~\eqref{eq:BetterChi3}.
\begin{proposition}\label{pr:KontsevichSoibelman}
Consider the factorization
\[
    \sum_{\gamma\in\Gamma_{D'},d\le0}(-\bL^\frac12)^{-\chi(\gamma)}[\Conn_{\gamma,d}^{prtl,-}(\divisor,\divisor',\zeta)]w^\gamma z^d=
    \prod_{\tau\in\R}\left(1+\sum_{\substack{\gamma\in\Gamma_{D'},d\le0\\\gamma\ne0, d+\gamma\star\sigma=\tau\rk\gamma}}C_{\gamma,d} w^\gamma z^d\right).
\]
Then for all $\gamma'\in\Gamma_{\divisor',\zeta}$ we have $[\Conn_{\gamma',d'}^{prtl,\sigma-ss,-}(\divisor,\divisor',\zeta)]=(-\bL^\frac12)^{\chi(\gamma')}C_{\gamma',d'}$.
\end{proposition}
\begin{proof}
  This follows from the general formalism of~\cite{KontsevichSoibelman08}. In more detail, let $\gamma'\in\Gamma_{\divisor',\zeta}$. For $\gamma_1,\ldots,\gamma_m\in\Gamma_{D'}$ with $\gamma_1+\ldots+\gamma_m=\gamma'$ and $d_1,\ldots,d_m\in\Z_{\le0}$, let $\Conn_{\gamma_1,d_1\ldots,\gamma_m,d_m}$ denote the constructible subset of $\Conn^{prtl,-}(\divisor,\divisor',\zeta)$ consisting of stacks such that the Harder--Narasimhan filtration has length~$m$ and the $i$-th quotient has class $(\gamma_i,d_i)$ for $i=1,\ldots,m$. Then the motivic class of the stack $\Conn_{\gamma,d}^{prtl,-}(\divisor,\divisor',\zeta)$ is equal to the sum of motivic classes $[\Conn_{\gamma_1,d_1\ldots,\gamma_m,d_m}]$, where $m$ ranges over positive integers and $\gamma_i$ and $d_i$ range over all sequences with $\sum_i\gamma_i=\gamma'$, $\sum_id_i=d$, and such that the sequence $\tau_i:=(d_i+\gamma_i\star\sigma)/\rk\gamma_i$ is strictly decreasing. Since $\chi(\gamma_1+\ldots+\gamma_m)-\chi(\gamma_1)-\ldots-\chi(\gamma_m)$ is always even, the statement of the proposition follows from the identity
  \[
    [\Conn_{\gamma_1,d_1\ldots,\gamma_m,d_m}]=\bL^{\frac12(\chi(\gamma_1+\ldots+\gamma_m)-\chi(\gamma_1)-\ldots-\chi(\gamma_m))}
    \prod_{i=1}^m[\Conn_{\gamma_i,d_i}^{prtl,\sigma-ss,-}(\divisor,\divisor',\zeta)].
  \]
  This statement is proved as in~\cite[Prop.~3.6.1]{FedorovSoibelmans} with the help of Lemma~\ref{lm:Ext}(i),~\ref{lm:EulerChar}, and~\ref{lm:Summand}.
\end{proof}

\begin{corollary}\label{cor:ExplNonnegativeSemistable}
(i) Assume that $\divisor'<\divisor$, $\zeta\in\FNF(\divisor,\divisor')$, and $\gamma'\in\Gamma_{\divisor',\zeta}$. Set $\delta:=\deg\divisor-|D'|$. Then the motivic class of $\Conn^{prtl,\sigma-ss,-}_{\gamma',d'}(\divisor,\divisor',\zeta)$ in $\cMot(\kk)$ is equal to the coefficient at $w^{\gamma'}z^{-d'}$ of
\[
    (-\bL^\frac12)^{\chi(\gamma')} \Exp\Bigl(\left(\Log\Omega^{Sch,mot}_{X,D',\delta}\right)_{-d+\gamma\star\sigma=\tau\rk\gamma}\Bigr),
\]
where $\tau=\frac{d'+\gamma'\star\sigma}{\rk\gamma'}$.

(ii) Assume that $\zeta\in\FNF(\divisor)$ and $\gamma'\in\Gamma_{\divisor,\zeta}$. Let $\delta:=\deg\divisor-|D|$ be the irregularity of $\divisor$. Then the class of $\Conn^{\sigma-ss,-}_{\gamma',d'}(\divisor,\zeta)$ in $\cMot(\kk)$ is equal to the coefficient at $w^{\gamma'}z^{-d'}$ of
\[
    (-\bL^\frac12)^{\chi(\gamma')} \Exp\Biggl(\left(\bL\Log\Omega^{Sch,mot}_{X,D,\delta}\right)_{\substack{-d+\gamma\star\sigma=\tau\rk\gamma\\-\epsilon d+\gamma\star\zeta=0}}\Biggr),
\]
where $\tau=\frac{d'+\gamma'\star\sigma}{\rk\gamma'}$.
\end{corollary}
\begin{proof}
We apply the homomorphism $\Mot(\kk)\to\cMot(\kk)$ to the identity of Proposition~\ref{pr:KontsevichSoibelman}. Then we apply the plethystic logarithm. It remains to use Proposition~\ref{pr:ExplicitPartial-} in part~(i) and Proposition~\ref{pr:ExplicitFull-} in part~(ii). The requirement that $\gamma'\in\Gamma_{\divisor',\zeta}$ carries through all these steps thanks to Lemma~\ref{lm:GoodExp}.
\end{proof}

\section{Negative degree limit of parabolic \texorpdfstring{$\epsilon$}{epsilon}-connections}\label{sect:Stabilization} Recall that we have divisors $\divisor'=\sum_{x\in D'}n'_xx\le\divisor=\sum_{x\in D}{n_xx}$ on $X$, and we are studying the moduli stacks $\Conn_{\gamma,d}^{prtl}(\epsilon,\divisor,\divisor',\zeta)$ of parabolic $\epsilon$-connections with partially fixed normal forms, where $(\gamma,d)\in\Gamma_{D'}\times\Z$, $\epsilon\in\kk$, and $\zeta\in\FNF(\divisor,\divisor')$, see Section~\ref{sect:IrregFNF}.

So far, we have calculated motivic classes of various moduli stacks of parabolic $\epsilon$-connections with nonpositive underlying vector bundles. In this section, we get rid of this technical assumption, calculating motivic classes of the stacks $\Conn_{\gamma,d}(\epsilon,\divisor,\zeta)=\Conn_{\gamma,d}^{prtl}(\epsilon,\divisor,\divisor,\zeta)$, whenever $\epsilon\ne0$, as well as motivic classes of semistable loci of the stacks $\Conn_{\gamma,d}^{prtl}(\epsilon,\divisor,\divisor',\zeta)$ without any restriction on $\epsilon$; see Theorems~\ref{th:AnswerConnInMot} and~\ref{th:AnswerInCmot}. Our argument is somewhat different for Higgs bundles ($\epsilon=0$) and connections ($\epsilon\ne0$), so we will write $\epsilon$ explicitly in the notation for the stacks of $\epsilon$-connections.

\subsection{Twisting parabolic $\epsilon$-connections}\label{sect:TwistConn} Fix a point $x\in D$. Let $\zeta_x^\epsilon$ be the canonical section of $\Omega_X(x)/\Omega_X$ with $\res\zeta_x^\epsilon=\epsilon$. For $N\in\Z$ and $\zeta=(\zeta_{x,j})\in\FNF(\divisor)$, we define $\zeta(N)\in\FNF(\divisor)$ by the formulas $(\zeta(N))_{y,j}=\zeta_{y,j}$ whenever $y\ne x$ and $(\zeta(N))_{x,j}=\zeta_{x,j}+N\zeta_x^\epsilon$.

\begin{lemma}\label{lm:tensorisation}
  Let $x\in D$ be a point.

  (i) Assume that $n'_x<n_x$. Then the tensorisation with $\cL:=\cO_X(-x)$ gives an isomorphism
  \[
    \Conn^{prtl}_{\gamma,d}(\epsilon,\divisor,\divisor',\zeta)\to
    \Conn^{prtl}_{\gamma,d-\rk\gamma}(\epsilon,\divisor,\divisor',\zeta).
  \]
  If $\gamma\in\Gamma_{\divisor',\zeta}$ and $\sigma$ is a sequence of parabolic weights of type $(X,\divisor')$, then this isomorphism takes the substack of $\sigma$-semistable parabolic $\epsilon$-connections onto the substack of $\sigma$-semistable parabolic $\epsilon$-connections.

  (ii) The tensorisation with $\cL:=\cO_X(-x)$ gives an isomorphism
  \[
    \Conn_{\gamma,d}(\epsilon,\divisor,\zeta)\to
    \Conn_{\gamma,d-\rk\gamma}(\epsilon,\divisor,\zeta(1)).
  \]
  If $\gamma\in\Gamma_{\divisor,\zeta}$ and $\sigma$ is a sequence of parabolic weights of type $(X,\divisor)$, then this isomorphism takes the substack of $\sigma$-semistable parabolic $\epsilon$-connections onto the substack of $\sigma$-semistable parabolic $\epsilon$-connections.
\end{lemma}
\begin{proof}
  (i) Let $(E,E_{\bullet,\bullet},\nabla)\in\Conn^{prtl}(\epsilon,\divisor,\divisor',\zeta)(\kk)$. Since $E|_{X-x}=E(-x)|_{X-x}$, the $\epsilon$-connection~$\nabla$ gives rise to a meromorphic $\epsilon$-connection~$\nabla'$ on $E(-x)$. Since $\nabla$ has a pole of order at most $n_x$ at~$x$ and $n_x>0$, it is easy to see that $\nabla'$ has a pole of order at most $n_x$ at $x$. Thus $\nabla'$ is an irregular $\epsilon$-connection on $E(-x)$ with poles bounded by $\divisor$. We need to show that $(E(-x),E_{\bullet,\bullet}(-x),\nabla')\in\Conn^{prtl}(\epsilon,\divisor,\divisor',\zeta)(\kk)$. By Definition~\ref{def:ModSpacePartial}, trivializing $E$ near $x$ compatibly with $E_{x,\bullet}$, we can write:
  \[
        \nabla=\epsilon d+A+B,
  \]
  where $A$ is a block upper triangular matrix with scalar diagonal blocks with eigenvalues $\zeta_{x,j}$, and $B$ has a pole of order at most $n_x-n'_x$ at $x$. Choose a coordinate $z$ near $x$, in the induced local trivialization of $E(-x)$, we get
  \begin{equation}\label{eq:TwistedFNF}
    \nabla'=\epsilon d+A+B+\epsilon\frac{dz}z\Id.
  \end{equation}
  Since $n'_x<n_x$, the matrix-valued 1-form $B+\epsilon\frac{dz}z\Id$ has still a pole of order at most $n_x-n'_x$, and the claim follows. It is clear that the construction works in families, so we get a morphism of stacks. We construct the inverse morphism similarly. The statement about semistable loci follows from the definitions.

  (ii) Similar to part~(i). In more detail, consider~\eqref{eq:TwistedFNF} and note that $A+\epsilon\frac{dz}z\Id$ is also block upper triangular with diagonal blocks being scalar matrices with eigenvalues $\zeta(1)_{x,j}$.
\end{proof}

\subsection{Negative degree limit of connections with fully fixed formal normal forms}\label{sect:StabConn} We first consider the case $\epsilon\ne0$ and $\divisor'=\divisor$, that is, we work with connections (as opposed to Higgs bundles) and we fix the full formal normal form at each singular point. In this case, the stack $\Conn_{\gamma,d}(\epsilon,\divisor,\zeta)$ is already of finite type, at least, when $\gamma\in\Gamma_{\divisor,\zeta}$; see Corollary~\ref{cor:StabConn}. Thus, we do not need to add stability conditions (though we may, see Theorem~\ref{th:AnswerInCmot}(ii)). The remaining cases of Higgs bundles or connections with partially fixed normal forms will be considered in Section~\ref{sect:RemainingCase}.

We first estimate the gaps in the Harder--Narasimhan spectra of the underlying vector bundles of indecomposable parabolic bundles (cf.~\cite[Lemma~8.4]{FedorovSoibelmansParabolic}) underlying parabolic $\epsilon$-connections.
\begin{lemma}\label{lm:l}
  Assume that $(\bE,\nabla)\in\Conn^{prtl}_{\gamma,d}(\epsilon,\divisor,\divisor',\zeta)(\kk)$, where $\gamma\in\Gamma_{\divisor',\zeta}$. Assume also that~$\bE$ is indecomposable. Let $\tau_1>\tau_2>\ldots>\tau_m$ be the Harder--Narasimhan spectrum of the underlying vector bundle of $\bE$. Then for all $i$ we have
  \[
    \tau_i-\tau_{i+1}\le2g-2+\deg\divisor,
  \]
  where $g$ is the genus of $X$.
\end{lemma}
\begin{proof}
 Write $\bE=(E,E_{\bullet,\bullet})$. Let $0=E_0\subset E_1\subset\ldots\subset E_m=E$ be the Harder--Narasimhan filtration of they underlying vector bundle $E$ of $\bE$. Assume for a contradiction that $\tau_i>\tau_{i+1}+2g-2+\deg\divisor$ for some~$i$. We claim that $E_i$ is preserved by $\nabla$. Indeed, the composition
 \[
    \Phi\colon E_i\hookrightarrow E\xrightarrow{\nabla}E\otimes\Omega_X(\divisor)\twoheadrightarrow(E/E_i)\otimes\Omega_X(\divisor)
\]
is $\cO_X$-linear. On the other hand, the Harder--Narasimhan spectrum of $E_i$ is contained in $[\tau_i,\infty)$, while that of $(E/E_i)\otimes\Omega_X(\divisor)$ is contained in $(-\infty,\tau_{i+1}+2g-2+\deg\divisor]\subset(-\infty,\tau_i)$.
Therefore $\Hom(E_i,(E/E_i)\otimes\Omega_X(\divisor))=0$, so $\Phi=0$. Thus $\nabla$ preserves $E_i$. By Corollary~\ref{cor:Subbundle} we see that $\bE_i:=(E_i,E_i\cap E_{\bullet,\bullet})$ is a parabolic subbundle of $\bE$, so we have an exact sequence $0\to\bE_i\to\bE\to\bE/\bE_i\to0$. According to Lemma~\ref{lm:SerreDuality} the vector space $\Ext^1(\bE/\bE_i,\bE_i)$ is dual to a vector subspace of $\Hom(E_i,(E/E_i)\otimes\Omega_X(\divisor))$. However, as we have already explained, this space is zero, so $\Ext^1(\bE/\bE_i,\bE_i)=0$ as well. Hence, the exact sequence is split and $\bE$ is decomposable.
\end{proof}

The above lemma has the following corollary (cf.~\cite[Lemma~2.8]{MozgovoySchiffmann2020} and~\cite[Lemma~7.2]{FedorovSoibelmansParabolic}):
\begin{lemma}\label{lm:EstimateTau}
  Assume that $(\bE,\nabla)\in\Conn^{prtl}_{\gamma,d}(\epsilon,\divisor,\divisor',\zeta)(\kk)$, where $\gamma\in\Gamma_{\divisor',\zeta}$. Assume also that~$\bE$ is indecomposable. Let $\tau_1>\tau_2>\ldots>\tau_m$ be the Harder--Narasimhan spectrum of the underlying vector bundle of $\bE$. Then
  \[
    \tau_1\le\frac dr+\frac{(r-1)l}2 ,
  \]
  where $l=\max(2g-2+\deg\divisor,0)$, $r=\rk\bE$.
\end{lemma}
\begin{proof}
  By the previous lemma we have $\tau_i\ge\tau_1-(i-1)l$. Let $r_i$ be the rank of the $i$-th subquotient in the Harder--Narasimhan filtration of the underlying vector bundle of $\bE$. Therefore
  \[
    \frac dr=\frac{\sum_i\tau_ir_i}r\ge\frac{\sum_i\tau_1r_i}r-\frac{\sum_i(i-1)lr_i}r
    \ge\tau_1-\frac{(r-1)l}2,
  \]
  where we used that $\sum_ir_i=r$ and $\sum_i(i-1)r_i\le r(r-1)/2$. The statement follows.
\end{proof}

We are now assuming that $\epsilon\ne0$ and $\divisor'=\divisor$. Recall that $\zeta=(\zeta_{x,j})\in\FNF(\divisor)$. Let $|\bullet|$ be any norm on the $\Q$-vector subspace of $\kk$ generated by $\epsilon$ and~$\res_x\zeta_{x,j}$, where $x$ ranges over $D$ and $j$ ranges over positive integers. If~$\kk$ is embedded into~$\C$, we can take the usual absolute value for $|\bullet|$. For $\gamma=(r,r_{\bullet,\bullet})\in\Gamma_D$, we set $|\zeta|_\gamma:=\sum_{x\in D}(\max_{i\colon r_{x,i}\ne0}|\res_x\zeta_{x,i}|)$. We note that if $(E,E_{\bullet,\bullet})\in\Bun_{\gamma,d}^{par}(X,\divisor)(\kk)$ has $(\epsilon,\zeta)$-degree zero, that is, $\epsilon d+\gamma\star\zeta=0$, then
\begin{equation}\label{eq:norm}
  \frac dr\le\frac{|\zeta|_\gamma}{|\epsilon|}.
\end{equation}

The following lemma shows that the Harder--Narasimhan spectra of vector bundles underlying parabolic $\epsilon$-connections in $\Conn_{\gamma,d}(\epsilon,\divisor,\zeta)$ are uniformly bounded above (independently of $d$), provided that $\gamma\in\Gamma_{\divisor,\zeta}$.
\begin{lemma}\label{lm:ConnBounded}
  Let $(E,E_{\bullet,\bullet},\nabla)\in\Conn_{\gamma,d}(\epsilon,\divisor,\zeta)(\kk)$, where $\gamma\in\Gamma_{\divisor,\zeta}$. Let $\tau_1>\tau_2>\ldots>\tau_m$ be the Harder--Narasimhan spectrum of $E$. Then
  \[
    \tau_1\le\frac{|\zeta|_\gamma}{|\epsilon|}+\frac{(r-1)l}2,
  \]
  where $r=\rk E$ and $l$ as in Lemma~\ref{lm:EstimateTau}.
\end{lemma}
\begin{proof}
  We write $(E,E_{\bullet,\bullet})=\bigoplus_i(E^i,E^i_{\bullet,\bullet})$, where $(E^i,E^i_{\bullet,\bullet})$ are indecomposable. Let $\tau_1^i>\ldots>\tau_{m_i}^{i}$ be the HN-spectrum of $E^i$. Then for some $i$ we have $\tau_1^i\ge\tau_1$. Next, $\nabla$ induces an $\epsilon$-connection $\nabla^i$ on $E^i$ and $(E^i,E^i_{\bullet,\bullet},\nabla^i)\in\Conn_{\gamma_i,d_i}(\epsilon,\divisor,\zeta)(\kk)$ for some $\gamma_i\in\Gamma_{\divisor,\zeta}$ and $d_i\in\Z$. Hence by the previous lemma,
  \[
    \tau_1\le\tau_1^i\le\frac{d_i}{r_i}++\frac{(r_i-1)l}2,
  \]
  where $r_i=\rk\gamma_i$. By Corollary~\ref{cor:indecomposable summand}, the $(\epsilon,\zeta)$-degree of $(E^i,E^i_{\bullet,\bullet})$ is zero. Thus, by~\eqref{eq:norm} we have $\frac{d_i}{r_i}\le\frac{|\zeta|_{\gamma_i}}{|\epsilon|}\le\frac{|\zeta|_\gamma}{|\epsilon|}$. It remains to notice that $r_i\le r$.
\end{proof}

\begin{corollary}\label{cor:StabConn}
Assume that $\epsilon\ne0$ and $\gamma\in\Gamma_{\divisor,\zeta}$. Then for $N\gg0$ we have
\[
      \Conn_{\gamma,d}(\epsilon,\divisor,\zeta)\simeq
      \Conn_{\gamma,d-N\rk\gamma}^-(\epsilon,\divisor,\zeta(N)).
\]
In particular, the stack $\Conn_{\gamma,d}(\epsilon,\divisor,\zeta)$ is of finite type.
\end{corollary}
\begin{proof}
  Let $(E,E_{\bullet,\bullet},\nabla)\in\Conn_{\gamma,d}(\epsilon,\divisor,\zeta)(K)$, where $K\supset\kk$ is a field extension. Since twisting by $\cO_X(-x)$ shifts the HN-spectrum of $E$ by $-1$, we see from Lemma~\ref{lm:ConnBounded} (upon base changing to $K$) that $E(-Nx)$ is HN-nonpositive, once $N\ge\frac{|\zeta|_\gamma}{|\epsilon|}+\frac{(r-1)l}2$. It remains to iterate $N$ times the isomorphism of Lemma~\ref{lm:tensorisation}(ii).
\end{proof}

Now we are ready for our first main result. We refer the reader to Section~\ref{sect:ConnFnf} for the notations.
\begin{theorem}\label{th:AnswerConnInMot}
Let, as before, $\divisor=\sum_{x\in D}{n_xx}$, where for all $x$ we have $n_x>0$, be a divisor on $X$, $\zeta\in\FNF(\divisor)$. Assume that $\epsilon\ne0$ and $\gamma'\in\Gamma_{\divisor,\zeta}$. Set $\delta:=\deg\divisor-|D|$ and let as in~\eqref{eq:BetterChi2} for $\gamma\in\Gamma_D$
\[
    \chi(\gamma)=\chi(r,r_{\bullet,\bullet}):=(2g-2)r^2-\delta r+2\sum_{x\in D,i<j}n_xr_{x,i}r_{x,j}.
\]
Then for all $N\gg0$ the motivic class of $\Conn_{\gamma',d'}(\epsilon,\divisor,\zeta)$ in $\cMot(\kk)$ is equal to the coefficient of
\[
    (-\bL^\frac12)^{\chi(\gamma')}\Exp\left(\Bigl(\bL\Log\Omega^{Sch,mot}_{X,D,\delta}\Bigr)_{-\epsilon d+\gamma\star\zeta=-\epsilon N\rk\gamma}\right)
\]
at $w^{\gamma'}z^{-d'+N\rk\gamma'}$.
\end{theorem}
We remark that $[\Conn_{\gamma',d'}(\epsilon,\divisor,\zeta)]=0$ unless $\epsilon d'+\gamma'\star\zeta=0$.
\begin{proof}
According to Corollary~\ref{cor:StabConn}, for $N\gg0$, we have
\[
    [\Conn_{\gamma',d'}(\epsilon,\divisor,\zeta)]=[\Conn_{\gamma',d'-N\rk\gamma'}^-(\epsilon,\divisor,\zeta(N))].
\]
By Proposition~\ref{pr:ExplicitFull-} this motivic class is equal to the coefficient at $w^{\gamma'}z^{-d'+N\rk\gamma'}$ in
\[
    (-\bL^\frac12)^{\chi(\gamma')}\Exp\left(\Bigl(\bL\Log\Omega^{Sch,mot}_{X,D,\delta}\Bigr)_{-\epsilon d+\gamma\star\zeta(N)=0}\right).
\]
We have
\begin{equation}\label{eq:d+zeta(N)}
    -\epsilon d+\gamma\star\zeta(N)=-\epsilon d+\gamma\star\zeta+\sum_{j>0}r_{x,j}(N\epsilon)=
    -\epsilon d+\gamma\star\zeta+\epsilon N\rk\gamma.
\end{equation}
Now the claim follows.
\end{proof}

\subsection{Negative degree limit of semistable parabolic $\epsilon$-connections in the remaining cases}\label{sect:RemainingCase} We now add stability conditions, so we can consider the cases when $\epsilon=0$ (that is, the case of Higgs bundles) or $\divisor'<\divisor$. We study the stacks $\Conn^{prtl,\sigma-ss}_{\gamma,d}(\epsilon,\divisor,\divisor',\zeta)$ defined in Section~\ref{sect:DefStab}. We recall that these stacks are only defined when $\gamma\in\Gamma_{\divisor',\zeta}\subset\Gamma_{D'}$, see Definition~\ref{def:Semistable}. The following is an irregular version of~\cite[Lemma~7.3]{FedorovSoibelmansParabolic}.

\begin{proposition}\label{pr:Stabilization}
  Let $\sigma$ be a sequence of parabolic weights of type $(X,\divisor')$ as in Section~\ref{sect:DefStab}. Assume that $\gamma\in\Gamma_{\divisor',\zeta}$. Then there is $d_0=d_0(g,\divisor,\rk\gamma,\sigma)$ such that for $d\le d_0$ we have
  \[
    \Conn^{prtl,\sigma-ss}_{\gamma,d}(\epsilon,\divisor,\divisor',\zeta)=
    \Conn^{prtl,\sigma-ss,-}_{\gamma,d}(\epsilon,\divisor,\divisor',\zeta).
  \]
\end{proposition}
\begin{proof}
    Take $(\bE,\nabla)\in\Conn^{prtl}_{\gamma,d}(\epsilon,\divisor,\divisor',\zeta)(K)$, where $K$ is an extension of $\kk$. Let $\tau_1>\tau_2>\ldots>\tau_m$ be the Harder--Narasimhan spectrum of $E$, where $E$ is the underlying vector bundle of $\bE$. Write $\gamma=(r,r_{\bullet,\bullet})$. Replacing $\sigma_{x,j}$ with 0 for all $x$ and $j$ with $r_{x,j}=0$ does not change our stacks, so we may assume that only finitely many components of $\sigma$ are non-zero. Let $|\sigma|:=\sum_{x\in D'}(\max_i|\sigma_{x,i}|)$. We show that if $(\bE,\nabla)$ is $\sigma$-semistable, then for all $i$ we have $\tau_i-\tau_{i+1}\le l:=\max(2g-2+\deg\divisor,2|\sigma|)$. Assume the contrary. By Lemma~\ref{lm:l} and its proof, we can decompose $\bE=\bE_i\oplus\bF$, where the underlying vector bundle $E_i$ of $\bE_i$ is the $i$-th term of the Harder--Narasimhan filtration of $E$ and, moreover, $E_i$ is preserved by $\nabla$. The slope of $E_i$ is at least $\tau_i$, so the $\sigma$-slope of $\bE_i$ is at least $\tau_i-|\sigma|$. Similarly, the $\sigma$-slope of $\bF$ is at most $\tau_{i-1}+|\sigma|$ and it follows from our assumption on $\tau_i-\tau_{i+1}$ that $\bE_i\subset\bE$ is a~destabilizing subbundle. This contradiction shows that $\tau_i-\tau_{i+1}\le l$ for all $i$.

    Next, take $d_0$ such that $\frac{d_0}r+\frac{(r-1)l}2<-2|\sigma|$ and let $d\le d_0$. Exactly as in the proof of Lemma~\ref{lm:EstimateTau} we obtain
    \[
        \tau_1\le\frac dr+\frac{(r-1)l}2<0
    \]
    so that $(\bE,\nabla)\in\Conn^{prtl,\sigma-ss,-}_{\gamma,d}(\epsilon,\divisor,\divisor',\zeta)(K)$.

    Conversely, assume that $(\bE,\nabla)\in\Conn^{prtl,\sigma-ss,-}_{\gamma,d}(\epsilon,\divisor,\divisor',\zeta)(K)$. We need to show that it is $\sigma$-semistable. Assume the contrary, and let $\bE'$ be the minimal destabilizing quotient bundle. Let $E'$ be the underlying vector bundle of $\bE'$. Denote by $d'$ and $r'$ the degree and the rank of $E'$ and by $\tau'_1>\tau'_2>\ldots$ the Harder--Narasimhan spectrum of the $E'$. Note that the $\sigma$-slope of $\bE'$ is smaller than that of $\bE$, so $d'/r'<d/r+2|\sigma|$. Then, since $(\bE',\nabla|_{\bE'})$ is $\sigma$-semistable, we obtain similarly to the above
    \[
        \tau'_1\le\frac{d'}{r'}+\frac{(r'-1)l}2\le\frac dr+2|\sigma|+\frac{(r-1)l}2<0
    \]
    so that $\bE'$ is nonpositive. This, however, contradicts the assumption that $\bE$ is nonpositive-semistable.
\end{proof}

Combining Proposition~\ref{pr:Stabilization} with Lemma~\ref{lm:tensorisation}, and noting that $\zeta(1)=\zeta$ when $\epsilon=0$, we get the following corollary.
\begin{corollary}\label{cor:StabHiggs}
Assume that $\divisor'<\divisor$ or that $\epsilon=0$. Let $\gamma\in\Gamma_{\divisor',\zeta}$. Then there is $N_0=N_0(g,\divisor,\rk\gamma,\sigma)$ such that for $N\ge N_0$ we have
\[
      \Conn^{prtl,\sigma-ss}_{\gamma,d}(\epsilon,\divisor,\divisor',\zeta)\simeq
      \Conn^{prtl,\sigma-ss,-}_{\gamma,d-N\rk\gamma}(\epsilon,\divisor,\divisor',\zeta).
\]
In particular, the stack $\Conn^{prtl,\sigma-ss}_{\gamma,d}(\epsilon,\divisor,\divisor',\zeta)$ is of finite type.
\end{corollary}
Now we are ready for our second main result, where the notation is from Definition~\ref{def:ModSpacePartial}.
\begin{theorem}\label{th:AnswerInCmot}
Let, as before, $\divisor'=\sum_{x\in D'}n'_xx\le\divisor=\sum_{x\in D}{n_xx}$ be effective divisors on $X$ with supports $D'$ and $D$ respectively. Let $\zeta\in\FNF(\divisor,\divisor')$ and let $\sigma$ be a sequence of parabolic weights of type $(X,\divisor')$. Assume that $\gamma'\in\Gamma_{\divisor',\zeta}$ and $d'\in\Z$. Set $\delta:=\deg\divisor-|D'|$ and as in~\eqref{eq:BetterChi3} set
\[
    \chi(\gamma)=\chi(r,r_{\bullet,\bullet}):=
(2g-2+\deg(\divisor-\divisor'))r^2+r\sum_{x\in D'}(1-n'_x)+2\sum_{x\in D',i<j}n'_xr_{x,i}r_{x,j}.
\]
Let $\tau:=\frac{d'+\gamma'\star\sigma}{\rk\gamma'}$.

(i) Assume that $\divisor'<\divisor$. Then for all $N\gg0$ the motivic class of $\Conn^{prtl,\sigma-ss}_{\gamma',d'}(\epsilon,\divisor,\divisor',\zeta)$ in $\cMot(\kk)$ is equal to the coefficient at $w^{\gamma'}z^{-d'+N\rk\gamma'}$ in
\[
    (-\bL^\frac12)^{\chi(\gamma')} \Exp\left(\Bigl(\Log\Omega^{Sch,mot}_{X,D',\delta}\Bigr)_{-d+\gamma\star\sigma=(\tau+N)\rk\gamma}\right).
\]

(ii) Assume that $\divisor'=\divisor$. Then for all $N\gg0$ the motivic class of $\Conn^{\sigma-ss}_{\gamma',d'}(\epsilon,\divisor,\zeta)$ in $\cMot(\kk)$ is equal to the coefficient at $w^{\gamma'}z^{-d'+N\rk\gamma'}$ in
\[
    (-\bL^\frac12)^{\chi(\gamma')}\Exp\left(\Bigl(\bL\Log\Omega^{Sch,mot}_{X,D,\delta}\Bigr)_%
    {\substack{-\epsilon d+\gamma\star\zeta=-\epsilon N\rk\gamma\\-d+\gamma\star\sigma=(\tau+N)\rk\gamma}}\right).
\]
\end{theorem}
\begin{proof}
  (i) Combine the previous corollary with Corollary~\ref{cor:ExplNonnegativeSemistable}(i) and~\eqref{eq:d+zeta(N)}.

  (ii) If $\epsilon=0$, we combine the previous corollary with Corollary~\ref{cor:ExplNonnegativeSemistable}(ii) applied with $\epsilon=0$ and use~\eqref{eq:d+zeta(N)}. If $\epsilon\ne0$, we first write the following version of Kontsevich--Soibelman factorization formula (whose proof is completely similar to~Proposition~\ref{pr:KontsevichSoibelman}). Consider the factorization
\begin{equation}\label{eq:factor1}
    \sum_{\gamma\in\Gamma_D,d\in\Z}(-\bL^\frac12)^{-\chi(\gamma)}[\Conn_{\gamma,d}(\epsilon,\divisor,\zeta)]w^\gamma z^d=
    \prod_{\tau\in\R}\left(1+\sum_{\substack{\gamma\in\Gamma_D,d\in\Z\\ \gamma\ne0,-d+\gamma\star\sigma=\tau\rk\gamma}}C_{\gamma,d} w^\gamma z^d\right).
\end{equation}
Then for all $\gamma\in\Gamma_{\divisor,\zeta}$ and $d\in\Z$, we have $[\Conn_{\gamma,d}^{\sigma-ss}(\epsilon,\divisor,\zeta)]=(-\bL^\frac12)^{\chi(\gamma)}C_{\gamma,-d}$.

Next, let us make a substitution $\Omega^{Sch,mot}_{X,D,\delta}=\Omega^{Sch,mot}_{X,D,\delta}(z,w,w_{\bullet,\bullet})\mapsto \Omega^{Sch,mot}_{X,D,\delta}(z,wz^N,w_{\bullet,\bullet})$ and note that the coefficient of the new series at $w^\gamma z^{-d}$ is equal to the coefficient of $\Omega^{Sch,mot}_{X,D,\delta}$ at $w^\gamma z^{-d+N\rk\gamma}$. Since the substitution $w\mapsto wz^N$ commutes with plethystic $\Exp$ and $\Log$, we can re-interpret Theorem~\ref{th:AnswerConnInMot} in the following way: for all $\gamma''\in\Gamma_{\divisor,\zeta}$, there is $N_0(\gamma'')$ such that for all $d''\in\Z$ and $N\ge N_0(\gamma'')$ the class $(-\bL^\frac12)^{-\chi(\gamma'')}[\Conn_{\gamma'',d''}(\epsilon,\divisor,\zeta)]$ in $\cMot(\kk)$ is equal to the coefficient of
\[
    \Exp\left(\Bigl(\bL\Log\Omega^{Sch,mot}_{X,D,\delta}(z,wz^N,w_{\bullet,\bullet})\Bigr)_{-\epsilon d+\gamma\star\zeta=0}\right)
\]
at $w^{\gamma''}z^{-d''}$. Take $N$ such that for all $\gamma''\le\gamma'$ we have $N\ge N(\gamma'')$. We can factorize
\begin{multline}\label{eq:factor2}
    \Exp\left(\Bigl(\bL\Log\Omega^{Sch,mot}_{X,D,\delta}(z,wz^N,w_{\bullet,\bullet})\Bigr)_{-\epsilon d+\gamma\star\zeta=0}\right)\\=
    \prod_{\tau\in\R}\Exp\left(\Bigl(\bL\Log\Omega_{X,D,\delta}^{Sch,mot}(z,wz^N,w_{\bullet,\bullet})\Bigr)_%
    {\substack{-\epsilon d+\gamma\star\zeta=-\epsilon N\rk\gamma\\-d+\gamma\star\sigma=(\tau+N)\rk\gamma}}\right).
\end{multline}
Let us compare the left hand sides of~\eqref{eq:factor1} and~\eqref{eq:factor2}. By our choice of $N$, the coefficients at $w^{\gamma''}z^{-d}$ coincide whenever $\gamma''\le\gamma'$. Thus, the coefficients in the right hand side at $w^{\gamma'} z^{-d'}$ coincide as well. The required identity follows.
\end{proof}

\begin{remark}
When $\divisor$ is a reduced divisor, that is, when $\delta=0$, Theorems~\ref{th:AnswerConnInMot} and~\ref{th:AnswerInCmot}(ii) above become~\cite[Thm.~8.8,~7.5, and~8.16]{FedorovSoibelmansParabolic}. The twist by $\bL^{(g-1)\langle\lambda,\lambda\rangle}$ was omitted in the definition of Donaldson--Thomas invariants in~\cite[Sect.~1.3]{FedorovSoibelmansParabolic}. We use this opportunity to make a correct definition. When $\divisor=0$, we recover~\cite[Thm.~1.3.3, Cor.~1.3.4]{FedorovSoibelmans}.
\end{remark}

\section{Simplifying Schiffmann's generating function}\label{sect:MellitSimplification}
Recall from Section~\ref{sect:MotClasses} that $\sMot(\kk)$ is the universal $\lambda$-ring quotient of the pre-$\lambda$-ring $\cMot(\kk)$. The formulas given in Theorems~\ref{th:AnswerConnInMot} and~\ref{th:AnswerInCmot} can be simplified if we work in $\sMot(\kk)$: we can replace the complicated generating function $\Omega^{Sch,mot}_{g,D,\delta}$ with a much simpler generating function $\Omega^{mot}_{g,D,\delta}$. In this section we compare the two generating functions (see Theorem~\ref{th:MellitIdentity}). We will work in ``a~universal'' $\lambda$-ring $V_g$ (to be defined below in Corollary~\ref{cor:Kg}). Applying a homomorphism from~$V_g$ to $\sMot(\kk)[[z]]$ (see Proposition~\ref{pr:HomFromUniv}(iii)), we obtain the needed relation (see Corollary~\ref{cor:MellitSimpl}). In the case of Higgs bundles with regular singularities, the simplification goes back to Mellit, see~\cite{MellitIntegrality,MellitNoPunctures,MellitPunctures}. Our argument is based on Mellit's results. We remark that the notion of admissibility, which is central for this simplification, goes back to~\cite[Sect.~6, Thm.~9]{KontsevichSoibelman10}.

\subsection{$\lambda$-rings and plethystic operations}\label{Sect:SpecialRings}
Concerning $\lambda$-rings, we will use the terminology of~\cite[Exp.~V]{SGA6} and of the book~\cite{knutson1973rings}. In particular, we will speak about pre-$\lambda$-rings and $\lambda$-rings rather than about $\lambda$-rings and special $\lambda$-rings. The latter terminology is used, e.g., in the introductory exposition of Grothendieck~\cite[Exp.~0]{SGA6} or in a more recent paper~\cite{LarsenLuntsRational}.

In more detail, as we already discussed in Section~\ref{sect:MotClasses}, a \emph{pre-$\lambda$-ring\/} is a ring $R$ endowed with a group homomorphism $\lambda\colon R\to(1+zR[[z]])^\times$ such that $\lambda\equiv1+\Id_R\pmod{z^2R[[z]]}$. We write $\lambda=\sum_{n=0}^\infty\lambda_nz^n$, so that the previous conditions becomes $\lambda_0=1$, $\lambda_1=\Id_R$. A \emph{$\lambda$-ring\/} is a pre-$\lambda$-ring~$R$ such that $\lambda(1)=1+z$ and for all $x,y\in R$ and $m,n\in\Z_{\ge0}$, $\lambda_n(xy)$ and $\lambda_n(\lambda_m(x))$ can be expressed in terms of $\lambda_i(x)$ and $\lambda_j(y)$ using the ring operations in a standard way (see~\cite[Exp.~5, Def.~2.4, (2.4.1)]{SGA6} and~\cite[Ch.~I, Sect.~1]{knutson1973rings}).

Let $F\in\Sym_\Z[x_\bullet]$ be any symmetric function. Then we can uniquely write $F=P(e_1,e_2,\ldots)$, where~$P$ is a polynomial and $e_n$ are the elementary symmetric functions. If $R$ is a pre-$\lambda$-ring, we define the \emph{plethystic action of $F$ on $r\in R$\/} by $F[r]:=P(\lambda_1(r),\lambda_2(r),\ldots)$. In other words, for a fixed~$r$, the assignment $F\mapsto F[r]$ is the unique ring homomorphism $\Sym_\Z[x_\bullet]\to R$ sending $e_n$ to $\lambda_n(r)$. In particular, we have the plethystic Adams operations $\psi_n\colon r\mapsto p_n[r]$, where $p_n=\sum_ix_i^n$ is the $n$-th power sum symmetric function. If $R$ is a $\lambda$-ring, then $\psi_n\colon R\to R$ is a homomorphism. Moreover,
\begin{equation}\label{eq:LambdaPsi}
    \psi_0=1,\quad\psi_1=\Id_R,\quad\text{and for all $n,m\in\Z_{\ge0}$ we have }\psi_{nm}=\psi_n\circ\psi_m,
\end{equation}
see~\cite[Ch.~I, Sect.~4]{knutson1973rings}. Conversely, if $\Q\subset R$, then, since $p_n$ freely generate $\Sym_\Q[x_\bullet]$, given operations $\psi_n$ as above, we can define the plethystic action of $\Sym_\Q[x_\bullet]$ and, in particular, the operations~$\lambda_n$. If $\psi_n$ satisfy~\eqref{eq:LambdaPsi}, we get a $\lambda$-ring structure on $R$. Thus, if $\Q\subset R$, then the $\lambda$-ring structure on $R$ is the same as a sequence of homomorphisms $\psi_n$ satisfying~\eqref{eq:LambdaPsi}, see loc.~cit.~or~\cite[Sect.~2.1]{MellitPunctures}. It follows from $\lambda(1)=1+z$ that every $\lambda$-ring $R$ has characteristic zero. In particular, $R$ injects into $R\otimes\Q$, so if we want to define a $\lambda$-ring structure on $R$, we can first define it on $R\otimes\Q$ as above and then show that the $\lambda$-operation preserves $R$.

\begin{lemma}\label{lm:1dim}
Let $R$ be a $\lambda$-ring and $r\in R$ be such that $\psi_i(r)=r^i$ for all $i$. Then

(i) $\lambda(r)=1+rz$;

(ii) for all $r'\in R$ we have $\lambda_n(rr')=r^n\lambda_n(r')$.
\end{lemma}
\begin{proof}
  (i) Since the maps $\Sym_\Z[x_\bullet]\to R\colon F\mapsto F[r]$ is a homomorphism, it is enough to prove a~similar statement in $\Sym_\Z[x_\bullet]$, where it follows from~\cite[(2.10')]{macdonald1998symmetric}.

  (ii) For any partition $\mu=(\mu_1\ge\mu_2\ge\ldots)$ define $e_\mu(x):=e_{\mu_1}[x]e_{\mu_2}[x]\ldots$. Also, recall symmetric functions $m_\mu$ (see, e.g.,~\cite[Ch.~1, Sect.~2]{macdonald1998symmetric}). The formal identity $\prod_{i,j}(1+x_iy_j)=\sum_\mu e_\mu(x)m_\mu(y)$ (see~\cite[(4.2')]{macdonald1998symmetric}) shows that for any $x,y\in R$ we have $e_n[xy]=\sum_{\mu\dashv n}e_\mu[x]m_\mu[y]$ (see the remark following~\cite[(4.2')]{macdonald1998symmetric}). However, by part~(i) we have $e_\mu[r]=0$ unless $\mu=1^n$ in which case $e_\mu[r]=r^n$. Thus, the only non-zero term in the sum
  \[
    \lambda_n(rr')=e_n[rr']=\sum_{\mu\dashv n}e_\mu[r]m_\mu[r']
  \]
  is the term $r^nm_{1^n}[r']=r^n\lambda_n(r')$.
\end{proof}

\subsection{Plethystic exponentials and logarithms}\label{sect:ExpLog} Let $R$ be a ring with an ideal $I$ such that $R$ is $I$-adically complete. Assume that $R$ is given a pre-$\lambda$-ring structure such that for all $n$ we have $h_n[I]\subset I^n$, where $h_n$ is the full symmetric polynomial of degree $n$. Equivalently, $\lambda_n(I)\subset I^n$ for all $n$. In this case we define the plethystic exponential $\Exp(r):=\sum_{n=0}^\infty h_n[r]$ whenever $r\in I$ (cf.~\cite[Sect.~2]{MellitIntegrality}). It is easy to check that $\Exp$ is an isomorphism of abelian groups $\Exp\colon I\xrightarrow{\sim}(1+I)^\times$. We denote by $\Log\colon(1+I)^\times\to I$ the inverse isomorphism.

In particular, if $A$ is any pre-$\lambda$-ring, then we can take $R=A[[z]]$, $R=A[[\Gamma_D]]$, or $R=A[[\Gamma_D,z]]$ with the induced pre-$\lambda$-ring structure. Then $I$ consists of formal power series with zero constant terms (denoted by $A[[z]]^0$, $A[[\Gamma_D]]^0$, and $A[[\Gamma_D,z]]^0$ respectively). This setup was applied in Sections~\ref{sect:Plethystic} and~\ref{sect:FullyFixed}, when $A$ is a ring of motivic classes.

\subsection{$L$-functions and $\lambda$-rings}\label{sect:Rg}
Recall that we denote the symmetric group of permutations of $n$ elements by $\Symm_n$.
Let $g\ge1$ be an integer. Consider the ring of Laurent polynomials in $g+1$ variables $\Z[q^{\pm1},\alpha_1^{\pm1},\ldots\alpha_g^{\pm1}]$. The semi-direct product $W_g:=\Symm_g\rtimes(\Z/2\Z)^g$ acts on this ring as follows: $W_g$ fixes $q$, $\Symm_g$ permutes $\alpha_i$, while the $i$-th copy of $\Z/2\Z$ maps $\alpha_i$ to $q\alpha_i^{-1}$ and $\alpha_j$ to itself whenever $j\ne i$. We denote by $R_g$ the ring of invariants:
\[
    R_g:=\Z[q^{\pm1},\alpha_1^{\pm1},\ldots,\alpha_g^{\pm1}]^{W_g}.
\]
It will be convenient for us to set $\alpha_{i+g}:=q\alpha_i^{-1}$ for $i=1,\ldots,g$. We define a \emph{universal $L$-function}
\[
   L^{univ}:=\prod_{i=1}^{2g}(1-\alpha_iz)=\prod_{i=1}^g(1-\alpha_iz)\prod_{i=1}^g(1-q\alpha_i^{-1}z)\in R_g[z].
\]
We note the ``functional equation'': $L^{univ}(q^{-1}z^{-1})=q^{-g}z^{-2g}L^{univ}(z)$.
\begin{proposition}\label{pr:GeneratorsOfRg}
(i) As a $\Z[q,q^{-1}]$-algebra, $R_g$ is generated by the coefficients of $L^{univ}$, the ideal of relations is generated by the coefficients of the functional equation together with the relation that the constant term of $L^{univ}$ is 1.

(ii) The ring $\Z[q^{\pm1},\alpha_1^{\pm1},\ldots,\alpha_g^{\pm1}]$ has a unique structure of a $\lambda$-ring given by $\psi_n(q^{\pm1})=q^{\pm n}$, $\psi_n(\alpha_i^{\pm1})=\alpha_i^{\pm n}$. Moreover, this structure preserves $R_g\subset\Z[q^{\pm1},\alpha_1^{\pm1},\ldots,\alpha_g^{\pm1}]$.

(iii) We have $\lambda(\alpha_1+\ldots+\alpha_{2g})=L^{univ}(-z)$. Thus $R_g$ is generated by $q$, $q^{-1}$, and $\alpha_1+\ldots+\alpha_{2g}$ as a $\lambda$-ring.
\end{proposition}
\begin{proof}
(i) For $i=0,\ldots,2g$, denote by $b_i$ the $i$-th coefficient of $L^{univ}$. Then the functional equations reads $b_{2g-i}=q^{g-i}b_i$. Thus, the statement can be reformulated as follows: $R_g$ is freely generated by $b_1$, \ldots, $b_g$ as a $\Z[q^{\pm1}]$-algebra.

Consider any $f\in R_g\subset\Z[q^{\pm1},\alpha_1^{\pm1},\ldots,\alpha_g^{\pm1}]$. Since $f$ is invariant under the substitution $\alpha_g\mapsto q\alpha_g^{-1}$, it can be uniquely written as $P(\alpha_g+q\alpha_g^{-1})$, where $P$ is a polynomial with coefficients in $\Z[q^{\pm1},\alpha_1^{\pm1},\ldots,\alpha_{g-1}^{\pm1}]$. Repeating this procedure, we can uniquely write
$f=P(x_1,\ldots,x_g)$, where $x_i=\alpha_i+q\alpha_i^{-1}$ and $P$ is a polynomial with coefficients in $\Z[q^{\pm1}]$. Moreover, this polynomial is symmetric, so we can uniquely write
\[
    f=P(x_1,\ldots,x_g)=Q(y_1,\ldots,y_g),
\]
where $y_i=e_i(x_1,\ldots,x_g)$ are the elementary symmetric polynomials applied to $x_i$. We have shown that $y_i$ freely generate $R_g$ as a $\Z[q^{\pm1}]$-algebra. We have
\[
    L^{univ}=(qz^2-x_1z+1)\ldots(qz^2-x_gz+1),
\]
so it is easy to see that for $i=1,\ldots,g$, the polynomial $b_i$ is the sum of $(-1)^iy_i$ and a symmetric polynomial in $x_1$, \ldots, $x_g$ of degree smaller than $i$. It follows that $b_i$ with $i=1,\ldots,g$ also freely generate $R_g$. This proves part~(i).

Next, from the description of $\lambda$-rings containing $\Q$ given in Section~\ref{Sect:SpecialRings}, we see that there is a~unique structure of a $\lambda$-ring on $\Q[q^{\pm1},\alpha_1^{\pm1},\ldots,\alpha_g^{\pm1}]$ such that $\psi_n(q^{\pm1})=q^{\pm n}$, $\psi_n(\alpha_i^{\pm1})=\alpha_i^{\pm n}$ for all $i$ and~$n$. By Lemma~\ref{lm:1dim}(i) for every monomial $m=q^n\alpha_1^{k_1}\ldots\alpha_g^{k_g}$, where $n,k_1,\ldots,k_g\in\Z$, we have $\lambda(m)=1+mz$. Thus the $\lambda$-operation preserves $\Z[q^{\pm1},\alpha_1^{\pm1},\ldots,\alpha_g^{\pm1}]$. This proves the first part of~(ii). On the other hand, the $\psi_n$-operations commute with the action of $W_g$, so they preserve $R_g\otimes\Q$. Thus the
$\lambda$-operation preserves
\[
    R_g=(R_g\otimes\Q)\cap\Z[q^{\pm1},\alpha_1^{\pm1},\ldots,\alpha_g^{\pm1}].
\]
This proves part~(ii). To prove part~(iii), we note that by Lemma~\ref{lm:1dim}(i) we have $\lambda(\alpha_i)=1+\alpha_iz$. Therefore
\[
    \lambda(\alpha_1+\ldots+\alpha_{2g})=\prod_i\lambda(\alpha_i)=L^{univ}(-z).
\]
\end{proof}

\begin{remark}
There is a slightly different way to define $R_g$ (cf.~\cite[Sect.~1.1]{SchiffmannIndecomposable}). Let $(V,\omega)$ be a symplectic vector space of dimension $2g$ over a field $K$. Let $\GSp(V)$ be the subgroup of $\GL(V)$ consisting of operators $A$ such that $\omega(Av_1,Av_2)=c\cdot\omega(v_1,v_2)$ for any $v_1,v_2\in V$ and some $c=c(A)\in K$ (a.k.a similitude group of $\omega$). Let $T\subset\GSp(V)$ be the standard maximal torus. Consider the ring of Laurent polynomials in $2g$ variables $\Z[\alpha_1^{\pm1},\ldots\alpha_{2g}^{\pm1}]$. Then its quotient by the ideal generated by the relations $\alpha_i\alpha_{i+g}=\alpha_j\alpha_{j+g}, 1\le i,j\le g$ is identified with the group ring $\Z[X^*(T)]$ of the character lattice of $T$. Next, the group $W_g=\Symm_g\rtimes(\Z/2\Z)^g$ is identified with the Weyl group of $\GSp(V)$. It acts on $X^*(T)$ as follows: $\Symm_g$ permutes pairs $(\alpha_i,\alpha_{i+g})$, while the $i$-th copy of $\Z/2\Z$ switches $\alpha_i$ and $\alpha_{i+g}$. Then $R_g$ is canonically isomorphic to $\Z[X^*(T)]^{W_g}$.

For example, if our base field $\kk$ is equal to $\C$, we can take $K=\C$ and $V=H^1(X,\C)$. If $\kk$ is algebraically closed, we can take $K=\Q_l$ and $V=H^1_{et}(X,\Q_l)$. If $\kk$ is an algebraic closure of a~finite field (which is not allowed in our setup), then $\alpha_i$ should be viewed as the Frobenius eigenvalues (see~\cite[Sect.~1.2]{schiffmann2018kac}).
\end{remark}

\subsection{Localizations of $\lambda$-rings}
If $R$ is a $\lambda$-ring, then the polynomial ring $R[t]$ becomes a $\lambda$-ring, where the $\lambda$-operation is given by $\lambda(\sum_i{a_it^i})(z)=\prod_i\lambda(a_i)(t^iz)\in R[t][[z]]$ (see~\cite[Ch.~I, Sect.~2]{knutson1973rings}). Similarly, we can extend the $\lambda$-ring structure to $R[[t]]$. We also get an induced $\lambda$-ring structure on the ring of symmetric functions $\Sym_R(w_\bullet)$. We now explain how to extend the $\lambda$-ring structures to certain localizations.

\begin{lemma}\label{lm:lambda_r'r}
  Let $R$ be a $\lambda$-ring, and $r\in R$ be such that for all $i>0$ we have $\psi_i(r)=r^i$ and $r^i\ne1$. Then

  (i) we can uniquely extend the $\lambda$-ring structure to the localization $S^{-1}R$, where the multiplicative subset $S$ is generated by all elements of the form $r^i-1$;

  (ii) let $R\to R'$ be a homomorphism of $\lambda$-rings such that the image of $S$ is contained in the set $(R')^\times$ of invertible elements of $R'$. Then the unique extension of this homomorphism to a homomorphism $S^{-1}R\to R'$ is a homomorphism of $\lambda$-rings.
\end{lemma}
\begin{proof}

  (i) The homomorphisms $\psi_i\otimes\Id_\Q\colon R\otimes\Q\to R\otimes\Q$ extend to $S^{-1}R\otimes\Q$ because we have $(\psi_i\otimes\Id_\Q)(r^j-1)=r^{ij}-1$, so $\psi_i\otimes\Id_\Q$ takes the multiplicative set to itself. Thus, we get a $\lambda$-ring structure on $S^{-1}R\otimes\Q$. We show that $S^{-1}R$ is preserved by the $\lambda$-operations. Consider $y\in S^{-1}R$. An inductive argument shows that we may assume that $y=x/(r^i-1)$, where $\lambda(x)\in S^{-1}R[[z]]$. Set $f(z):=\lambda(y)\in(S^{-1}R\otimes\Q)[[z]]$. We have $r^iy-y=x$, so by Lemma~\ref{lm:1dim}(ii) we obtain $f(r^iz)/f(z)=\lambda(x)\in S^{-1}R[[z]]$. However, since $r^j\ne1$ for all $j$, for every $g(z)\in S^{-1}R[[z]]$ there is a unique $f(z)\in S^{-1}R[[z]]$ such that $f(r^iz)/f(z)=g(z)$. We see that $f(z)\in S^{-1}R[[z]]$ so that the $\lambda$-operation preserves $S^{-1}R$. The uniqueness also follows.

  (ii) This is enough to prove after the tensorization with $\Q$. Thus, we may assume that both $R$ and $R'$ contain $\Q$, and we only need to check that the morphism $S^{-1}R\to R'$ commutes with the $\psi$-operations. This follows because a homomorphism from $S^{-1}R$ is fully determined by its restriction to $R$.
\end{proof}

In the next corollary we define a $\lambda$-ring $K_g$, which is an integral form of the ring $K_g$ defined in~\cite[Sect.~1.1]{SchiffmannIndecomposable} (that is, Schiffmann's $K_g$ corresponds to $K_g\otimes\Q$ in our notation). We also define the $\lambda$-ring~$V_g$, which is a localization of $K_g[z]$.
\begin{corollary}\label{cor:Kg}
(i) Let $K_g$ be the localization of $R_g$ with respect to the multiplicative set generated by $q^i-1$, where $i\ge1$. Then there is a unique extension of the $\lambda$-ring structure on $R_g$ to $K_g$.

(ii) Let $V_g$ be the localization of $R_g[z]$ with respect to the multiplicative set generated by $q^i-z^j$, where $i,j\ge0$, $(i,j)\ne(0,0)$. Then there is a unique extension of the $\lambda$-ring structure on $R_g[z]$ to $V_g$.
\end{corollary}
\begin{proof}
  In part~(i) we apply Lemma~\ref{lm:lambda_r'r}(i) with $r=q$. In part~(ii) we note that $V_g$ is the localization with respect to the multiplicative set generated by $z^jq^{-i}-1$, where $i,j\ge0$, $(i,j)\ne(0,0)$. Since $\psi_n(z^jq^{-i}-1)=z^{nj}q^{-ni}-1$, we can apply Lemma~\ref{lm:lambda_r'r}(i) iteratively with $r=z^jq^{-i}$, where $(i,j)$ range over all pairs of coprime nonnegative integers.
\end{proof}

\subsection{Homomorphisms from universal $\lambda$-rings to rings of motivic classes}\label{sect:MotHomo} In Section~\ref{sect:MotClasses} we recalled the ring $\Mot(\kk)$ and its dimensional completion $\cMot(\kk)$ (see~\cite[Sect.~1.1]{FedorovSoibelmans} for precise definitions). We now recall the pre-$\lambda$-ring structure on $\cMot(\kk)$. For a quasi-projective $\kk$-variety $Y$ we define, following~\cite{KapranovMotivic}, its motivic $\zeta$-function by
\[
    \zeta_Y(z):=\sum_i[Y^i/\Symm_i]z^i\in\cMot(\kk)[[z]],
\]
where $\Symm_i$ is the permutation group. Every class $A\in\cMot(\kk)$ can be written as the limit of a sequence $([Y_i]-[Z_i])/\bL^{n_i}$, where $Y_i,Z_i$ are quasi-projective $\kk$-varieties. Then
\[
    \lambda(A):=\lim_{i\to\infty}\frac{\zeta_{Z_i}(-\bL^{n_i}z)}{\zeta_{Y_i}(-\bL^{n_i}z)}.
\]
We note that this pre-$\lambda$-ring structure is opposite to that defined in~\cite[Sect.~1.3.1]{FedorovSoibelmans}. The reason for working with the opposite structure, is that, by definition, in any $\lambda$-ring we have $\lambda(1)=1+z$. Thus, since in the pre-$\lambda$-ring structure of~\cite[Sect.~1.3.1]{FedorovSoibelmans} we have $\lambda(1)=(1-z)^{-1}$, this pre-$\lambda$-ring has no non-trivial pre-$\lambda$-ring homomorphisms to a $\lambda$-ring. On the other hand, \emph{our\/} pre-$\lambda$-ring has a non-trivial $\lambda$-ring quotient, as we are going to see in Remark~\ref{rm:sMot_nontrivial}(i). We denote the universal quotient of $\cMot(\kk)$ that is a $\lambda$-ring by $\sMot(\kk)$, see~\cite{LarsenLuntsRational}.

Recall that in Section~\ref{sect:IntroHLV} we defined the motivic $L$-function of a smooth projective geometrically connected $\kk$-curve $X$ by
\[
    L_X(z)=\zeta_X(z)(1-z)(1-\bL z)\in\Mot(\kk)[z].
\]

\begin{proposition}\label{pr:HomFromUniv}
  (i) Let $X$ be as above. Then there is a unique ring homomorphism $K_g\to\Mot(\kk)$, sending $q^{\pm1}$ to $\bL^{\pm1}$ and the coefficient of $L^{univ}$ at $z^i$ to the coefficient of $L_X$ at $z^i$ for $i=1,\ldots,2g$.

  (ii) Moreover, the  composition of the projection $\Mot(\kk)\to\sMot(\kk)$ with the homomorphism of part~(i) is a homomorphism of $\lambda$-rings.

  (iii) The homomorphism from part (ii) extends uniquely to a homomorphism of $\lambda$-rings $V_g\to\sMot(\kk)[[z]]$ sending $z\in V_g$ to $z\in\sMot(\kk)[[z]]$.
\end{proposition}
\begin{proof}
  (i) Since $\bL^i-1$ is invertible in $\Mot(\kk)$, it is enough to prove a similar statement with $R_g$ instead of $K_g$. In view of Proposition~\ref{pr:GeneratorsOfRg}(i), we just need to show that $L_X$ has constant term~1 and satisfies the functional equation. The former is obvious, the latter is an easy consequence of~\cite[Thm.~1.1.9]{KapranovMotivic} (cf.~\cite[Prop.~1.3.1(iii)]{FedorovSoibelmans}).

  (ii) By Lemma~\ref{lm:lambda_r'r}(ii) it is enough to prove a similar statement with $R_g$ instead of $K_g$. For a~ring homomorphism of $\lambda$-rings, to check that it is a $\lambda$-ring homomorphism, we just need to check the corresponding identity on the generators. In view of Proposition~\ref{pr:GeneratorsOfRg}(iii), we only need to check the following:
  \[
  \begin{split}
    \lambda(q)=1+qz&\mapsto1+Lz=\lambda(\bL),\\
    \lambda(\alpha_1+\ldots+\alpha_{2g})=L^{univ}(-z)&\mapsto
    L_X(-z)=\lambda(-[X]+1+\bL).
  \end{split}
  \]
  These identities are clear.

  (iii) Note that $1-z^j\bL^{-i}$ is invertible in $\sMot(\kk)[[z]]$ when $j>0$: the inverse is $\sum_{n=0}^\infty z^{jn}\bL^{-jn}$. It is also invertible when $j=0$, $i\ne0$. Thus, the existence and uniqueness of the homomorphism follow. To check that it is a homomorphism of $\lambda$-rings, we apply Lemma~\ref{lm:lambda_r'r}(ii).
\end{proof}
\begin{remark}\label{rem:NotHomo}
   Note that $\Mot(\kk)$ has a pre-$\lambda$-ring structure, so one can ask whether the homomorphism of part~(i) is a homomorphisms of pre-$\lambda$-rings. However, to check that a homomorphism of rings is a homomorphisms of pre-$\lambda$-rings, it is \emph{not enough\/} to check this on the generators. In fact, we expect that the homomorphism $R_g\to\Mot(\kk)$ is \emph{not\/} a homomorphism of pre-$\lambda$-rings.
\end{remark}

\subsection{Universal generating function}
Recall that in Section~\ref{sect:IntroHLV} we defined the motivic generating function~\eqref{eq:HLV}. Now we define its universal version. Let $g,\delta\ge0$ be integers and $D$ be a finite set. Define
\begin{equation}\label{eq:HLV_Univ}
    \Omega^{univ}_{g,D,\delta}:=\sum_{\mu\in\cP}w^{|\mu|}(-q^\frac12)^{(2g+\delta)\langle\mu,\mu\rangle}
    z^{2\delta n(\mu')}
    \prod_{\Box\in\mu}\frac{L^{univ}(z^{2a(\Box)+1}q^{-l(\Box)-1})}
    {(z^{2a(\Box)+2}-q^{l(\Box)})(z^{2a(\Box)}-q^{l(\Box)+1})}
    \prod_{x\in D}\tilde H_\mu(w_{x,\bullet};z^2,q).
\end{equation}
Here $\cP$, $|\mu|$, $\langle\mu,\mu\rangle$, $\mu'$, $n(\mu')$, $a(\Box)$, $l(\Box)$ are as in Section~\ref{sect:IntroHLV}, while $\tilde H_\mu(w_{x,\bullet};z^2,q)$ are modified Macdonald polynomials, cf.~Section~\ref{sect:Macdonald}. We have
\[
    \Omega^{univ}_{g,D,\delta}\in\Sym_{V_g}(w_{\bullet,\bullet})[q^{\frac12}][[w]],
\]
where $V_g$ was defined in Corollary~\ref{cor:Kg}(ii). Next, the homomorphism of Proposition~\ref{pr:HomFromUniv}(iii), induces a homomorphism of $\lambda$-rings
\begin{equation}\label{eq:UnivHomo}
    \Sym_{V_g}(w_{\bullet,\bullet})[q^{\frac12}][[w]]\to
    \Sym_{\sMot(\kk)[[z]]}(w_{\bullet,\bullet})[\bL^{\frac12}][[w]];
\end{equation}
the image of $\Omega^{univ}_{g,D,\delta}$ under this homomorphism is $\Omega^{mot}_{X,D,\delta}$, where $g=g(X)$.

\subsection{Comparing the universal generating function with the universal Schiffmann's generating function}\label{sect:Schiffmann}
We now define another generating function. Recall that in~\cite[Sect.~1.4]{SchiffmannIndecomposable}, Schiffmann defined rational functions $J_\mu,H_\mu\in V_g$ indexed by partitions. We define
\[
    \Omega^{Sch,univ}_{g,D,\delta}:=\sum_{\mu\in\cP}w^{|\mu|}(-q^\frac12)^{(2g-2+\delta)\langle\mu,\mu\rangle}
    z^{\delta n(\mu')}
    J_\mu H_\mu\prod_{x\in D}\tilde H_\mu(w_{x,\bullet};z,q)\in\Sym_{V_g}(w_{\bullet,\bullet})[q^{\frac12}][[w]].
\]
Recall that in Section~\ref{sect:TwistedPairs} we introduced $J_{\mu,X}^{mot}(z), H_{\mu,X}^{mot}(z)\in\cMot(\kk)[[z]]$ (defined in~\cite[Sect.~1.3.2]{FedorovSoibelmans}). It is clear from the formulas that the images of $J_{\mu,X}^{mot}(z)$ and $H_{\mu,X}^{mot}(z)$ in $\sMot(\kk)[[z]]$ are equal to the images of $J_\mu$ and $H_\mu$ under the homomorphism of Proposition~\ref{pr:HomFromUniv}(iii). Consequently, the image of $\Omega^{Sch,mot}_{X,D,\delta}$ defined in~\eqref{eq:OmegaSchiffmann} in $\Sym_{\sMot(\kk)[[z]]}(w_{\bullet,\bullet})[\bL^{\frac12}][[w]]$ is equal to the image of $\Omega^{Sch,univ}_{g,D,\delta}$ under the homomorphism~\eqref{eq:UnivHomo}. We have a plethystic exponent
\[
 \Exp\colon\Sym_{V_g}(w_{\bullet,\bullet})[q^{\frac12}][[w]]^0\to
 1+\Sym_{V_g}(w_{\bullet,\bullet})[q^{\frac12}][[w]]^0\colon
 \sum_{\gamma,d} A_\gamma q^{n_\gamma}w^\gamma\mapsto\prod_{\gamma,d}\lambda^{-1}_{A_\gamma}(-q^{n_\gamma}w^\gamma),
\]
where the superscript ``$^0$'' stands for series with zero constant term. Comparing this with~\eqref{eq:ExpLog}, we see that the homomorphism~\eqref{eq:UnivHomo} commutes with taking plethystic exponents. Next, $\Exp$ is a~bijection, and we denote the inverse bijection by $\Log$.

\begin{theorem}\label{th:MellitIdentity}
Set $\HH^{univ}_{g,D,\delta}:=(1-z^2)\Log(\Omega^{univ}_{g,D,\delta})$ and $\HH^{Sch,univ}_{g,D,\delta}:=(1-z)\Log(\Omega^{Sch,univ}_{g,D,\delta})$. Then

(i) for $\gamma\in\Gamma_D$, the coefficients of $\HH^{univ}_{g,D,\delta}$ and of $\HH^{Sch,univ}_{g,D,\delta}$ at $w^\gamma$ are Laurent polynomials in $z$. More precisely, they belong to $K_g[q^\frac12,z^{\pm1}]$;

(ii) we have
\[
    \HH^{univ}_{g,D,\delta}\Biggl|_{z=1}=\HH^{Sch,univ}_{g,D,\delta}\Biggl|_{z=1}.
\]
\end{theorem}
\begin{proof}
  Note that $K_g\subset\Q(q)[\alpha_1^{\pm1},\ldots,\alpha_g^{\pm1}]$, so it is enough to prove the statement in this bigger ring. We introduce some notation in order to re-write $\Omega^{Sch,univ}_{g,D,\delta}$, cf.~\cite[(4)]{MellitNoPunctures}. For a partition $\mu$ set
  \[
    N_\mu(u)=N_\mu(u;z,q):=\prod_{\Box\in\mu}(z^{a(\Box)}-uq^{1+l(\Box)})(z^{a(\Box)+1}-u^{-1}q^{l(\Box)}).
  \]
  Next, following~\cite[(5)]{MellitNoPunctures}, set
  \begin{multline*}
    f(z_1,\ldots,z_n)
    \\:=\prod_i\prod_{k=1}^g\frac{1-\alpha_k^{-1}}{1-\alpha_k^{-1}z_i}
    \sum_{\sigma\in\Symm_n}\sigma\left\{
    \prod_{i>j}\left(
    \frac1{1-\frac{z_i}{z_j}}
    \prod_{k=1}^g\frac{1-\alpha_k^{-1}\frac{z_i}{z_j}}{1-q\alpha_k^{-1}\frac{z_i}{z_j}}
    \right)
    \prod_{i>j+1}
    \left(1-q\frac{z_i}{z_j}\right)\prod_{i\ge2}(1-z_i)
    \right\}.
  \end{multline*}
  When $\mu$ is a partition, set $f_\mu:=f(z_1(\mu),\ldots,z_{\ell(\mu)}(\mu))$, where $z_i(\mu):=q^{-\ell(\mu)+i}z^{\mu_i}$, $\ell(\mu)$ is the length of $\mu$. By~\cite[Prop.~3.1]{MellitNoPunctures} (see also~\cite[Prop.~3.4]{OGormanMozgovoy}), we obtain
  \[
    \Omega^{Sch,univ}_{g,D,\delta}=
    \sum_{\mu\in\cP}w^{|\mu|}(-q^\frac12)^{\delta\langle\mu,\mu\rangle}
    z^{\delta n(\mu')}
    f_\mu\frac{\prod_{i=1}^gN_\mu(\alpha_i^{-1})}{N_\mu(1)}
    \prod_{x\in D}\tilde H_\mu(w_{x,\bullet};z,q),
  \]
  Substituting $w$ with $wq^{-\delta/2}$, we obtain:
  \[
    \Omega^{Sch,univ}_{g,D,\delta}(wq^{-\delta/2})=
    \sum_{\mu\in\cP}w^{|\mu|}(-1)^{\delta|\mu|}q^{\delta n(\mu)}
    z^{\delta n(\mu')}
    f_\mu\frac{\prod_{i=1}^gN_\mu(\alpha_i^{-1})}{N_\mu(1)}
    \prod_{x\in D}\tilde H_\mu(w_{x,\bullet};z,q).
  \]
  We will need one more series, which is obtained from the previous one by replacing in the previous series $f_\mu$ with 1:
  \begin{equation}\label{eq:Sch2}
    \Omega^{HLV,univ}_{g,D,\delta}:=
    \sum_{\mu\in\cP}w^{|\mu|}(-1)^{\delta|\mu|}q^{\delta n(\mu)}
    z^{\delta n(\mu')}
    \frac{\prod_{i=1}^gN_\mu(\alpha_i^{-1})}{N_\mu(1)}
    \prod_{x\in D}\tilde H_\mu(w_{x,\bullet};z,q).
  \end{equation}
  This series is an irregular version of the HLV generating function introduced in~\cite[Sect.~7.1]{MellitPunctures}. We emphasize that the coefficients of this series belong to $\Q(z,q)[\alpha_1^{\pm1},\ldots,\alpha_g^{\pm1}]^{\Symm_g}$ rather than to a smaller ring $R_g\otimes\Q(z,q)=\Q(z,q)[\alpha_1^{\pm1},\ldots,\alpha_g^{\pm1}]^{W_g}$ (in other words, the series is not $(\Z/2\Z)^g$-invariant). Thus, we cannot apply the homomorphism of Proposition~\ref{pr:HomFromUniv}. As a result, the series has no motivic interpretation. To relate the two series, we need a lemma.
\begin{lemma}\label{lm:Admissible} The series $\Omega^{HLV,univ}_{g,D,\delta}$ is \emph{admissible\/} in the sense that
\[
    (q-1)(1-z)\Log\Omega^{HLV,univ}_{g,D,\delta}\in
    \Q[z,q,\alpha_i^{\pm1},w_{\bullet,\bullet}][[w]].
\]
\end{lemma}
\begin{proof}
  When $\delta=0$ this follows from~\cite[Thm.~7.1]{MellitIntegrality} (where $q$, $t$, and $u_i$ should be replaced with $z$, $q$, and $\alpha_i^{-1}$ respectively). The case $\delta>0$ is proved as follows. Recall the operator $\nabla$ on $\Sym_{\Q[z,q]}[w_\bullet]$ defined by
  \[
    \nabla\tilde H_\mu=(-1)^{|\mu|}q^{n(\mu)}z^{n(\mu')}\tilde H_\mu,
  \]
  see~\cite[(15)]{MellitIntegrality}. Then $\nabla$ is admissible by~\cite[Cor.~6.3]{MellitIntegrality}. However, $\Omega^{HLV,univ}_{g,D,\delta}$ is obtained from $\Omega^{HLV,univ}_{g,D,0}$ by applying $\nabla^\delta$ to any of the sequences of variables $w_{x,\bullet}$. This operation preserves admissibility because $\nabla$ is admissible. We are done.
\end{proof}

We return to the proof of Theorem~\ref{th:MellitIdentity}. We would like to apply~\cite[Lemma~5.1]{MellitNoPunctures} with $R=\Q(q^{\frac12})[z^{\pm1},\alpha_1^{\pm1},\ldots,\alpha_g^{\pm1}]$, and~\eqref{eq:Sch2} in place of $\Omega[X]$. In fact, we need a multipoint version of this lemma, where the modified Macdonald polynomials are replaced by products of modified Macdonald polynomials. We note that the multipoint version is implicit in the omitted proof of~\cite[Thm.~7.2]{MellitPunctures}. From this, we see that for all $\gamma\in\Gamma_D$ the coefficients of $(1-z)\Log\Omega^{HLV,univ}_{g,D,\delta}$ and of
$(1-z)\Log\Omega^{Sch,univ}_{g,D,\delta}(wq^{-\delta/2})$ at $w^\gamma$ belong to $R$ and that
\begin{equation}\label{eq:Mellit}
    \Biggr((1-z)\Log\Omega^{Sch,univ}_{g,D,\delta}(wq^{-\delta/2})\Biggl)\Biggl|_{z=1}=
    \Biggr((1-z)\Log\Omega^{HLV,univ}_{g,D,\delta}\Biggl)\Biggl|_{z=1}
\end{equation}
(Cf.~the proof of~\cite[Thm.~4.4]{OGormanMozgovoy}.) Now we see that the coefficients of $(1-z)\Log\Omega^{Sch,univ}_{g,D,\delta}$ also belong to $R$ and thus to $K_g[q^\frac12,z^{\pm1}]=R\cup K_g[q^\frac12](z)$. This proves the second half of part~(i).

Finally, we relate $\Omega^{univ}_{g,D,\delta}$ with $\Omega^{HLV,univ}_{g,D,\delta}$:
\begin{lemma}\label{lm:SubstIrr}
    \[
        \Omega^{univ}_{g,D,\delta}(\alpha_1,\ldots,\alpha_g,z,q,w_{\bullet,\bullet},wq^{-\delta/2})=
        \Omega^{HLV,univ}_{g,D,\delta}((\alpha_1z),\ldots,(\alpha_gz),z^2,q,w_{\bullet,\bullet},w).
    \]
\end{lemma}
\begin{proof}
Recall that we have $\sum_{\Box\in\mu}(2l(\Box)+1)=\langle\mu,\mu\rangle=2n(\mu)+|\mu|$. The lemma follows from the identities
\begin{multline*}
    \prod_{i=1}^gN_\mu((\alpha_iz)^{-1};z^2,q)=
    \prod_{i=1}^g\prod_{\Box\in\mu}(z^{2a(\Box)}-(\alpha_iz)^{-1}q^{l(\Box)+1})(z^{2a(\Box)+2}-(\alpha_iz)q^{l(\Box)})\\=
    \prod_{i=1}^g\prod_{\Box\in\mu}(\alpha_i z^{2a(\Box)+1}-q^{l(\Box)+1})(\alpha_i^{-1}z^{2a(\Box)+1}-q^{l(\Box)})\\=
    q^{g\langle\mu,\mu\rangle}\prod_{\Box\in\mu}L^{univ}(\alpha_1,\ldots,\alpha_g,z^{2a(\Box)+1}q^{-l(\Box)-1}),
\end{multline*}
and
\[
    (-q^\frac12)^{(2g+\delta)\langle\mu,\mu\rangle}q^{-\frac\delta2|\mu|}=(-1)^{\delta|\mu|}q^{\delta n(\mu)}q^{g\langle\mu,\mu\rangle}.
\]
\end{proof}
We return to the proof of Theorem~\ref{th:MellitIdentity}. By Lemma~\ref{lm:Admissible}, the coefficients of $(1-z)\Log\Omega^{HLV,univ}_{g,D,\delta}$ belong to $R=\Q(q^{\frac12})[z^{\pm1},\alpha_1^{\pm1},\ldots,\alpha_g^{\pm1}]$.
Note that the substitutions $w\mapsto wq^{-\delta/2}$ and $z\mapsto z^2$ commute with $\Log$. Now by Lemma~\ref{lm:SubstIrr} the coefficients of $\HH^{univ}_{g,D,\delta}$ belong to $R$ and thus to $K_g[q^\frac12,z^{\pm1}]$. This completes the proof of part~(i) of the theorem. We also obtain
\[
    \HH^{univ}_{g,D,\delta}(wq^{-\delta/2})\Bigl|_{z=1}=
    \Biggl((1-z)\Log\Omega^{HLV,univ}_{g,D,\delta}\Biggl)\Biggl|_{z=1}.
\]
Combining this with~\eqref{eq:Mellit}, we complete the proof of Theorem~\ref{th:MellitIdentity}.
\end{proof}

Recall from~\eqref{eq:HH_mot} the series $\HH^{mot}_{X,D,\delta}=(1-z^2)\Log\Omega^{mot}_{X,D,\delta}$, and let $\HH^{Sch,mot}_{X,D,\delta}$ be the image of the series $(1-z)\Log\Omega^{Sch,mot}_{X,D,\delta}$ under the projection $\cMot(\kk)\to\sMot(\kk)$ ($\Omega^{Sch,mot}_{X,D,\delta}$ is defined by~\eqref{eq:OmegaSchiffmann}). Applying to the identity of Theorem~\ref{th:MellitIdentity} the homomorphism $K_g\to\sMot(\kk)$ from Proposition~\ref{pr:HomFromUniv}(ii), we get the following corollary.
\begin{corollary}\label{cor:MellitSimpl}
For all $\gamma\in\Gamma_D$, the coefficients of $\HH^{mot}_{X,D,\delta}$ and of $\HH^{Sch,mot}_{X,D,\delta}$ at $w^\gamma$ belong to the ring $\sMot(\kk)[\bL^\frac12,z^{\pm1}]$, and we have
\[
    \HH^{mot}_{X,D,\delta}\Biggl|_{z=1}=\HH^{Sch,mot}_{X,D,\delta}\Biggl|_{z=1}.
\]
\end{corollary}

\section{Proofs of main theorems} In this section we will prove Theorems~\ref{th:IntroMain},~\ref{th:IntroMain2}, and~\ref{th:IntroMain3} by simplifying the formulas obtained in Theorems~\ref{th:AnswerConnInMot} and~\ref{th:AnswerInCmot} using Corollary~\ref{cor:MellitSimpl}. This requires working in $\sMot(\kk)$. The curve~$X$ is fixed in this section, so we omit it from the notation for the stacks. Recall also that we have a subset $\Gamma_{\divisor',\zeta}\subset\Gamma_{D'}$ consisting of the classes $\gamma$ such that $\gamma$ is full at all $x\in D'$ with $n'_x\ge2$ and $\zeta$ is non-resonant for $\gamma$ at all such $x$.

\subsection{Stabilization of coefficients of Schiffmann's generating series}

\begin{proposition}\label{pr:DT_Stabilize}
  For all $\gamma\in\Gamma_D$ there is $d_0\in\Z$ such that the coefficients of $(1-z)^{-1}\HH^{Sch,mot}_{X,D,\delta}$ at $w^\gamma z^d$ are independent of $d$ as long as $d\ge d_0$. Moreover, these coefficients are equal to the coefficient at $w^\gamma$ in $\HH^{mot}_{X,D,\delta}\bigl|_{z=1}$.
\end{proposition}
\begin{proof}
  This is essentially a reformulation of Corollary~\ref{cor:MellitSimpl}. Indeed, keeping the notations of the corollary, write
  \[
        \frac{\HH^{Sch,mot}_{X,D,\delta}}{(1-z)}=
        \frac{\HH^{Sch,mot}_{X,D,\delta}-\HH^{Sch,mot}_{X,D,\delta}\Bigl|_{z=1}}{(1-z)}+
        \frac{\HH^{Sch,mot}_{X,D,\delta}\Bigl|_{z=1}}{(1-z)}.
  \]
  By Corollary~\ref{cor:MellitSimpl}, the coefficients of $\HH^{Sch,mot}_{X,D,\delta}$ are Laurent polynomials in $z$, so the coefficients of the first summand above are also polynomials in $z$. This means that for every $\gamma\in\Gamma_D$ there is $d_0$ such that the coefficient of the first summand at $w^\gamma z^d$ vanishes for $d>d_0$. According to Corollary~\ref{cor:MellitSimpl} the second summand is equal to
  \[
    \frac{\HH^{mot}_{X,D,\delta}|_{z=1}}{(1-z)}=\left(\HH^{mot}_{X,D,\delta}\Bigl|_{z=1}\right)\sum_{d=0}^\infty z^d.
  \]
  The statement follows.
\end{proof}

\begin{remark}
Similarly to~\cite[Sect.~1.4]{FedorovSoibelmansParabolic} and~\cite[(1)]{MozgovoySchiffmann2020}, one can call the coefficients of $\Log\Omega^{Sch,mot}_{X,D,\delta}$ the \emph{motivic Donaldson--Thomas invariants} of the moduli stacks of parabolic $\epsilon$-connections. According to Theorems~\ref{th:AnswerConnInMot} and~\ref{th:AnswerInCmot} the motivic classes of the stacks considered in this paper can be calculated in terms of these invariants.
\end{remark}

\subsection{Proof of Theorem~\ref{th:IntroMain}} We need some notation. Let $H_n\in\sMot(\kk)$, $n\ge0$, be a~sequence of motivic classes. We will write $\lim_{n\to\infty} H_n=H$ if there is $n_0$ such that $H_n=H$ whenever $n\ge n_0$. In other words, we consider discrete topology on $\sMot(\kk)$. More generally, let $H_n=\sum_{\gamma\in\Gamma_D}H_{\gamma,n}w^\gamma\in\sMot(\kk)[[\Gamma_D]]$ be a sequence of motivic series indexed by $\Z_{\ge0}$. We write $\lim_{n\to\infty}H_n=H$, where $H=\sum_{\gamma\in\Gamma_D}H_\gamma w^\gamma$, if for all $\gamma\in\Gamma_D$ there is $n_0\in\Z$ such that $H_{n,\gamma}=H_\gamma$ when $n\ge n_0$. In other words, we consider coefficientwise convergence, where coefficients are contained in a discrete topological space.

\begin{lemma}\label{lm:LimExp}
  Assume that $H_n\in\sMot(\kk)[[\Gamma_D]]^0$ and $\lim_{n\to\infty} H_n=H$. Then $\lim_{n\to\infty}\Exp H_n=\Exp H$.
\end{lemma}
\begin{proof}
  Note that a coefficient of $\Exp H_n$ at $w^\gamma$ depends only on coefficients of $H_n$ at $w_{\gamma'}$ such that $\gamma-\gamma'\in\Gamma_D$. However, for a fixed $\gamma$ there are only finitely many such $\gamma'$.
\end{proof}

We are now ready to prove our first theorem from the Introduction. Recall that $\divisor=\sum_{x\in D}{n_xx}$ is an effective divisor on $X$ with support $D\subset X(\kk)$, $\epsilon\in\kk$, $\zeta\in\FNF(\divisor)$ (see~\eqref{eq:FNF_full}) and we are studying moduli stacks $\Conn_{\gamma,d}(\epsilon,\divisor,\zeta)$ of parabolic $\epsilon$-connections defined in Definition~\ref{def:ModSpace}. Here $\gamma\in\Gamma_D$, $d\in\Z$.
\begin{theorem}[See Theorem~\ref{th:IntroMain}]\label{th:Main1}
Assume that $\epsilon\ne0$, $\zeta\in\FNF(\divisor)$, and $\gamma'\in\Gamma_{\divisor,\zeta}$. Assume that $d'\in\Z$ is such that $\epsilon d'+\gamma'\star\zeta'=0$. Then the motivic class of $\Conn_{\gamma',d'}(\epsilon,\divisor,\zeta)$ in $\sMot(\kk)$ is equal to the coefficient at~$w^{\gamma'}$ in
\[
    (-\bL^\frac12)^{\chi(\gamma')} \Exp\left(\bL\left(\HH^{mot}_{X,D,\delta}\Bigl|_{z=1}\right)_{\gamma\star\zeta\in\epsilon\Z}\right).
\]
\end{theorem}
\begin{proof} According to Theorem~\ref{th:AnswerConnInMot} for all $N\gg0$ the motivic class of $\Conn_{\gamma',d'}(\epsilon,\divisor,\zeta)$ in $\cMot(\kk)$ is equal to the coefficient of
\[
    (-\bL^\frac12)^{\chi(\gamma')}\Exp\left(\Bigl(\bL\Log\Omega^{Sch,mot}_{X,D,\delta}\Bigr)_{-\epsilon d+\gamma\star\zeta=-\epsilon N\rk\gamma}\right)
\]
at $w^{\gamma'} z^{-d'+N\rk\gamma'}$. Applying the homomorphism $\cMot(\kk)\to\sMot(\kk)$, by definition of $\HH^{Sch,mot}_{X,D,\delta}$, we see that for all $N\gg0$ the motivic class of $\Conn_{\gamma',d'}(\epsilon,\divisor,\zeta)$ in $\sMot(\kk)$ is equal to the coefficient of
\[
    (-\bL^\frac12)^{\chi(\gamma')}\Exp\left(\left(\bL\frac{\HH^{Sch,mot}_{X,D,\delta}}{1-z}\right)_{-\epsilon d+\gamma\star\zeta=-\epsilon N\rk\gamma}\right)
\]
at $w^{\gamma'} z^{-d'+N\rk\gamma'}$. We note that every non-zero monomial $A_{\gamma,d}w^\gamma z^d$ in the above series is such that $-\epsilon d+\gamma\star\zeta=-\epsilon N\rk\gamma$. Since $\epsilon\ne0$ and $\epsilon d'+N\rk\gamma'=0$, the only non-zero monomial of the form $Aw^{\gamma'}z^d$ is the monomial $[\Conn_{\gamma',d'}(\epsilon,\divisor,\zeta)]w^{\gamma'} z^{-d'+N\rk\gamma'}$.
We see that $[\Conn_{\gamma',d'}(\epsilon,\divisor,\zeta)]$ is the coefficient at $w^{\gamma'}$ in the series obtained by substituting $z=1$:
\[
    (-\bL^\frac12)^{\chi(\gamma')}\Exp\left(\left(\left(\bL\frac{\HH^{Sch,mot}_{X,D,\delta}}{1-z}\right)_{-\epsilon d+\gamma\star\zeta=-\epsilon N\rk\gamma}\right)_{z=1}\right).
\]
We have also used that $\Exp$ commutes with the substitution $z=1$.

\begin{lemma}\label{lm:Step}
\[
  \lim_{N\to\infty}\left(\left(\bL\frac{\HH^{Sch,mot}_{X,D,\delta}}{1-z}\right)_{-\epsilon d+\gamma\star\zeta=-\epsilon N\rk\gamma}\right)_{z=1}=
  \left(\bL\left(\HH^{mot}_{X,D,\delta}\right)_{z=1}\right)_{\gamma\star\zeta\in\epsilon\Z}.
\]
\end{lemma}
\begin{proof}
  Write
  \[
    \left(\left(\bL\frac{\HH^{Sch,mot}_{X,D,\delta}}{1-z}\right)_{-\epsilon d+\gamma\star\zeta=-\epsilon N\rk\gamma}\right)_{z=1}=\sum_{\gamma\in\Gamma_D}C_{\gamma,N}w^\gamma.
  \]
  and
  \[
    \left(\bL\left(\HH^{mot}_{X,D,\delta}\right)_{z=1}\right)_{\gamma\star\zeta\in\epsilon\Z}=
    \sum_{\gamma\in\Gamma_D}C_\gamma w^\gamma.
  \]
  If $\gamma\star\zeta\notin\epsilon\Z$, then $\lim_{N\to\infty}C_{\gamma,N}=\lim_{N\to\infty}0=0=C_\gamma$. If $\gamma\star\zeta\in\epsilon\Z$, then $C_{\gamma,N}$ is equal to the coefficient of $(1-z)^{-1}\bL\;\HH^{Sch,mot}_{X,D,\delta}$ at $w^\gamma z^{\epsilon^{-1}\gamma\star\zeta+N\rk\gamma}$. Thus, we obtain that $\lim_{N\to\infty}C_{\gamma,N}=C_\gamma$ from Proposition~\ref{pr:DT_Stabilize}.
\end{proof}
Combining this with Lemma~\ref{lm:LimExp}, we get the statement of the theorem.
\end{proof}

\subsection{Proof of Theorem~\ref{th:IntroMain2}} Recall that if $\sigma$ is a sequence of parabolic weights of type $(X,\divisor)$ (see Definition~\ref{def:ParWeights}), $\zeta\in\FNF(\divisor)$, and $\gamma\in\Gamma_{\divisor,\zeta}$, then we have open substacks $\Conn^{\sigma-ss}_{\gamma,d}(\epsilon,\divisor,\zeta)\subset\Conn_{\gamma,d}(\epsilon,\divisor,\zeta)$ classifying $\sigma$-semistable parabolic $\epsilon$-connections.
\begin{theorem}[see Theorem~\ref{th:IntroMain2}]\label{th:Main2}
Assume that $\zeta\in\FNF(\divisor)$ and $\gamma'\in\Gamma_{\divisor,\zeta}$. Assume that $d'\in\Z$ is such that $\epsilon d'+\gamma'\star\zeta=0$. Let $\sigma$ be a sequence of parabolic weights of type $(X,\divisor)$. Then the motivic class of $\Conn^{\sigma-ss}_{\gamma',d'}(\epsilon,\divisor,\zeta)$ in $\sMot(\kk)$ is equal to the coefficient at $w^{\gamma'}$ in
\[
    (-\bL^\frac12)^{\chi(\gamma')} \Exp\left(\bL\left(\HH^{mot}_{X,D,\delta}\Bigl|_{z=1}\right)_{\substack{\gamma\star\zeta\in\epsilon\Z\\ \gamma\star\sigma-\tau\rk\gamma\in\Z }}\right),
\]
where $\tau:=(d'+\gamma'\star\sigma)/\rk\gamma'$.
\end{theorem}
\begin{proof}
  Similarly to the proof of Theorem~\ref{th:Main1}, the statement is derived form Theorem~\ref{th:AnswerInCmot}(ii). Lemma~\ref{lm:Step} should be replaced with
  \[
    \lim_{N\to\infty}\left(\left(\bL\frac{\HH^{Sch,mot}_{X,D,\delta}}{1-z}\right)_%
    {\substack{-\epsilon d+\gamma\star\zeta=-\epsilon N\rk\gamma\\-d+\gamma\star\sigma=(\tau-N)\rk\gamma}}\right)_{z=1}=
    \left(\bL\left(\HH^{mot}_{X,D,\delta}\right)_{z=1}\right)_{\substack{\gamma\star\zeta\in\epsilon\Z\\ \gamma\star\sigma-\tau\rk\gamma\in\Z}}.
  \]
\end{proof}

\subsection{Proof of Theorem~\ref{th:IntroMain3}} Finally, consider the stacks $\Conn^{prtl,\sigma-ss}_{\gamma',d'}(\epsilon,\divisor,\divisor',\zeta)$ classifying $\sigma$-semistable parabolic $\epsilon$-connections with partially fixed normal forms. Here $\divisor'$ is any effective divisor such that $\divisor'<\divisor$, and we denote its support by $D'$.
\begin{theorem}[see Theorem~\ref{th:IntroMain3}]\label{th:Main3}
Assume that $\divisor'<\divisor$, that is, $n'_x<n_x$ for at least one $x\in D$. Assume that $\zeta\in\FNF(\divisor,\divisor')$ and that $\gamma'\in\Gamma_{\divisor',\zeta}$. Let $\sigma$ be a sequence of parabolic weights of type $(X,\divisor')$. Then the motivic class of $\Conn^{prtl,\sigma-ss}_{\gamma',d'}(\epsilon,\divisor,\divisor',\zeta)$ in $\sMot(\kk)$ is equal to the coefficient at $w^{\gamma'}$ in
\[
    (-\bL^\frac12)^{\chi(\gamma')} \Exp\left(\left(\HH^{mot}_{X,D',\delta}\Bigl|_{z=1}\right)_{\gamma\star\sigma-\tau\rk\gamma\in\Z}\right),
\]
where $\tau:=(d'+\gamma'\star\sigma)/\rk\gamma'$.
\end{theorem}
\begin{proof}
    Similarly to the proof of Theorem~\ref{th:Main1}, the statement is derived from Theorem~\ref{th:AnswerInCmot}(i). Lemma~\ref{lm:Step} should be replaced with
  \[
    \lim_{N\to\infty}\left(\left(\frac{\HH^{Sch,mot}_{X,D',\delta}}{1-z}\right)_%
    {-d+\gamma\star\sigma=(\tau-N)\rk\gamma}\right)_{z=1}=
    \left(\left(\HH^{mot}_{X,D',\delta}\right)_{z=1}\right)_%
    {\gamma\star\sigma-\tau\rk\gamma\in\Z}.
  \]
\end{proof}

\section{E-polynomials and virtual Poincar\'e polynomials of stacks of parabolic \texorpdfstring{$\epsilon$}{epsilon}-connections}\label{sect:Deligne polynomials}
In this section we use Theorems~\ref{th:Main1},~\ref{th:Main2}, and~\ref{th:Main3} to calculate the E-polynomials (also known as Hodge--Deligne polynomials) of moduli stacks of parabolic $\epsilon$-connections. Specializing E-polynomials, we obtain the virtual Poincar\'e polynomials. The main result is Corollary~\ref{cor:Poincare}. In Section~\ref{sect:Pantev}, we use this result to give evidence for a conjecture of~\cite{DiaconescuDonagiPantev}.

\subsection{E-polynomials of a stack} Let $\Z[u,v][[(uv)^{-1}]]$ denote the completion of $\Z[u^{\pm1},v^{\pm1}]$ with respect to the $(uv)^{-1}$-adic topology.

\begin{proposition}\label{pr:HodgeDeligne}
  (i) There is a unique continuous ring homomorphism
  \[
    E\colon\cMot(\kk)\to\Z[u,v][[(uv)^{-1}]]
  \]
  sending the class $[Y]$ of a smooth projective variety $Y$ to
  \begin{equation}\label{eq:HodgeDeligne}
    \sum_{p,q=0}^{\dim Y}(-1)^{p+q}h^{p,q}(Y)u^pv^q,
  \end{equation}
  where $h^{p,q}(Y)$ are the Hodge numbers of $Y$.

  (ii) This homomorphism is a homomorphism of pre-$\lambda$-rings.
\end{proposition}
It follows from part~(ii) that $E$ factors through the universal $\lambda$-ring quotient, so we obtain a~homomorphism of $\lambda$-rings $\sMot(\kk)\to\Z[u,v][[(uv)^{-1}]]$, which we also denote by $E$.
\begin{proof}
  (i) Let $\Mot_{var}(\kk)$ denote the Grothendieck group of varieties, that is, the group generated by symbols $[Y]$, where $Y$ is a $\kk$-variety with relations $[Y]=[Y-Y']+[Y']$ whenever $Y'$ is a closed subvariety of $Y$. This has a ring structure given by $[Y_1][Y_2]=[Y_1\times_\kk Y_2]$ (the product is a reduced scheme because $k$ has characteristic zero). We put $\bL=[\A^1_\kk]$ and let $\cMot_{var}(\kk)$ be the completion of $\Mot_{var}(\kk)[\bL^{-1}]$ with respect to the filtration $F^\bullet\Mot_{var}(\kk)[\bL^{-1}]$, where $F^m\Mot_{var}(\kk)[\bL^{-1}]$ is generated by $\bL^{-n}[Y]$ with $\dim Y-n\le-m$. We have a canonical isomorphism of topological rings $\cMot_{var}(\kk)\simeq\cMot(\kk)$, see~\cite[Sect.~2.5]{FedorovSoibelmans}.

  Now uniqueness of $E$ follows because smooth projective varieties additively generate $\Mot_{var}(\kk)$ (we are using resolution of singularities). For the existence recall the Hodge--Deligne homomorphism
  \[
    E\colon\Mot_{var}(\kk)\to\Z[[u^{-1},v^{-1}]],
  \]
  satisfying~\eqref{eq:HodgeDeligne}. For a smooth projective geometrically integral $\kk$-curve $X$ we have
  \begin{equation}
    E([X])=1-gu-gv+uv.
  \end{equation}
  so that
  \[
    E(\bL)=E([\P^1_\kk])-E([\Spec\kk])=1+uv-1=uv.
  \]
  Now it is clear that $E$ extends by continuity to $\cMot_{var}(\kk)\simeq\cMot(\kk)$.

  (ii) Follows from~\cite[Cor.~2]{MaximSaitoEtAl2011}.
\end{proof}
By Section~\ref{sect:ExpLog}, we have a plethystic exponential $\Exp\colon\Z[u,v][[(uv)^{-1},t]]^0\to1+\Z[u,v][[(uv)^{-1},t]]^0$. Since $E$ is a homomorphism of pre-$\lambda$-rings, we have
\[
    E(\zeta_X(t))=\Exp(tE[X])=\Exp(t(1-gu-gv+uv))=\frac{(1-ut)^g(1-vt)^g}{(1-t)(1-uvt)}
\]
and
\[
    E(L_X(t))=(1-ut)^g(1-vt)^g.
\]

\begin{remark}\label{rm:sMot_nontrivial}
(i) It follows from the above proposition that in $\sMot(\kk)$ we have $[\cX]\ne0$ whenever $\cX$ is a non-empty stack. In particular, $\sMot(\kk)\ne0$.

(ii) We have a surjective homomorphism $\Z[q^{\pm1},\alpha_1^{\pm1},\ldots,\alpha_g^{\pm1}]\to\Z[u^{\pm1},v^{\pm1}]$, sending~$q$ to $uv$, and $\alpha_i$ to $u$. Recall from Section~\ref{sect:Rg} the ring $R_g=\Z[q^{\pm1},\alpha_1^{\pm1},\ldots,\alpha_g^{\pm1}]^{W_g}$. Restricting the above homomorphism to $R_g$, we obtain the unique homomorphism $E_g\colon R_g\to\Z[u^{\pm1},v^{\pm1}][[(uv)^{-1}]]$ sending~$q$ to $uv$ and $L^{univ}$ to $(1-ut)^g(1-vt)^g$, where uniqueness follows from Proposition~\ref{pr:GeneratorsOfRg}(i). It is clear that the composition of the homomorphism of Proposition~\ref{pr:HomFromUniv}(ii) with $E$ is equal to $E_g$.
\end{remark}

\subsection{E-polynomials of moduli stacks of parabolic $\epsilon$-connections}\label{sect:EPoly} Applying the homomorphism $E$ to $\tilde H_\mu^{mot}(w_\bullet;z)$, we get $\tilde H_\mu(w_\bullet;uv,z)$, see Remark~\ref{rm:sMot_nontrivial}(ii). Thus, applying $E$ to the coefficients of the motivic generating function $\Omega^{mot}_{X,D,\delta}$ (see~\eqref{eq:HLV}), we get the following series:
\begin{multline}\label{eq:EPoly}
    \Omega^E_{g,D,\delta}:=E(\Omega^{mot}_{X,D,\delta})=\sum_{\mu\in\cP}w^{|\mu|}(-u^{1/2}v^{1/2})^{(2g+\delta)\langle\mu,\mu\rangle}z^{2\delta n(\mu')}
    \\ \times\prod_{\Box\in\mu}
    \frac{(1-z^{2a(\Box)+1}u^{-l(\Box)}v^{-l(\Box)-1})^g(1-z^{2a(\Box)+1}u^{-l(\Box)-1}v^{-l(\Box)})^g}
    {(z^{2a(\Box)+2}-(uv)^{l(\Box)})(z^{2a(\Box)}-(uv)^{l(\Box)+1})}
    \prod_{x\in D}\tilde H_\mu(w_{x,\bullet};z^2,uv)\\=
    \sum_{\mu\in\cP}w^{|\mu|}(u^{1/2}v^{1/2})^{\delta|\mu|}\prod_{\Box\in\mu}
    \frac{(-uv)^{\delta l(\Box)}z^{2\delta a(\Box)}(u^{l(\Box)}v^{l(\Box)+1}-z^{2a(\Box)+1})^g(u^{l(\Box)+1}v^{l(\Box)}-z^{2a(\Box)+1})^g}
    {(u^{l(\Box)+1}v^{l(\Box)+1}-z^{2a(\Box)})(u^{l(\Box)}v^{l(\Box)}-z^{2a(\Box)+2})}\\ \times\prod_{x\in D}\tilde H_\mu(w_{x,\bullet};z^2,uv).
\end{multline}
We note that $\Omega^E_{g,D,\delta}$ only depends on the genus of $X$.

If we plug $u=v=t$ in the E-polynomial, we get the virtual Poincar\'e polynomial, which we denote by $P$. Thus $P$ is a ring homomorphism $P\colon\cMot(\kk)^{\lambda-ring}\to\Z((t^{-1}))$ sending the class of a smooth projective variety $Y$ to $\sum_i(-1)^i h^i(X,\Q)t^i$. Making this substitution in~\eqref{eq:EPoly}, we obtain
\[\label{eq:PPoly}
    \Omega^P_{g,D,\delta}:=P(\Omega^{mot}_{X,D,\delta})=
    \sum_{\mu\in\cP}(wt^\delta)^{|\mu|}\prod_{\Box\in\mu}
    \frac{(-t^{2l(\Box)}z^{2a(\Box)})^\delta(t^{2l(\Box)+1}-z^{2a(\Box)+1})^{2g}}
    {(t^{2l(\Box)+2}-z^{2a(\Box)})(t^{2l(\Box)}-z^{2a(\Box)+2})}\prod_{x\in D}\tilde H_\mu(w_{x,\bullet};z^2,t^2).
\]
Set
\[
    \HH^E_{g,D,\delta}:=(1-z^2)\Log\Omega^E_{g,D,\delta}\;\text{ and }\;
    \HH^P_{g,D,\delta}:=(1-z^2)\Log\Omega^P_{g,D,\delta},
\]
where the plethystic logarithm is defined as in Section~\ref{sect:ExpLog}. Clearly, $\HH^E_{g,D,\delta}$ and $\HH^P_{g,D,\delta}$ are the images of $\HH^{mot}_{g,D,\delta}$ under $E$ and $P$ respectively.

In Theorems~\ref{th:Main1},~\ref{th:Main2}, and~\ref{th:Main3} we calculated the motivic classes of various stacks of parabolic $\epsilon$-connections. We can calculate the $E$-polynomials of these stacks simply by replacing $\HH^{mot}_{X,D,\delta}$ with $\HH^E_{g,D,\delta}$ and $\bL$ with $uv$. Similarly, we can calculate the virtual Poincar\'e polynomials by replacing $\HH^{mot}_{X,D,\delta}$ with $\HH^P_{g,D,\delta}$ and $\bL$ with $t^2$. For example, let us do this for Theorem~\ref{th:Main2}.

\begin{corollary}[of Theorem~\ref{th:Main2}]\label{cor:Poincare}
Assume that $\zeta\in\FNF(\divisor)$ and $\gamma'\in\Gamma_{\divisor,\zeta}$. Assume that $d'\in\Z$ is such that $\epsilon d'+\gamma'\star\zeta=0$. Let $\sigma$ be a sequence of parabolic weights of type $(X,\divisor)$. Set $\tau:=(d'+\gamma'\star\sigma)/\rk\gamma'$. Then the E-polynomial of $\Conn^{\sigma-ss}_{\gamma',d'}(\epsilon,X,\divisor,\zeta)$ is equal to the coefficient at~$w^{\gamma'}$ in
\[
    (-(uv)^{1/2})^{\chi(\gamma')} \Exp\left(uv\left(\HH^E_{g,D,\delta}\Bigl|_{z=1}\right)_{\substack{\gamma\star\zeta\in\epsilon\Z\\ \gamma\star\sigma-\tau\rk\gamma\in\Z }}\right).
\]
Similarly, the virtual Poincar\'e polynomial of $\Conn^{\sigma-ss}_{\gamma',d'}(\epsilon,X,\divisor,\zeta)$ is equal to the coefficient at $w^{\gamma'}$ in
\[
    (-t)^{\chi(\gamma')}
    \Exp\left(t^2\left(\HH^P_{g,D,\delta}\Bigl|_{z=1}\right)_{\substack{\gamma\star\zeta\in\epsilon\Z\\ \gamma\star\sigma-\tau\rk\gamma\in\Z }}\right).
\]
\end{corollary}

\subsection{An example}\label{sect:Pantev} In~\cite{DiaconescuDonagiPantev} the authors consider the case when $\divisor'=\divisor=np$, $p\in X(\kk)$,
\[
    \gamma'=(r,(\underbrace{1,1,\ldots,1}_{r\text{ ones}},0,\ldots))=(r,1^r),
\]
$d'\in\Z$, and $\sigma=(\sigma_i)\in\R^{\Z_{>0}}$ is generic in the sense that for any non-empty proper subset $S\subsetneq\{1,\ldots,r\}$ and any $d\in\Z$ we have
\[
    \frac{d+\sum_{i\in S}\sigma_i}{|S|}\ne\frac{d'+\sum_{i=1}^r\sigma_i}r.
\]
Let us take $\epsilon\in\kk$ and $\zeta\in\FNF(\divisor)=(\Omega_X(np)/\Omega_X)^{\Z_{>0}}$. We will calculate the virtual Poincar\'e polynomial of $\Conn_{\gamma',d'}^{\sigma-ss}(\epsilon,\divisor,\zeta)$ in this case. We need to assume that $\epsilon d'+\sum_j\res\zeta_{p,j}=0$, else the moduli stack is trivial. We have $\delta=n-1$, $D=\{p\}$ so
\begin{multline*}
    \HH^P_{g,D,\delta}=\HH^P_{g,\{p\},n-1}=(1-z^2)\Log\Omega^P_{g,\{p\},n-1}\\=(1-z^2)\Log\left(
    \sum_{\mu\in\cP}(wt^{n-1})^{|\mu|}\prod_{\Box\in\mu}
    \frac{(-t^{2l(\Box)}z^{2a(\Box)})^{n-1}(t^{2l(\Box)+1}-z^{2a(\Box)+1})^{2g}}
    {(t^{2l(\Box)+2}-z^{2a(\Box)})(t^{2l(\Box)}-z^{2a(\Box)+2})}\tilde H_\mu(w_{\bullet};z^2,t^2)\right).
\end{multline*}
Next, for a partition $\mu$ set as in~\cite[(1.5)]{DiaconescuDonagiPantev}
  \[
    d(\mu,n,g)=2|\mu|^2(g-1)+2n\sum_{i<j}\mu_i\mu_j+2.
  \]
Write
  \[
    (t^2-1)\HH^P_{g,\{p\},n-1}=\sum_{\mu\in\cP}(-1)^{(n-1)|\mu|}w^{|\mu|}t^{d(\mu,n,g)+|\mu|(n-1)}\HH_{\mu,n}(z,t)m_{\mu}(w_\bullet),
  \]
where $m_\mu$ are the standard monomial symmetric functions.
\begin{proposition}\label{pr:Pantev}
Assume that $\epsilon d'+\sum_j\res\zeta_{p,j}=0$. Then

(i)
\[
    P(\Conn_{(r,1^r),d'}^{\sigma-ss}(\epsilon,np,\zeta))=
    t^{2d(1^r,n,g)}\frac{\HH_{1^r,n}(1,t)}{t^2-1}=
    \frac{\HH_{1^r,n}(1,t^{-1})}{t^2-1}.
\]

(ii) The stack $\Conn_{(r,1^r),d'}^{\sigma-ss}(\epsilon,np,\zeta)$ has a smooth coarse moduli space whose virtual Poincar\'e polynomial is
\begin{equation}\label{eq:Poincare}
    t^{2d(1^r,n,g)}\HH_{1^r,n}(1,t)=\HH_{1^r,n}(1,t^{-1}).
\end{equation}
\end{proposition}
\begin{proof}
We would like to apply Corollary~\ref{cor:Poincare}. Take $\tau:=\frac{d'+\sum_{i=1}^r\sigma_i}r$. We have
\[
t^2\Bigl(\HH^P_{g,\{p\},n-1}\Bigr|_{z=1}\Bigr)_{\substack{\gamma\star\zeta\in\epsilon\Z\\ \gamma\star\sigma-\tau\rk\gamma\in\Z }}=\left(\frac{(-1)^{(n-1)r}t^{d(1^r,n,g)+r(n-1)+2}\HH_{1^r,n}(1,t)}{t^2-1}\right)w^{(r,1^r)}+\ldots,
\]
where $w^{(r,1^r)}=w^rw_{p,1}\ldots w_{p,r}$. It follows from the genericity assumption on $\sigma$ that among the omitted terms there are no non-zero monomials $A_\gamma w^\gamma$ with $\gamma\le(r,1^r)$. Thus, when we exponentiate this expression, we will have the same coefficient at $w^{(r,1^r)}$. Now, Corollary~\ref{cor:Poincare} shows that
\[
    P(\Conn_{(r,1^r),d}^{\sigma-ss}(\epsilon,np,\zeta))=
    \frac{(-t)^{\chi((r,1^r))}(-1)^{r(n-1)}t^{d(1^r,n,g)+r(n-1)+2}\HH_{1^r,n}(1,t)}{t^2-1}.
\]
Note that $\chi((r,1^r))$ and $r(n-1)$ have the same parity and that $\chi((r,1^r))+r(n-1)+2=d(1^r,n,g)$. Now the first equality of part~(i) follows.

It is known that the coarse moduli space of the stack $\Conn_{\gamma,d}^{\sigma-ss}(\epsilon,np,\zeta)$ exists as a variety. The stack is a~$\mathbb G_m$-gerbe over its coarse moduli space so their motivic classes differ by a factor of $[\mathbb G_m]=\bL-1$. Hence, the virtual Poincar\'e polynomials differ by a factor of $t^2-1$. We see that the virtual Poincar\'e polynomial of the coarse moduli space is given by the LHS of~\eqref{eq:Poincare}. Moreover, the coarse moduli space is a smooth  quasi-projective variety (thanks to the genericity assumption on $\sigma$) of dimension $d(\mu,n,g)$ and $t^{2d(\mu,n,g)}\HH_{1^r,n}(1,t)$ is its compactly supported Poincar\'e polynomial. Hence, the equality in~\eqref{eq:Poincare} follows from Poincar\'e duality. Now the second equality in part~(i) follows as well.
\end{proof}

\begin{remark}\label{rem:Pantev}
  Let us compare the above formula with predictions of~\cite{DiaconescuDonagiPantev} when $\epsilon=0$ (that is, for Higgs bundles). In~\cite[Sect.~1.6]{DiaconescuDonagiPantev} the authors consider the coarse moduli spaces $\cH^s_{\underline\xi}(C,D;\underline\alpha,d,\underline m)$ of Higgs bundles on a curve $C$. In our notation, they are the coarse moduli spaces of the stacks $\Conn^{\underline\alpha-ss}_{\underline m,d}(0,C,D,\underline\xi)$. Via the wild non-abelian Hodge correspondence these moduli spaces have Betti counterparts, denoted $\cS_{Q,M}$ (the irregular type $Q$ and the formal monodromy $M$ are determined by $\underline\alpha$, $\underline m$, and $\underline\xi$). The varieties $\cS_{Q,M}$ are homeomorphic to $\cH^s_{\underline\xi}(C,D;\underline\alpha,d,\underline m)$, so they have the same compactly supported Poincar\'e polynomials.

  Let $Z_{HMV}(z,t,w_\bullet)$ be obtained from $\Omega^P_{g,\{p\},n-1}$ by substituting $w=t^{1-n}$. This coincides with $Z_{HMV}$ defined in~\cite[Sect.~1.2]{DiaconescuDonagiPantev} except that our $t$ should be replaced by $w$ and our $w_\bullet$ should be replaced with $\mathrm x$. Then our $\HH_{\mu,n}(z,t)$ coincides with the one defined in~\cite[(1.7)]{DiaconescuDonagiPantev} except that our $t$ should be replaced with $w$.

  Next,~\cite[(1.6) and (1.8)]{DiaconescuDonagiPantev} predict (following~\cite{HauselMerebWong}) that
  \[
    WP(\cS_{Q,M};u,t)=\HH_{1^r,n}(u^{1/2},-u^{-1/2}t^{-1}),
  \]
  where $WP$ is the weighted Poincar\'e polynomial. Substituting $u=1$ and replacing $t$ with $-t$, we get the usual Poincar\'e polynomial, and the prediction becomes our Proposition~\ref{pr:Pantev}(ii).
\end{remark}

\pdfbookmark[1]{References}{ref}

\bibliographystyle{alphanum}
\bibliography{irregular}
\end{document}